\newif\ifarxived
\arxivedtrue
\ifx\Arxived\Kaba\else\arxivedtrue\fi
\newif\ifpcversion
\ifx\Pcversion\Kaba\else\pcversiontrue\fi
\newif\ifextended
\ifx\Extended\Kaba\else\extendedtrue\fi
\newif\ifprivate
\ifx\Private\Kaba\else\privatetrue\extendedtrue\fi
\newif\ifJapanese
\newif\iftesting
\ifextended\testingtrue\fi
\ifarxived
\extendedtrue
\testingfalse
\privatefalse
\fi
\documentclass[uplatex, 12pt]{article}
\ifextended
\fi
\ifarxived
\usepackage[whole]{bxcjkjatype} 
\usepackage{setspace}
\fi
\ifarxived
\usepackage[bookmarks=true,backref=page,bookmarksnumbered=true, 
  bookmarkstype=toc,
  pdftitle={Extendible cardinals, and Laver-generic large cardinal axioms for extendibility},%
  pdfsubject={},%
  pdfkeywords={Laver-generic large cardinal, Maximality Principles, Recurrence Axioms, 
    generic absoluteness}]{hyperref}
\else
\usepackage[dvipdfmx,backref=page,bookmarks=true,bookmarksnumbered=true,
  bookmarkstype=toc,
]{hyperref}
\fi
\ifarxived
\else
\ifextended
\usepackage[dvipdfmx]{pxjahyper}
\fi
\fi
\ifarxived
\hypersetup{
  colorlinks=true,
  linkcolor=blue,
  filecolor=blue,      
  citecolor=cyan,
  urlcolor=cyan,
  pdftitle={Extendible cardinals, and Laver-generic large cardinal axioms for extendibility},%
  pdfauthor={Sakaé Fuchino},%
  pdfsubject={},%
  pdfkeywords={Laver-generic large cardinal, Maximality Principles, Recurrence Axioms, 
    generic absoluteness}
}
\else
\ifextended
\hypersetup{
  colorlinks=true,
  linkcolor=blue,
  filecolor=blue,      
  citecolor=cyan,
  urlcolor=cyan,
  pdftitle={Extendible cardinals, and Laver-generic large cardinal axioms for extendibility},%
  pdfauthor={Sakaé Fuchino (渕野 昌)},%
  pdfsubject={},%
  pdfkeywords={Laver-generic large cardinal, Maximality Principles, Recurrence Axioms, 
    generic absoluteness}
}
\else
\hypersetup{
  colorlinks=true,
  linkcolor=black,
  filecolor=black,      
  citecolor=black,
  urlcolor=black,
  pdftitle={Extendible cardinals, and Laver-generic large cardinal axioms for extendibility},%
  pdfauthor={Sakaé Fuchino},%
  pdfsubject={},%
  pdfkeywords={Laver-generic large cardinal, Maximality Principles, Recurrence Axioms, 
    generic absoluteness}
}
\fi
\fi

\ifarxived
\usepackage{graphicx, color}
\else
\usepackage[dvipdfmx]{graphicx, color}
\fi
\usepackage{amsmath, amssymb}
\usepackage{bbm}
\usepackage{dsfont}
\usepackage{bbold}
\usepackage{accents} 
\usepackage{marginnote}
\usepackage[many]{tcolorbox}
\usepackage{mathtools} 
\usepackage[normalem]{ulem}
\usepackage{setspace}
\definecolor{darkelectricblue}{rgb}{0.33, 0.39, 0.52}
\definecolor{darkgreen}{rgb}{0.31, 0.47, 0.26}
\ifextended
\newcommand{\extendedcolor}{\ifarxived\else\color{darkelectricblue}\fi}

\newcommand{\privatecolor}{\color{darkgreen}}
\newcommand{\darkred}{\color[rgb]{0.8,0.1,0.1}}

\newcommand{\It}{\it\darkred{}}
\else

\newcommand{\darkred}{}

\newcommand{\It}{\it{}}
\fi
\usepackage{theorem}
\makeatletter
\def\amsbb{\use@mathgroup \M@U \symAMSb}
\makeatother
\newcommand{\bbd}[1]{{\amsbb{#1}}}
\theorembodyfont{\ifJapanese\rm\else\it\fi}
\newcount\minute	
\newcount\hour		
\newcount\hourMins  
\ifJapanese
\def\today%
{
  \the\year 年\,\zeroPadTwo{\the\month}月\,\zeroPadTwo{\the\day}日%
}
\fi
\def\now%
{
  \minute=\time    
  \hour=\time \divide \hour by 60 
  \hourMins=\hour \multiply\hourMins by 60
  \advance\minute by -\hourMins 
  \zeroPadTwo{\the\hour}:\zeroPadTwo{\the\minute}%
}
\def\zeroPadTwo#1%
{
  \ifnum #1<10 0\fi    
  #1
}
\title{\bf \scalebox{1}[1.2]{Extendible cardinals, and Laver-generic}\\
  \scalebox{1}[1.2]{large cardinal axioms for extendibility}}
\author{\protect\scalebox{1}[1.2]{\quad\bf Sakaé Fuchino$^{\ast}$}\quad\ifarxived
  \protect\scalebox{1}[1.2]{\bf(渕野\ 昌)}\medskip\\
  \else\ifextended\protect\scalebox{1}[1.2]{\bf(渕野\ 昌)}\medskip\\ 
  \else
  {\protect\scalebox{1}[1.2]{\bf(渕野\ 昌)}}
  \fi
  \fi
  }
\date{}

\ifarxived
\else
\ifextended
\fi
\fi
\renewcommand{\baselinestretch}{1.2}
\renewcommand{\thefootnote}{(\arabic{footnote})\,}
\ifpcversion
\newcommand{\pdfpage}[1]{\#page=#1}
\else
\newcommand{\pdfpage}[1]{\#page#1}
\fi
\iftesting
\fi
\iftesting
\newcommand{\Label}[1]{\label{#1}\marginnote{{\color{cyan}%
      \renewcommand{\baselinestretch}{0.4}\tiny 
		  \rlap{#1}}}}
\newcommand{\LabelR}[1]{\label{#1}\marginnote{\color{cyan}{%
      \renewcommand{\baselinestretch}{0.4}\tiny 
		  \hfill\rlap{#1}}}}

\else
\newcommand{\Label}[1]{\label{#1}}
\newcommand{\LabelR}[1]{\label{#1}}

\fi

\def\memo#1{\ifprivate\marginnote{{\privatecolor\normalsize%
      \renewcommand{\baselinestretch}{0.4}\tiny\mbox{}\vspace{-0.2ex}
      \par\relax%
			#1\par\mbox{}}}\else\fi}%
\def\memox#1{}
\def\imemo#1{\ifprivate{\privatecolor\footnotesize\tt
			#1}\else\fi}%
\def\imemox#1{}

\newcounter{frml}[section]
\newcounter{frmla}[section]
\def\thefrml{{\arabic{section}.\arabic{frml}}}
\def\thefrmla{{$\aleph$\arabic{section}.\arabic{frmla}}}
\def\frmlabel#1{\refstepcounter{frml}{\def\baka{#1}\ifx\baka\empty\else\label{#1}\fi}%
{\rm({\thefrml})\hfill\hfill\hfill}}
\def\ifrmlabel#1{\refstepcounter{frml}{\def\baka{#1}\ifx\baka\empty\else\label{#1}\fi}%
{\iftesting\darkred\fi\rm({\thefrml})\,:\hspace{0.6em}}}
\def\frmlabela#1{\refstepcounter{frmla}{\def\baka{#1}\ifx\baka\empty\else\label{#1}\fi}%
{\rm({\thefrmla})\hfill\hfill\hfill}}
\def\ifrmlabela#1{\refstepcounter{frmla}{\def\baka{#1}\ifx\baka\empty\else\label{#1}\fi}%
{\iftesting\darkred\fi\rm({\thefrmla})\,:\hspace{0.6em}}}
\def\xitem[#1]{\item[\frmlabel{#1}]\mbox{}%
	\iftesting\marginnote{{\renewcommand{%
				\baselinestretch}{0.6}\color{cyan}\tiny#1}}\fi\ignorespaces}
\def\xitemq[#1]{\item[\frmlabel{#1}]\mbox{}%
	\ignorespaces}
\def\xitemd[#1]#2{\item[{\rm(\ref{#1})$#2$}\hfill\hfill\hfill]}
\def\xitemA[#1]{\item[\frmlabela{#1}]\mbox{}%
	\iftesting\marginnote{{\renewcommand{%
				\baselinestretch}{0.6}\tiny#1}}\fi\ignorespaces}
\def\xitemx[#1]{\item[]}
\def\smashx#1{{#1}}
\def\xitemsub[#1]#2{\item[\frmlabel{#1}$_{#2}$]\mbox{}%
	\iftesting\marginnote{{\renewcommand{%
				\baselinestretch}{0.6}\tiny#1}}\fi\ignorespaces}
\def\xxitem[#1][#2]{\item[{\rm(\ref{#1}{\makebox[1.6ex][r]{#2}})}]\mbox{}%
	\iftesting\marginnote{{\renewcommand{%
				\baselinestretch}{0.6}\tiny\{#1\}\{#2\}}}\fi\ignorespaces}
\def\xitemof#1{{\rm({\ref{#1}})}}
\def\Xitem[#1]{\item[{\makebox[7ex][l]{\rm(\ref{#1})}}]\iftesting\marginnote{{\renewcommand{%
				\baselinestretch}{0.6}\tiny#1}}\fi\ignorespaces}

\newenvironment{xitemize}{\begin{list}{}{\parsep=0.5\smallskipamount%
			\itemindent=-0.4ex%
			\itemsep=0.5\smallskipamount\leftmargin=4em\labelwidth=3em\labelsep=0.7em}}%
							 {\end{list}}
\iftesting
\def\ixitem[#1]{\ifrmlabel{#1}\marginnote{{\color{cyan}\renewcommand{%
				\baselinestretch}{0.6}\quad\tiny\rlap{#1}\mbox{}}}\ignorespaces}
\else
\def\ixitem[#1]{\ifrmlabel{#1}\ignorespaces}
\fi
\iftesting
\def\ixitemr[#1]{\ifrmlabel{#1}\marginnote{{\color{cyan}\renewcommand{%
				\baselinestretch}{0.6}\mbox{}\tiny\hfill{}\rlap{\qquad #1}{}\hspace{10ex}}}\ignorespaces}
\else
\def\ixitemr[#1]{\ifrmlabel{#1}\ignorespaces}
\fi
\iftesting
\def\ixitemc[#1]{\ifrmlabel{#1}\marginnote{{\color{cyan}\renewcommand{%
				\baselinestretch}{0.6}\mbox{}\tiny\qquad\qquad\qquad\qquad{}\rlap{#1}{}}}\ignorespaces}
\else
\def\ixitemc[#1]{\ifrmlabel{#1}\ignorespaces}
\fi
\iftesting
\def\ixitemq[#1]{\frmlabel{#1}\marginnote{{\color{cyan}\renewcommand{
				\baselinestretch}{0.6}\qquad\qquad\tiny\rlap{#1}\mbox{}}}\ignorespaces}
\else
\def\ixitemq[#1]{\frmlabel{#1}\ignorespaces}
\fi
\iftesting
\def\ixitemd[#1]#2{\mbox{(\ref{#1})#2\,:\ }\marginnote{{\color{cyan}\renewcommand{%
				\baselinestretch}{0.6}{}\qquad\qquad\tiny\rlap{see: #1}\mbox{}}}\ignorespaces}
\else
\def\ixitemd[#1]#2{(\mbox{\ref{#1})#2\,:\ }\ignorespaces}
\fi
\iftesting
\def\ixitema[#1]{\ifrmlabela{#1}\marginnote{{\color{cyan}\renewcommand{%
				\baselinestretch}{0.6}\quad\tiny\rlap{#1}}}\ignorespaces}
\else
\def\ixitema[#1]{\ifrmlabela{#1}\ignorespaces}
\fi
\iftesting
\def\ixitemac[#1]{\ifrmlabela{#1}\marginnote{{\color{cyan}\renewcommand{%
				\baselinestretch}{0.6}\mbox{}\tiny\qquad\qquad\qquad\qquad{}\rlap{#1}{}}}\ignorespaces}
\else
\def\ixitemac[#1]{\ifrmlabela{#1}\ignorespaces}
\fi
\iftesting
\def\ixitemar[#1]{\ifrmlabela{#1}\marginnote{{\color{cyan}\renewcommand{%
				\baselinestretch}{0.6}\mbox{}\tiny\hfill{}\rlap{\qquad #1}{}\hspace{10ex}}}\ignorespaces}
\else
\def\ixitemar[#1]{\ifrmlabela{#1}\ignorespaces}
\fi

\def\assert#1{\noindent\makebox[4.8ex][r]{\rm(\makebox[2.2ex][c]{#1})}\ \ \ignorespaces}
\def\wassert#1{\assert{#1}}
\def\wassertof#1{\makebox[4.8ex][r]{\rm(\makebox[2.2ex][c]{#1})\ }}%
\def\assertof#1{{(#1)}}%

\def\daimaru#1{\makebox[1em][c]{\mbox{\leavevmode\lower.144ex\hbox{
        \rlap{\hbox to 
          0.76em{\hfil\mbox{}\hfill{}\raisebox{0.054ex}{\scalebox{1.2}{○}}\hfil}}
        \raise0.342ex\hbox to 1em{\hfil{\hspace{0.16em}\footnotesize#1}\hfil}}}}\,}

\newtheorem{Thm}{\ifJapanese{\bf 定理}\else {\bf Theorem}\fi}[section]
\ifextended\tcolorboxenvironment{Thm}{
  colback=blue!4!white,
  boxrule=0pt,
  boxsep=0pt,
  left=8pt,right=4pt,top=2pt,bottom=4pt,
  oversize=4pt,
  sharp corners,
  before skip=\topsep,
  after skip=\topsep,
  breakable
}\fi

\newtheorem{ThmA}{\ifJapanese{\bf 定理\,A\!}\else{\bf Theorem\,A\!}\fi}[section]
\ifextended\tcolorboxenvironment{ThmA}{
  colback=blue!4!white,
  boxrule=0pt,
  boxsep=0pt,
  left=8pt,right=4pt,top=2pt,bottom=4pt,
  oversize=4pt,
  sharp corners,
  before skip=\topsep,
  after skip=\topsep,
  breakable
}\fi

\ifextended\tcolorboxenvironment{Ex}{
  colback=blue!3!white,
  boxrule=0pt,
  boxsep=0pt,
  left=8pt,right=4pt,top=2pt,bottom=4pt,
  oversize=4pt,
  sharp corners,
  before skip=\topsep,
  after skip=\topsep,
  breakable
}\fi
\newtheorem{Prop}[Thm]{\ifJapanese{\bf 命題}\else{\bf Proposition}\fi}
\ifextended\tcolorboxenvironment{Prop}{
  colback=blue!4!white,
  boxrule=0pt,
  boxsep=0pt,
  left=8pt,right=4pt,top=2pt,bottom=4pt,
  oversize=4pt,
  sharp corners,
  before skip=\topsep,
  after skip=\topsep,
  breakable
}\fi

\newtheorem{Lemma}[Thm]{\ifJapanese{\bf 補題}\else{\bf Lemma}\fi}
\ifextended\tcolorboxenvironment{Lemma}{
  colback=blue!4!white,
  boxrule=0pt,
  boxsep=0pt,
  left=8pt,right=4pt,top=2pt,bottom=4pt,
  oversize=4pt,
  sharp corners,
  before skip=\topsep,
  after skip=\topsep,
  breakable
}\fi
\newtheorem{LemmaA}[ThmA]{\ifJapanese{\bf 補題\,A\!}\else{\bf Lemma\,A\!}\fi}
\ifextended\tcolorboxenvironment{LemmaA}{
  colback=blue!4!white,
  boxrule=0pt,
  boxsep=0pt,
  left=8pt,right=4pt,top=2pt,bottom=4pt,
  oversize=4pt,
  sharp corners,
  before skip=\topsep,
  after skip=\topsep,
  breakable
}\fi
\newtheorem{PropA}[ThmA]{\ifJapanese{\bf 命題\,A\!}\else{\bf Proposition\,A\!}\fi}
\ifextended\tcolorboxenvironment{PropA}{
  colback=blue!4!white,
  boxrule=0pt,
  boxsep=0pt,
  left=8pt,right=4pt,top=2pt,bottom=4pt,
  oversize=4pt,
  sharp corners,
  before skip=\topsep,
  after skip=\topsep,
  breakable
}\fi

\ifextended\tcolorboxenvironment{CorA}{
  colback=blue!4!white,
  boxrule=0pt,
  boxsep=0pt,
  left=8pt,right=4pt,top=2pt,bottom=4pt,
  oversize=4pt,
  sharp corners,
  before skip=\topsep,
  after skip=\topsep,
  breakable
}\fi
\newtheorem{Cor}[Thm]{\ifJapanese{\bf 系}\else{\bf Corollary}\fi}
\ifextended\tcolorboxenvironment{Cor}{
  colback=blue!4!white,
  boxrule=0pt,
  boxsep=0pt,
  left=8pt,right=4pt,top=2pt,bottom=4pt,
  oversize=4pt,
  sharp corners,
  before skip=\topsep,
  after skip=\topsep,
  breakable
}\fi

\ifextended\tcolorboxenvironment{Remark}{
  colback=blue!4!white,
  boxrule=0pt,
  boxsep=0pt,
  left=8pt,right=4pt,top=2pt,bottom=4pt,
  oversize=4pt,
  sharp corners,
  before skip=\topsep,
  after skip=\topsep,
  breakable
}\fi

\newtheorem{Claim}{{\bf Claim}}[Thm]
\ifextended\tcolorboxenvironment{Claim}{
  colback=pink!5!white,
  boxrule=0pt,
  boxsep=0pt,
  left=8pt,right=4pt,top=2pt,bottom=4pt,
  oversize=4pt,
  sharp corners,
  before skip=\topsep,
  after skip=\topsep,
  breakable
}\fi

\ifextended\tcolorboxenvironment{ClaimA}{
  colback=pink!5!white,
  boxrule=0pt,
  boxsep=0pt,
  left=8pt,right=4pt,top=2pt,bottom=4pt,
  oversize=4pt,
  sharp corners,
  before skip=\topsep,
  after skip=\topsep,
  breakable
}\fi

\newcommand{\prf}{\noindent\ifJapanese{\bf 証明．\ }\ignorespaces\else{\bf 
		Proof.\ \ }\ignorespaces\fi}

\newcommand{\prfofClaim}{\raisebox{-.4ex}{\Large $\vdash$\ \ }}

\newcommand{\prfof}[1]{\ifJapanese{\bf #1 の証明．\ \ }%
	\ignorespaces\else{\bf Proof of #1:}\ \ \ignorespaces\fi}
\newcommand{\Thmof}[1]{\ifJapanese{定理\,\ref{#1}}\else{Theorem~\ref{#1}}\fi}

\newcommand{\bfThmof}[1]{\ifJapanese{\bf 定理\,\ref{#1}}\else{\bf Theorem~\ref{#1}}\fi}
\newcommand{\Lemmaof}[1]{\ifJapanese{補題\,\ref{#1}}\else{Lemma~\ref{#1}}\fi}

\newcommand{\LemmaAof}[1]{\ifJapanese{補題\,A\,\ref{#1}}\else{Lemma\,A\,\ref{#1}}\fi}

\newcommand{\Propof}[1]{\ifJapanese{命題\,\ref{#1}}\else{Proposition~\ref{#1}}\fi}
\newcommand{\PropAof}[1]{\ifJapanese{命題\,A\,\ref{#1}}\else{Proposition\,A\,\ref{#1}}\fi}

\newcommand{\Corof}[1]{\ifJapanese{系\,\ref{#1}}\else{Corollary~\ref{#1}}\fi}

\newcommand{\Claimof}[1]{{Claim \ref{#1}}}

\newcommand{\sectionof}[1]{\ifJapanese{第\ref{#1}節}\else{Section~\ref{#1}}\fi}

\newcommand{\Thmabove}{{\ifJapanese 定理\else Theorem\fi\ \number\theThm}}

\newcommand{\Lemmaabove}{{\ifJapanese 補題\else Lemma\fi\ \number\theThm}}

\newcommand{\Claimabove}{{Claim \number\theClaim}}

\newcommand{\ubecause}[3]{\underbrace{{}#1{}%
  \ifx\bakakaba#2\bakakaba\rule[-0.72ex]{0pt}{1pt}\else\rule[#2]{0pt}{1pt}\fi}_{\mbox{\footnotesize\clap{#3}}}}
\newcommand{\obecause}[3]{\overbrace{{}#1{}%
  \ifx\bakakaba#2\bakakaba\rule[1.62ex]{0pt}{1pt}\else\rule[#2]{0pt}{1pt}\fi}^{\mbox{\footnotesize\clap{#3}}}}
\newsavebox{\qedbox}\sbox{\qedbox}{
{\unitlength=0.05mm \begin{picture}(40,60)
\put(0,0){\framebox(30,44)[cc]{}}
\put(30,-7){\rule{7\unitlength}{44\unitlength}}
\put(10,-7){\rule{27\unitlength}{7\unitlength}}
\end{picture}}}
\newcommand{\qed}{\mbox{}\hfill\usebox{\qedbox}}
\newcommand{\smallqed}%
{\mbox{}\smallskip\hfill\raisebox{-.4ex}{\Large $\dashv$}}
\newcommand{\qedof}[1]%
{\mbox{} \hspace*{\fill}{\usebox{\qedbox}{\tiny~(#1)}}}
\newcommand{\Qedof}[1]%
{\mbox{} \hspace*{\fill}{\usebox{\qedbox}%
{\tiny~(#1~\number\theThm)}}}
\newcommand{\QedAof}[1]%
{\mbox{} \hspace*{\fill}{\usebox{\qedbox}%
{\tiny~(#1~\number\theThmA)}}}
\newcommand{\qedofThm}{\Qedof{\ifJapanese 定理\else Theorem\fi}}

\newcommand{\qedofCor}{\Qedof{\ifJapanese 系\else Corollary\fi}}
\newcommand{\qedofProp}{\Qedof{\ifJapanese 命題\else Proposition\fi}}
\newcommand{\qedofLemma}{\Qedof{\ifJapanese 補題\else Lemma\fi}}
\newcommand{\qedofLemmaA}{\QedAof{\ifJapanese 補題\,A\!\else Lemma\,A\!\fi}}

\newcommand{\qedskip}{\medskip}

\newcommand{\qedofClaim}%
{\mbox{}\hfill\raisebox{-.4ex}{\Large $\dashv$ }\nolinebreak%
\mbox{\tiny~(Claim~\number\theClaim)}}
\newcommand{\qedofClaimA}%
{\mbox{}\hfill\raisebox{-.4ex}{\Large $\dashv$ }\nolinebreak%
\mbox{\tiny~(Claim~A\,\number\theClaimA)}}
\newcommand{\qedofClaimAof}[1]%
{\mbox{}\hfill\raisebox{-.4ex}{\Large $\dashv$ }\nolinebreak%
\mbox{\tiny~(Claim~A\,\ref{#1})}}
\newcommand{\qedofSubclaim}%
{\mbox{}\hfill\raisebox{-.4ex}{\Large $\dashv$ }\nolinebreak%
\mbox{\tiny~(Subclaim~\number\theSubclaim)}}
\newcommand{\qedofSubsubclaim}%
{\mbox{}\hfill\raisebox{-.4ex}{\Large $\dashv$ }\nolinebreak%
\mbox{\tiny~(Subsubclaim~\number\theSubsubclaim)}}
\newcommand{\cardof}[1]{\mathopen{|\,}#1\mathclose{\,|}}

\newcommand{\rank}{rank\,}

\newcommand{\Card}{{\it Card\/}}

\newcommand{\setof}[2]{\{#1\,:\,#2\}}
\newcommand{\ssetof}[1]{\{#1\}}

\newcommand{\subseteqand}[1]{\mathrel{\mathop{\subseteq}%
		\limits_{\scriptscriptstyle\hbox to 14pt{$\scriptscriptstyle #1$\hss}}}}

\newcommand{\mapping}[3]{#1:#2\rightarrow #3}

\newcommand{\elembed}[3]{#1:#2\stackrel{\preccurlyeq\hspace{0.8ex}}{\rightarrow}#3}
\newcommand{\Elembed}[4]{#1:#2\stackrel{\prec\hspace{0.8ex}}{\rightarrow}_{#4}#3}

\newcommand{\fnsp}[2]{\mbox{}^{{#1}\hspace{-0.02em}}#2}
\newcommand{\imageof}{{}^{\,{\prime}{\prime}}}

\newcommand{\seqof}[2]{\langle#1\,:\,#2\rangle}
\newcommand{\pairof}[1]{\langle#1\rangle}
\newcommand{\psof}[1]{{\mathcal P}\/(#1)}

\newcommand{\forces}[2]{\,\|\hspace{-.35ex}\mbox{\sf--}_{\,#1\,}%
\mbox{\rm``}\,#2\,\mbox{\rm''}}

\newcommand{\modelof}[1]{\models\!\mbox{\rm``\,}#1\mbox{\rm\,''}}

\newcommand{\crit}{\mbox{\it crit\/}}

\newcommand{\bbone}{{\mathord{\mathbb{1}}}}

\newcommand{\circleq}{\mathrel{{\leqslant}%
		\hspace{-0.86ex}{\lower-0.53ex\hbox{$\scriptscriptstyle\circ$}}}}
\newcommand{\symb}[1]{{\mathord{\hspace{0.08em}\underbracket[0.6pt][2pt]{#1}}\hspace{0.08em}}}

\newcommand{\ol}[1]{\overline{#1}}
\newcommand{\ul}[1]{\underline{#1}}
\newcommand{\restr}{\restriction}

\newcommand{\cf}{\mathop{cf\/}}

\newcommand{\Col}{{\rm Col}}

\newcommand{\Fn}{{\rm Fn}}

\newcommand{\trcl}{\mathop{\mbox{\it trcl\/}}}

\newcommand{\natnums}{{\bbd{N}}}

\newcommand{\poP}{\bbd{P}}
\newcommand{\poQ}{\bbd{Q}}

\newcommand{\On}{{\rm On}}

\newcommand{\genG}{\mathbb{G}}

\newcommand{\utpoP}{\utilde{\mathbb{P}}}

\newcommand{\utpoQ}{\utilde{\mathbb{Q}}}

\newcommand{\genH}{\mathbb{H}}

\newcommand{\condp}{\mathbbm{p}}

\newcommand{\LT}{{<}\,}

\newcommand{\GT}{{>}\,}
\newcommand{\GE}{{\geq}\,}

\newcommand{\ctenten}{,\mbox{}\hspace{0.08ex}{.}{.}{.}\hspace{0.1ex}}
\newcommand{\ctentenc}{,{}\linebreak[0]\hspace{0.04ex}{{.}{.}{.}\hspace{0.1ex},\,}\linebreak[0]}

\newcommand{\xmbox}[1]{ $\relax{\rm #1}\relax$ }
\newcommand{\gmA}{\mathfrak{A}}
\newcommand{\gmB}{\mathfrak{B}}

\newcommand{\continuum}{2^{\aleph_0}}

\newcommand{\calC}{{\mathcal C}}
\newcommand{\calD}{{\mathcal D}}

\newcommand{\calH}{{\mathcal H}}

\newcommand{\calL}{{\mathcal L}}

\newcommand{\calO}{{\mathcal O}}
\newcommand{\calP}{{\mathcal P}}

\newcommand{\calS}{{\mathcal S}}

\newcommand{\calV}{{\mathcal V}}
\newcommand{\calW}{{\mathcal W}}

\newcommand{\uta}{\utilde{a}}
\newcommand{\utb}{\utilde{b}}

\newcommand{\utS}{\utilde{S}}

\newcommand{\Lin}{{\calL}_{\in}}

\newcommand{\ZF}{{\sf ZF}}
\newcommand{\ZFC}{{\sf ZFC}}
\newcommand{\CH}{{\sf CH}}

\newcommand{\SCH}{{\sf SCH}}
\newcommand{\SSH}{{\sf SSH}}

\newcommand{\LC}{{\it{LC}}}

\newcommand{\MM}{{\sf MM}}
\newcommand{\MMpp}{\mbox{\sf MM$^{++}$}}
\newcommand{\MApp}{\mbox{\sf MA$^{++}$}}
\newcommand{\PFApp}{\mbox{\sf PFA$^{++}$}}
\newcommand{\MP}{{\sf MP}}

\newcommand{\RcA}{{\sf RcA}}

\newcommand{\RcAp}{{\sf RcA}{$^+$}}

\newcommand{\BfRA}{{\mathbb{R}\mathbb{A}}}
\newcommand{\PFA}{{\sf PFA}}

\newcommand{\BFA}{{\sf BFA}}

\newcommand{\CCA}{{\sf CCA}}
\newcommand{\GA}{{\sf GA}}
\newcommand{\BA}{{\sf BA}}

\newcommand{\RP}{{\sf RP}}

\newcommand{\FRP}{{\sf FRP}}

\newcommand{\GRP}{{\sf GRP}}
\newcommand{\GRPp}{\mbox{\sf GRP$^+$}}
\newcommand{\RC}{{\sf RC}}

\newcommand{\refl}{{\mathfrak{r}\mathfrak{e}\mathfrak{f}\mathfrak{l}\,}}

\newcommand{\st}{such that}
\newcommand{\wrt}{with respect to}
\newcommand{\Wolog}{Without loss of generality}

\newcommand{\tfae}{the following are equivalent}

\newcommand{\po}{poset}
\newcommand{\pos}{posets}
\newcommand{\uniV}{\mathsf{V}}

\newcommand{\uniW}{\mathsf{W}}

\newcommand{\Pkl}[2]{\ifx\bakakaba#1\bakakaba\ifx\bakakaba#2\bakakaba{\mathcal 
    P}_\kappa(\lambda)\else{\mathcal P}_\kappa(#2)\fi\else{\mathcal P}_{#1}(#2)\fi}

\newcommand{\utildeT}[1]{%
  \underaccent{{\sim}}{#1}}
\newcommand{\utildeS}[1]{%
	\hbox to 0pt{\smash{$\mathop{\scriptstyle #1}\limits_{%
				\raisebox{0.6ex}[0pt]{$\scriptscriptstyle\sim$}}$}\hss}%
	\relax\phantom{\mathord{{#1}_{\rule[-0.6ex]{0pt}{1pt}}}}}
\newcommand{\utildeSS}[1]{%
	\hbox to 0pt{$\mathop{\scriptscriptstyle #1}%
		\limits_{\scriptscriptstyle\sim}$\hss}%
		\relax\phantom{\underline{#1}}}
\newcommand{\utilde}[1]{%
	\mathchoice{\utildeT{#1}}{\utildeT{#1}}{\utildeS{#1}}{\utildeSS{#1}}}

{\end{minipage}\end{trivlist}}

\begin{document}
\maketitle
\renewcommand{\thefootnote}{$\ast$\ }
  \footnotetext{Graduate School of System Informatics, Kobe University \\Rokko-dai 1-1, Nada, Kobe 657-8501 Japan
   \\
   \quad\scalebox{0.95}[1]{\tt fuchino@diamond.kobe-u.ac.jp}}
\ifextended
\phantomsection
\addcontentsline{toc}{section}{* Extendible cardinals, and Laver-generic large cardinal axioms for extendibility} 
\ifarxived
\addcontentsline{toc}{section}{**** by S.Fuchino}
\else
\addcontentsline{toc}{section}{**** by S.Fuchino (渕野昌)}
\fi
\fi

\ifextended
\addcontentsline{toc}{section}{Abstract}
\fi
\begin{abstract} {}
  We introduce (super-$C^{(\infty)}$-)Laver-generic large cardinal axioms 
  for extendibility  ((super-$C^{(\infty)}$-)LgLCAs for extendible, for short), and show 
  that most of the previously known consequences  
  of the (super-$C^{(\infty)}$-)LgLCAs for ultrahuge,
  in particular, general 
  forms of Resurrection Principles, Maximality Principles, and Absoluteness Theorems, 
  already follow from (super-$C^{(\infty)}$\mbox{-)}LgLCAs for extendible. 

  The consistency of LgLCAs for extendible (for transfinitely 
  iterable $\Sigma_2$-definable classes of \pos)  
  follows from an extendible cardinal while the  
  consistency of super-$C^{(\infty)}$-LgLCAs for extendible follows from
  a model with a strongly super-$C^{(\infty)}$-extendible cardinal. If $\mu$ is an almost-huge 
  cardinal, there are cofinally many $\kappa<\mu$ \st\
  $V_\mu\modelof{\kappa\xmbox{ is strongly super-}C^{(\infty)}\xmbox{ extendible}}$. 

  Most of the known reflection properties follow already from some of the LgLCAs for supercompact. We give a 
  survey on the related results. 

  We also show the separation between some of the LgLCAs as well as between LgLCAs and 
  their consequences.

  LgLCAs are generic large cardinal axioms in terms of generic elementary embeddings with 
  the critical point $\kappa_\refl=\max\{\aleph_2,2^{\aleph_0}\}$. We show  that Laver 
  generic large cardinal axioms for all \pos\ in terms of generic elementary embeddings 
  with the critical point $2^{\aleph_0}$ is also possible. We abbreviate this type of axiom 
  for the notion of extendibility as the LgLCAA for extendible and examine its 
  consequences. 
\end{abstract}
{}

\phantomsection
\addcontentsline{toc}{section}{Contents}
\newcommand{\myscalebox}[1]{\scalebox{0.88}[1.06]{#1}}
\begin{quotation}
	\footnotesize
	\noindent
	\centerline{
      \normalsize\tt\quad\ Contents\hspace{6em}\mbox{}}\mbox{}\\
  {\mbox{}\hspace{-1.6em}\tt\makebox[3.4ex][l]{\ref{intro}.}%
    \hyperref[intro]{\tt\myscalebox{Introduction}}}\ \ 
  \dotfill\ \ {\pageref{intro}}\\  
  {\mbox{}\hspace{-1.6em}\tt\makebox[3.4ex][l]{\ref{ext-super-ext}.}%
    \hyperref[ext-super-ext]{\tt\myscalebox{Extendible and super-$C^{(n)}$-extendible cardinals}}}\ \ \dotfill\ \ {\pageref{ext-super-ext}}\\  
  {\mbox{}\hspace{-1.6em}\tt\makebox[3.4ex][l]{\ref{ext}.}%
    \hyperref[ext]{\tt\myscalebox{Models with super-$C^{(\infty)}$-extendible cardinals}}}\ \ \dotfill\ \ {\pageref{ext}}\\  
  {\mbox{}\hspace{-1.6em}\tt\makebox[3.4ex][l]{\ref{everyth}.}%
    \hyperref[everyth]{\tt\myscalebox{%
        LgLCAs and 
        super-$C^{(\infty)}$-LgLCAs for extendiblity imply (almost) everything}}}\ \ \dotfill\ \ {\pageref{everyth}}\\  
  {\mbox{}\hspace{-1.6em}\tt\makebox[3.4ex][l]{\ref{Lg-ext}.}%
    \hyperref[Lg-ext]{\tt\myscalebox{Consistency of
        LgLCAs and 
        super-$C^{(\infty)}$-LgLCAs for extendibility}}}\ \ \dotfill\ \ {\pageref{Lg-ext}}\\  
  {\ifextended\extendedcolor
  {\mbox{}\hspace{-1.6em}\tt\makebox[3.4ex][l]{\ref{square}.}%
    \hyperref[square]{\tt\myscalebox{Simultaneous and diagonal reflections, and total failure of square principles}}\\
    \phantom{\mbox{}\hspace{-1.6em}\tt\makebox[3.4ex][l]{\ref{square}.}}%
    \hyperref[square]{\tt\myscalebox{under LgLCAs}}}\ \ \dotfill\ \ {\pageref{square}}\\  
  {\mbox{}\hspace{-1.6em}\tt\makebox[3.4ex][l]{\ref{CCA}.}%
    \hyperref[CCA]{\tt\myscalebox{Separation of some axioms under Continuum Coding Axiom}}}\ \ \dotfill\ \ {\pageref{CCA}}\\ 
  {\mbox{}\hspace{-1.6em}\tt\makebox[3.4ex][l]{\ref{LgLCAA}.}%
    \hyperref[LgLCAA]{\tt\myscalebox{Laver-generic Large Cardinal Axioms for all \pos}}}\ \ \dotfill\ \ {\pageref{LgLCAA}}\\\fi}
  {\mbox{}\hspace{-1.6em}\hyperref[ref]{\tt References}}\ \ \dotfill\ \ 
  {\pageref{ref}}\medskip\\ 
\end{quotation}
\renewcommand{\thefootnote}{}
\footnotetext{{\it Date:} January 23, 2025 \ifextended
  \qquad {\it Last update:} 
  \today\ (\now\ \ifarxived UTC\else JST\fi)\vspace{-1\smallskipamount}
  \else   \qquad {\it Main update:} January 31, 2025 (23:59 JST)
\fi
}
\footnotetext{{\it MSC2020 Mathematical Subject Classification:} 03E45, 03E50, 03E55, 03E57, 03E65
  	\vspace{-1\smallskipamount}}


\ifextended
\ifprivate
\footnotetext{\hspace{-1em}This is a private extended version of the note to appear in RIMS 
  2024 set theory meeting volume of RIMS Kôkyûroku. 
  \par \ifarxived\else All additional 
  details not contained in the submitted version of the paper are either typeset in 
  {\extendedcolor ``dark electric blue''} in case the text is also included in the 
  (non-private) extended 
  version of the paper, or in {\privatecolor ``dark green''} if the comment is thought only 
  for the author's private version. \fi
  \iffalse
  \sout{The numbering of the assertions is kept identical with the submitted version.}
  Since the changes from the submitted version are now quite extensive, I am not trying to keep the 
  numbering of the theorems and assertions identical with the numbering in the submitted 
  version. \else 
  The numbering of the assertions is kept identical with the submitted version.
  \fi

  The most recent version of this private extended edition is downloadable as:\\
  \ifprivate\href{https://fuchino.ddo.jp/papers/RIMS2024-extendible-xx.pdf}{\tt 
  https://fuchino.ddo.jp/papers/RIMS2024-extendible-xx.pdf}
\fi 
}
\else
\footnotetext{\extendedcolor This is an extended version of the note with the 
  title ``Extendible cardinals, and Laver-generic large cardinal axioms for extendibility'' 
  to appear in RIMS Kôkyûroku.
  \par\extendedcolor  All additional 
  details not contained in the submitted version of the paper are either typeset in 
  dark electric blue (the color in which this paragraph is typeset) or put in a separate appendices. 

  The most recent version of this paper is downloadable as:\medskip\\
  \qquad\qquad\href{https://fuchino.ddo.jp/papers/RIMS2024-extendible-\ifpcversion y\else x\fi.pdf}{\tt 
    https://fuchino.ddo.jp/papers/RIMS2024-extendible-\ifpcversion y\else x\fi.pdf} \medskip\\
  This is a version optimized for display on \ifpcversion a PC\else an iPad\fi. If your pdf 
  viewer is running on \ifpcversion an iPad\else a PC\fi, the following file might be 
  better:\medskip\\
  \qquad\qquad\href{https://fuchino.ddo.jp/papers/RIMS2024-extendible-\ifpcversion y\else x\fi.pdf}{\tt 
    https://fuchino.ddo.jp/papers/RIMS2024-extendible-\ifpcversion x\else y\fi.pdf} \medskip\\
}
\fi
\else
\footnotetext{The most recent and extended version of this paper 
  possibly icluding more details and 
  proofs than the present one is downloadable as:\\ \quad \href{https://fuchino.ddo.jp/papers/RIMS2024-extendible-x.pdf}{\tt https://fuchino.ddo.jp/papers/RIMS2024-extendible-x.pdf}}
\fi

\renewcommand{\thefootnote}{\arabic{footnote})\,}
\setcounter{footnote}{0}

\section{Introduction}
\Label{intro}
The present note is a short version of the more extensive \cite{DID} in preparation.

In \sectionof{ext-super-ext}, we begin with reviewing known characterizations of extendible 
cardinals (\Propof{p-ext-super-ext-1}). We then look into the super-$C^{(n)}$ and 
super-$C^{(\infty)}$ large cardinal versions of extendibility, and give their characterizations 
(\Propof{p-ext-super-ext-2}, \Thmof{p-ext-super-ext-3}).

In \sectionof{ext},
\memox{This paragraph should be rewritten! By characterization of super-$C^{(n)}$ extendible 
  cardinals as ...  }
we evaluate the consistency strength of super $C^{(\infty)}$-extendible cardinal:  
It is classical that if $\mu$ is almost-huge, then $V_\mu$ satisfies the Second-order Vopěnka 
Principle (\Lemmaof{p-ext-0}). We show that the Second-order Vopěnka Principle implies that 
there are cofinally many $\kappa<\mu$ \st\
$V_\mu\modelof{\kappa\xmbox{ is strongly super-}C^{(\infty)}\mbox{ extendible}}$ 
(\Propof{p-ext-1}).

In \sectionof{everyth}, we introduce Laver-generic large cardinal versions of these large 
cardinals, and the axioms asserting the existence of 
a/the Laver-generic large cardinal --- the (super-$C^{(\infty)}$-) $\calP$-Laver-generic 
large cardinal axioms for extendibility ((super-$C^{(\infty)}$-\nolinebreak)LgLCAs for extendible, for
short) for various classes $\calP$ of \pos, and show that 
most of the previously known consequences of the (super-$C^{(\infty)}$-)LgLCAs for 
ultrahugeness. In particular, the strong and general 
  forms of Resurrection Principles, Maximality Principles, and Absoluteness Theorems, 
  already follow from the (super-$C^{(\infty)}$-)LgLCAs for extendible. 

The consistency of the LgLCAs for extendible (for transfinitely iterable $\Sigma_2$-definable classes of \pos) 
follows from an extendible cardinal while a $V_\mu$ with 
strongly super-$C^{(\infty)}$ extendible  
cardinal can be generically extended to a model with the super-$C^{(\infty)}$- $\calP$-LgLCAs for extendible for transfinitely 
iterable classes $\calP$ of \pos\ (which may be defined by formulas more complex than
$\Sigma_2$, see \Thmof{p-Lg-ext-1}).  

In contrast, it is known that (super-$C^{(\infty)}$-)LgLCAs for hyperhugeness for 
transfinitely iterable class $\calP$ of \pos,  
axioms apparently stronger than the 
corresponding axioms for ultrahugeness, are equiconsistent with 
the existence of a genuine (super-$C^{(\infty)}$-)hyperhuge (\cite{recurrence}). 

Most of the known reflection properties are consequences of some of the LgLCAs for supercompact. In 
\sectionof{square}, we give a 
survey on this topic.  

In \sectionof{CCA}, we prove separation between some of the LgLCAs as well as between LgLCAs and 
their consequences.  

LgLCAs are generic large cardinal axioms in terms of generic elementary embeddings with the 
critical point
$\kappa_\refl =\max\ssetof{\aleph_2,\continuum}$. 
Laver generic large cardinal 
axioms for all posets in terms of generic elementary embeddings with the critical point
$\continuum$  is also possible. Such axioms are already discussed in \cite{recurrence} and \cite{janos}. 

In \sectionof{LgLCAA}, we abbreviate this type of axiom for the notion of extendibility as ``the LgLCAA for 
extendible'', and examine its soncistency and consequences.

Our notation is standard, and mostly compatible with that of \cite{millennium-book}, 
\cite{higher-inf}, and/or \cite{kunen-2011}, but with the following slight deviations: 
``$
\Elembed{j}{M}{V}{\kappa}$'' expresses the situation that $M$ and $N$ are transitive 
(sets or classes), $j$ is an elementary embedding of $M$ into $N$ and $\kappa$ is the 
critical point of $j$. 
We use letters with under-tilde to denote $\poP$-names for a \po\ $\poP$. 
Underline added to a symbol like $\ul{\alpha}$ emphasizes that the symbol is used to 
denote a variable in a language, mostly the 
language of \ZFC\ which is denoted by $\Lin$. A letter with under-bracket like
$\symb{c}$ emphasizes that the letter denotes a (new) constant symbol added to 
the language. 

In the following,  we always denote with $\calP$, a class of \pos.
We usually assume that the class $\calP$ of \pos\ satisfies the following properties which 
we call {\It iterablility}. A 
class $\calP$ of \pos\ is said to be (two-step) {\It iterable} if 
\begin{xitemize}
\xitem[x-theintro-0-0] 
  $\calP$ is 
  closed \wrt\ forcing equivalence, and $\ssetof{\bbone}\in\calP$; 
\xitem[x-theintro-0-0-a-0] $\calP$ is closed \wrt\ restriction. That is, for $\poP\in\calP$ 
  and $\condp\in\poP$, we always have $\poP\restr\condp\in\calP$; and 
\xitem[x-theintro-0-0-a-1] For any $\poP\in\calP$, and any $\poP$-name $\utpoQ$ of a \po\ with
  $\forces{\poP}{\utpoQ\in\calP}$, we have $\poP\ast\utpoQ\in\calP$.
\end{xitemize}

$\calP$ is {\It transfinitely iterable} if it is iterable and endowed with an iteration of 
arbitrary length for an appropriate notion of support, and the iteration satisfies preservation and 
factor lemmas. 

For a property $P$ of \pos, we shall say ``$\calP$ is $P$'' (e.g. ``$\poP$ is c.c.c.'') to 
indicate that all elements of  
$\calP$ have the property $P$. In contrast, if we say $\calP$ is the class of \pos\ with the 
property $P$, we mean $\calP=\setof{\poP}{\poP\models P}$.

We adopt the notation of \cite{bagaria-Cn} and denote
$C^{(n)}:=\setof{\alpha\in\On}{V_\alpha\prec_{\Sigma_n}\uniV}$ for $n\in\natnums$.
Intuitively we put $C^{(\infty)}:=\setof{\alpha\in\On}{V_\alpha\prec\uniV}$, though this is 
not a definable class in the language of \ZFC, due to undefinability of the truth. For each 
transitive set model $M$, however, $C:=\setof{\alpha\in\On\cap M}{{V_\alpha}^M\prec M}$ is 
a(n existent) set.  Note that 
the first order logic in this context on which the elementary submodel relation $\prec$ 
relies, is not the meta-mathematical one but rather the logic in the set theory whose 
formulas are the corresponding subset of $\omega$ (consisting  of Gödel numbers) in \ZFC. 
Also, $C$ is not a first-order definable subset of $M$, again because of the undefinability of 
the truth in $M$. 

The following (almost trivial) lemma is often used without mention: 
\begin{Lemma}{\rm(see e.g.\ Section 1 in Bagaria \cite{bagaria-Cn})}\Label{p-Lg} For an 
  uncountable cardinal $\alpha$, 
  $\calH(\alpha)=V_\alpha$ if and only if $V_\alpha\prec_{\Sigma_1}\uniV$. 

  If $V_\alpha\prec_{\Sigma_1}\uniV$ then $\alpha$ is an uncountable strong limit cardinal.  

  Thus, we have 
  \begin{xitemize}
  \xitemx[] 
    $C^{(1)}=\setof{\alpha}{\alpha\mbox{ is an uncountable limit cardinal with }V_\alpha=\calH(\alpha)}$.
    \ifextended\else\qed\fi
  \end{xitemize}
\end{Lemma}
{\ifextended\extendedcolor \prf Note that the equality  $\calH(\alpha)=V_\alpha$ only makes 
  sense if $\alpha$ is a cardinal. 

  Suppose first, that $\alpha$ is an uncountable cardinal and $\calH(\alpha)=V_\alpha$ holds. Then
  $V_\alpha=\calH(\alpha)\prec_{\Sigma_1}\uniV$.

  Suppose now that $V_\alpha\prec_{\Sigma_1}\uniV$ holds. Then
  for any $\beta<\alpha$, we have
  $V_\alpha\modelof{\exists\ul{\gamma}\ (\beta<\ul{\gamma}\ \land\ \ul{\gamma}
      \xmbox{ is a cardinal})}$. The witness of $\ul{\gamma}$ is then really a cardinal. 
  This shows that $\alpha$ is a limit cardinal $\GT\omega$, and hence in particular, uncountable. 

  $\calH(\alpha)\subseteq V_\alpha$ 
  holds always for a cardinal $\alpha$: if $a\in\calH(\alpha)$ then
\ixitem[x-Lg-ext-a] $\cardof{\trcl(a)}<\alpha$. Thus, 
  letting $\mapping{r}{\trcl(a)}{\On}$; $c\mapsto\rank(c)$, we have
  $r\imageof\trcl(a)=:\gamma<\alpha$ by \xitemof{x-Lg-ext-a}. Thus $a\subseteq V_\gamma$, 
  and hence $a\in V_{\gamma+1}\subseteq V_\alpha$.

  If $a\in V_\alpha$, we have
  $\uniV\modelof{\exists \ul{x}\exists\ul{\gamma}(a\subseteq \ul{x}
    \land \ul{x}\mbox{ is transitive}\land \cardof{\ul{x}}\leq\ul{\gamma})}$.
  It follows that
  $V_\alpha\modelof{\exists\ul{x}\exists\ul{\gamma}(a\subseteq\ul{x}\land\ul{x}\mbox{ is transitive}
    \land \cardof{\ul{x}}\leq\ul{\gamma})}$.
  Thus $a\in\calH(\alpha)$.

  For the next statement of the lemma, if $V_\alpha\prec_{\Sigma_1}\uniV$ then we already 
  have shown that $\alpha$ is an uncountable limit cardinal. 

  Also for any cardinal $\mu<\alpha$ we have $\psof{\mu}\in V_{\mu+2}\subseteq V_\alpha$ 
  and   $V_\alpha\modelof{\exists\ul{\nu}\ (\ul{\nu}=\cardof{\psof{\mu}})}$ 
  by $\Sigma_1$-elementarity. Thus $2^\mu<\alpha$.    
  \qedofLemma\qedskip
  \fi}

Note that there is a (first-order) sentence $\varphi$ \st\ $V_\alpha\models\varphi$ if and 
only if $V_\alpha=\calH(\alpha)$ for a cardinal $\alpha$.

\imemox{Thanks to Gabe Goldberg and Toshimichi Usuba. }

The following notes are results of an examination of what was suggested by Gabriel Goldberg in a 
discussion we had during 
his visit to Kobe after the RIMS Set Theory Workshop 2024. Toshimichi Usuba pointed out 
some elementary flows in early sketches of the note. I learned some known 
arguments used below in conversation with Hiroshi Sakai. 
I am grateful for their comments and advices. Also I would like to thank Andreas Leitz for giving 
me a permission to present an exposition of his \href{https://fuchino.ddo.jp/papers/RIMS2024-extendible-x.pdf\pdfpage{9}}{proof} of \Thmof{p-ext-super-ext-3} in the 
extended version of the present article.  

Back in the summer of 2015, I enjoyed a pleasant walk around the port 
of Yokohama with Joel Hamkins when we were together on the way to Kyoto starting from Tokyo and 
made a short stop in Yokohama. On the walk, Joel told me about his then recent 
researches and research projects, and one of them was about the Resurrection Axioms.

Now that his Resurrection Axioms are shown to be restricted versions of the LgLCAs
(see \Thmof{T-Lg-RA-0}), I notice that what I learned from him on that walk 
might have influenced me subliminally when I introduced the LgLCAs in the late 2010s. 
In that case, I have to to thank Joel again sincerely, also for the nice conversation we had in Yokohama.

\section{Extendible and super-\texorpdfstring{$C^{(\infty)}$}{C(∞)}-extendible cardinals}
\Label{ext-super-ext}
\memo{Characterization of Extendible and super-$C^{(n)}$-extendible cardinals\\
Scan\_2024-11-23--10.44 extendible - annotated
p.45〜
}In this section, we summarize some well-known and some other less well-known facts about 
extendible cardinals and introduce the notion of super-$C^{(n)}$-extendible cardinals.

It appears that the notion of super-$C^{(n)}$-extendible cardinals is equivalent to some 
other already known strong variants of extendibility, see \Thmof{p-ext-super-ext-3}. At the 
moment, it is 
yet unknown if similar equivalence is also available for super-$C^{(n)}$-ultrahuge 
cardinals, or super-$C^{(n)}$-hyperhuge cardinals.  

It is easy to see that the definition of an extendible cardinal in Kanamori 
\cite{higher-inf} is equivalent to its slight modification: a cardinal $\kappa$ is {\It 
  extendible} if 
\ixitem[x-ext-super-ext-a]
for any $\alpha>\kappa$ there are $\beta\in\On$ 
and $\Elembed{j}{V_\alpha}{V_\beta}{\kappa}$ \st\ \ixitem[x-ext-super-ext-0] $j(\kappa)>\alpha$. 

An extendible cardinal is supercompact (see e.g.\ Proposition 23.6 in 
\cite{higher-inf}). The following is easy to prove:  

\begin{Lemma}\Label{p-ext-super-ext-0}
  If $\kappa$ is extendible then there are class many 
  measurable cardinals.\ifextended\else\qed\fi
\end{Lemma}
{\ifextended\extendedcolor\prf If $\kappa$ is extendible then it is 
supercompact. Hence, in particular $\kappa$ 
is measurable. If $\Elembed{j_0}{V_\gamma}{V_\delta}{\kappa}$ with $j_0(\kappa)>\gamma$ 
then $V_\delta\modelof{\xmbox{\,there is a normal ultrafilter over }j_0(\kappa)}$ by 
elementarity. Since the normal ultrafilter over $j_0(\kappa)$ in $V_\delta$ is really a 
normal ultrafilter, $j_0(\kappa)$ is measurable.\qedofLemma\qedskip
\fi}

Since existence of a supercompact cardinal does not imply existence of any large cardinal  
above it (see Exercise 22.8 in \cite{higher-inf}), \Lemmaabove\ explains the transcendence 
of extendible cardinals above supercompact.

In Jech \cite{millennium-book}, extendibility is defined by \xitemof{x-ext-super-ext-a} without 
\xitemof{x-ext-super-ext-0}. We say in the following that $\kappa$ is {\It Jech-extendible} 
if it satisfies \xitemof{x-ext-super-ext-a} but not necessarily \xitemof{x-ext-super-ext-0}.  
The two definitions of extendibility are equivalent. In \Propof{p-ext-super-ext-1} below, 
we show the equivalence of these two together with some other characterizations of 
extendibility. 

The key fact to \Propof{p-ext-super-ext-1} is that 
the elementary embedding in \xitemof{x-ext-super-ext-a} can be often lifted to an 
elementary embedding of the whole universe $\uniV$. 

We call a mapping $\mapping{f}{M}{N}$ {\It cofinal} (in $N$) if, for all
$b\in N$, there is $a\in M$ \st\ $b\in f(a)$. 

{\ifprivate\privatecolor See also: {\tiny\href{https://math.stackexchange.com/questions/4402840/why-are-these-two-definitions-of-extendible-cardinal-equivalent}{https://math.stackexchange.com/questions/4402840/why-are-these-two-definitions-of-extendible-cardinal-equivalent}}\fi}

\begin{Lemma}
  \Label{p-ext-super-ext-0-0}{\rm(A special case of \href{https://fuchino.ddo.jp/papers/definability-of-glc-x.pdf\pdfpage{10}}{Lemma 6} in Fuchino and Sakai \cite{fuchino-sakai})}
  Suppose that $\theta$ is a cardinal and $\Elembed{j_0}{\calH(\theta)}{N}{}$ for a 
  transitive set $N$. 
  Let  $N_0:=\bigcup j_0\imageof{\calH(\theta)}$. Then 
  $\Elembed{j_0}{\calH(\theta)}{N_0}{}$ 
  and $j_0$ is cofinal in $N_0$.\qed 
\end{Lemma}
  
\begin{Lemma}
  \Label{p-ext-super-ext-0-1}{\rm(A special case of \href{https://fuchino.ddo.jp/papers/definability-of-glc-x.pdf\pdfpage{11}}{Lemma 7} in \cite{fuchino-sakai})}
  For any regular cardinal $\theta$ and any cofinal $\Elembed{j_0}{\calH(\theta)}{N}{}$, 
  there are $j$, $M\subseteq\uniV$ \st\ $\Elembed{j}{\uniV}{M}{}$, $N\subseteq M$ and
  $j_0\subseteq j$. \qed
\end{Lemma}

\begin{Prop}
  \Label{p-ext-super-ext-1} For a cardinal $\kappa$ \tfae:\smallskip

  \wassert{a} $\kappa$ is extendible. \smallskip

  \wassert{b} $\kappa$ is Jech-extendible.\smallskip

  \wassert{a$'$} For all $\lambda>\kappa$, there are $j$, $M\subseteq\uniV$ \st\ 
  $\Elembed{j}{\uniV}{M}{\kappa}$, $j(\kappa)>\lambda$ and $V_{j(\lambda)}\in M$.\smallskip

  \wassert{b$'$} For all $\lambda>\kappa$, there are $j$, $M\subseteq\uniV$ \st\ 
  $\Elembed{j}{\uniV}{M}{\kappa}$, 
  and $V_{j(\lambda)}\in M$.

\end{Prop}
\prf \assertof{a} $\Rightarrow$ \assertof{b}: is clear by definition. 

\assertof{b} $\Rightarrow$ \assertof{a}: This can be proved by an argument similar to that 
of the proof of \assertof{b} $\Rightarrow$ \assertof{a} of \Propof{p-ext-super-ext-2} below. 

\assertof{a} $\Rightarrow$ \assertof{a$'$}: 
follows from Lemmas \ref{p-ext-super-ext-0-0} and \ref{p-ext-super-ext-0-1}.\memo{Scan\_2024-11-23--10.44 
  extendible - annotated p58 }
{\ifextended\extendedcolor\smallskip

  Assume that $\kappa$ is extendible, and suppose $\lambda>\kappa$. We want to show 
  that there is $j$ as in \assertof{b} for this $\lambda$.

  Let $\lambda'>\lambda$ be a regular cardinal \st\ $V_\lambda\in\calH(\lambda')$, and let
  $\lambda''>\lambda$ be \st\ $\calH(\lambda'')=V_{\lambda''}$. 

  By assumption there is $\Elembed{j''_0}{V_{\lambda''}}{V_{\mu''}}{\kappa}$ for 
  some $\mu''$ with $j(\kappa)>\lambda''$ ($>\lambda$).  Letting $j'_0:=j''_0\restr\calH(\lambda')$, we have
  $\Elembed{j'_0}{\calH(\lambda')}{\calH(j(\lambda'))}{\kappa}$.

  Let $N_0:=\bigcup j\imageof{\calH(\lambda')}$. Then we have
  $\Elembed{j'_0}{\calH(\lambda')}{N_0}{\kappa}$, and $j'_0$ is cofinal in $N_0$ by 
  \Lemmaof{p-ext-super-ext-0-0}. By \Lemmaof{p-ext-super-ext-0-1}, $j'_0$ has a lifting
  $j\supseteq j'_0$ with $\Elembed{j}{\uniV}{M}{\kappa}$ for some transitive
  $M\subseteq\uniV$. Since $j(\lambda)=j''_0(\lambda)$, We have
  $j(V_\lambda)=V_{j(\lambda)}\in j(\ssetof{V_\lambda})\subseteq N_0\subseteq M$. Thus this 
  $j$ is as desired. 
\fi}\smallskip

\assertof{a$'$} $\Rightarrow$ \assertof{b$'$}: is trivial. \smallskip

\assertof{b$'$} $\Rightarrow$ \assertof{b}: is obtained by restricting elementary 
embeddings on $\uniV$ to $V_\lambda$'s.  
\qedofProp
\qedskip

The notion of super-$C^{(n)}$-large cardinal was introduced in Fuchino and Usuba 
\cite{recurrence}. \Propof{p-ext-super-ext-1} in mind, we define 
the super-$C^{(n)}$-extendibility as follows: For a natural number $n$, we call a cardinal
$\kappa$ {\It super-$C^{(n)}$-extendible} if 
for any $\lambda_0>\kappa$ there are $\lambda\geq\lambda_0$ with 
$V_\lambda\prec_{\Sigma_n}\uniV$, and $j$, $M\subseteq\uniV$ \st\
$\Elembed{j}{\uniV}{M}{\kappa}$, $j(\kappa)>\lambda$, $V_{j(\lambda)}\in M$ and
$V_{j(\lambda)}\prec_{\Sigma_n}\uniV$.

We call a cardinal $\kappa$ {\It super-$C^{(\infty)}$-extendible} if $\kappa$ is
super-$C^{(n)}$-extendible for all $n\in\omega$. In general, we cannot formulate the assertion
``$\kappa$ is super-$C^{(\infty)}$-extendible'' in the language of \ZF\ since we would need 
an infinitary logic to do this unless we are allowed to introduce a new constant symbol to 
refer the cardinal across the infinitely many formulas expressing the $C^{(n)}$-exitendibility 
of the cardinal for each $n\in\natnums$. However, there are 
certain situations in which we can say that a cardinal is super-$C^{(\infty)}$-extendible. One  
of them is when we are talking about a cardinal in a set model. In this case, being 
``super-$C^{(\infty)}$-extendible'' in the model is an $\calL_{\omega_1,\omega}$ sentence which is satisfied by 
the cardinal in the model. Another situation is when we are talking about a cardinal in an inner model 
and the cardinal is definable in $\uniV$ (e.g.\ as $\continuum$ in the outer 
model $\uniV$). Note that in the latter case, we can formulate the 
super-$C^{(\infty)}$-extendibility of the cardinal in 
infinitely many formulas, and hence $n$ in this case ranges only over metamathematical natural numbers. 

Similarly to \Propof{p-ext-super-ext-1}, we have the following equivalence:
\begin{Prop}\Label{p-ext-super-ext-2}
  For a cardinal $\kappa$ and $n\geq1$, the following are equivalent:\smallskip

  \wassert{a} For any $\lambda_0>\kappa$ there are $\lambda>\lambda_0$ with
  $V_\lambda\prec_{\Sigma_n}\uniV$, $j_0$, and $\mu$ \st\
  $\Elembed{j_0}{V_\lambda}{V_\mu}{\kappa}$, $j(\kappa)>\lambda$, and
  $V_\mu\prec_{\Sigma_n}\uniV$. \smallskip 

  \wassert{b} For any $\lambda_0>\kappa$ there are $\lambda>\lambda_0$ with
  $V_\lambda\prec_{\Sigma_n}\uniV$, $j_0$, and $\mu$ \st\
  $\Elembed{j_0}{V_\lambda}{V_\mu}{\kappa}$, and
  $V_\mu\prec_{\Sigma_n}\uniV$ (without the condition ``$j(\kappa)>\lambda$''). \smallskip 

  \wassert{a$'$} $\kappa$ is super-$C^{(n)}$-extendible. 

  \wassert{b$'$} for any $\lambda_0>\kappa$ there are $\lambda\geq\lambda_0$ with 
  $V_\lambda\prec_{\Sigma_n}\uniV$, and $j$, $M\subseteq\uniV$ \st\
  $\Elembed{j}{\uniV}{M}{\kappa}$, $V_{j(\lambda)}\in M$, and
  $V_{j(\lambda)}\prec_{\Sigma_n}\uniV$ (without the condition ``$j(\kappa)>\lambda$''). \qed
\end{Prop}
\prf The proof is similar to that of \Lemmaof{p-ext-super-ext-1}.
We only show \assertof{b} $\Rightarrow$ \assertof{a}. The following proof is a modification 
of the \href{https://math.stackexchange.com/questions/4402840/why-are-these-two-definitions-of-extendible-cardinal-equivalent}{proof} of \Lemmaof{p-ext-super-ext-1}, \assertof{b} $\Rightarrow$ \assertof{a} given 
by Farmer S in \cite{farmer}. 

Assume, toward a contradiction, that $\kappa$ satisfies \assertof{b} but not \assertof{a}. 
Then there is a $\gamma$ \st
\begin{xitemize}
\xitem[x-ext-super-ext-1]  for all sufficiently large
  $\lambda>\kappa$, 
if \ixitemr[x-ext-super-ext-2] $V_\lambda\prec_{\Sigma_n}\uniV$, and $\mu$, $j$ are 
  \st\ \ixitem[x-ext-super-ext-3] $\Elembed{j}{V_\lambda}{V_\mu}{\kappa}$ and 
\ixitemr[x-ext-super-ext-4] $V_\mu\prec_{\Sigma_n}\uniV$, \\then $j(\kappa)<\gamma$. 
\end{xitemize}
In the following, let $\gamma$ be the least such $\gamma$.

\begin{Claim}\Label{cl-ext-super-ext-0}
  $\gamma$ is a limit ordinal. For all sufficiently large $\lambda$ with \xitemof{x-ext-super-ext-2} and for 
  all $\xi<\gamma$, there are $\mu$, $j$ with \xitemof{x-ext-super-ext-3}, 
  \xitemof{x-ext-super-ext-4} \st\ $j(\kappa)>\xi$. 
\end{Claim}
\prfofClaim Suppose $\gamma$ is not a limit ordinal, say $\gamma=\xi+1$. Then there are 
cofinally many $\lambda\in\On$ 
\st\ $V_\lambda\prec_{\Sigma_n}\uniV$  (actually $\lambda\in\Card$, see \Lemmaof{p-Lg}), 
and there are $j$ and $\mu$ with  
\xitemof{x-ext-super-ext-3}, \xitemof{x-ext-super-ext-4} and $j(\kappa)=\xi$. 
By restricting of $j$'s as right above, it follows that, for all $\lambda>\xi$ with
$V_\lambda\prec_{\Sigma_n}\uniV$, there are $j$ and $\mu$ as above. 

Let 
$\lambda^*$ be a sufficiently large such $\lambda$ where ``sufficiently large'' is meant in 
terms of \xitemof{x-ext-super-ext-1}. Let $j^*$  and $\mu^*$ be \st\
$\Elembed{j^*}{V_{\lambda^*}}{V_{\mu^*}}{\kappa}$, $j^*(\kappa)=\xi$, and $V_{\mu^*}\prec_{\Sigma_n}\uniV$.

Since $\lambda^*\leq\mu^*$, there is also $\Elembed{k}{V_{\mu^*}}{V_{\nu^*}}{\kappa}$ \st\ 
$V_{\nu^*}\prec_{\Sigma_n}\uniV$ and $k(\kappa)=\xi$. But then we have
$\Elembed{k\circ j^*}{V_{\lambda^*}}{V_{\nu^*}}{\kappa}$ and
$k\circ j^*(\kappa)=k(\xi)>k(\kappa)=\xi$. This is a contradiction to 
\xitemof{x-ext-super-ext-1}.

The second assertion of the claim follows from this and the minimality of $\gamma$.\\
\qedofClaim

\begin{Claim}\Label{cl-ext-super-ext-1}
  For all sufficiently large $\mu>\kappa$ with $V_\mu\prec_{\Sigma_n}\uniV$, and $k$, $\nu$ with 
$V_\nu\prec_{\Sigma_n}\uniV$ and $\Elembed{k}{V_\mu}{V_\nu}{\kappa}$, we have
  $k\imageof{\gamma}\subseteq\gamma$. 
\end{Claim}
\prfofClaim Suppose otherwise. Then we find $\xi<\gamma$ \st, for cofinally 
many $\mu>\kappa$ with $V_\mu\prec_{\Sigma_n}\uniV$, there are $\nu$, $k$ \st\
$\Elembed{k}{V_\mu}{V_\nu}{\kappa}$, $V_\nu\prec_{\Sigma_n}\uniV$ and $k(\xi)\geq\gamma$. 
By considering restrictions of $k$'s as above, we conclude that for all $\mu>\xi$ with
$V_\mu\prec_{\Sigma_n}\uniV$, there are $\nu$ and $k$ as above.

Let $\lambda>\xi$ and $j$ (together with rechosen $\mu$ and $k$ for this $\lambda$) be \st\
$V_\lambda\prec_{\Sigma_n}\uniV$, $\Elembed{j}{V_\lambda}{ V_\mu}{\kappa}$ and $j(\kappa)>\xi$ 
(possible by the second half of \Claimof{cl-ext-super-ext-0}).
Then we have $\Elembed{k\circ j}{V_\lambda}{V_\nu}{\kappa}$ and
$k\circ j(\kappa)>k(\xi)\geq\gamma$. 

Since $\lambda$, $\mu$, $\nu$, $j$, $k$ can be chosen \st\ $\lambda$ is sufficiently large 
(in terms of \xitemof{x-ext-super-ext-1}), this is a contradiction.
\qedofClaim\qedskip

Now, let $\lambda>\kappa$ be sufficiently large with $\lambda\geq\gamma+2$, 
$V_\lambda\prec_{\Sigma_n}\uniV$, and $\Elembed{j}{V_\lambda}{V_\mu}{\kappa}$ with
$V_\mu\prec_{\Sigma_n}\uniV$. By \Claimabove, we have $j\imageof{\gamma}\subseteq\gamma$.

{\bf Case 1.} $\cf(\gamma)=\omega$. Then $j(\gamma)=\gamma$ and hence
$\Elembed{j\restr V_{\gamma+2}}{V_{\gamma+2}}{V_{\gamma+2}}{\kappa}$. This is a 
contradiction to Kunen's proof (see e.g.\ Kanamori \cite{higher-inf}, Corollary 23.14). 

{\bf Case 2.} $\cf(\gamma)>\omega$. then, letting $\kappa_0:=\kappa$,
$\kappa_{n+1}:=j(\kappa_n)$ for $n\in\omega$ 
and $\kappa_\omega:=\sup_{n\in\omega}\kappa_n$, we have $\kappa_\omega<\gamma$, and
$\Elembed{j\restr V_{\kappa_{\omega+2}}}{V_{\kappa_\omega+2}}{V_{\kappa_\omega+2}}{\kappa}$. 
This is again a contradiction to Kunen's proof.
\qedofProp\qedskip


Super-$C^{(n)}$-extendibility is actually equivalent to $C^{(n)}$-extendibility of Bagaria 
\cite{bagaria-Cn}. Konstantinos 
Tsaprounis proved the equivalence for a variant of super-$C^{(n)}$-extendibility which he 
called $C^{(n)+}$-extendibility in \cite{tsaprounis3}.

A cardinal $\kappa$ is {\It$C^{(n)}$-extendible} if, for any $\alpha>\kappa$, there is $\beta$
and $\Elembed{j}{V_\alpha}{V_\beta}{\kappa}$ \st\ $j(\kappa)>\alpha$ and
$V_{j(\kappa)}\prec_{\Sigma_n}\uniV$.

{\ifextended\extendedcolor
A cardinal $\kappa$ is {\It$C^{(n)+}$-extendible} if for any $\lambda_0>\kappa$, there are
$\lambda\geq\lambda_0$ with $V_\lambda\prec_{\Sigma_n}\uniV$ and $j$, $M\subseteq V$ \st\
$\Elembed{j}{\uniV}{M}{\kappa}$ $j(\kappa)>\lambda$, $V_{j(\lambda)}\in M$, 
and $V_{j(\kappa)}\prec_{\Sigma_n}V_{j(\lambda)}\prec_{\Sigma_n}\uniV$. 
\fi}

The following notion is introduced by Benjamin Goodman \cite{goodman}.

A cardinal $\kappa$ is {\It supercompact for $C^{(n)}$} if, for any $\lambda>\kappa$ there is
$\Elembed{j}{\uniV}{M}{\kappa}$ \st\ $\fnsp{\lambda}{M}\subseteq M$ and
$C^{(n)}\cap\lambda=(C^{(n)})^M\cap\lambda$.

Andreas Lietz recently found a short proof of the following \Thmof{p-ext-super-ext-3}. 
Goodman possibly proved the equivalence of \assertof{a} and \assertof{c} in the theorem, 
but mentioned only the case of $n=1$ in his \cite{goodman}. \ifextended\else His proof is 
given in the extended version of the present article.\fi

\begin{Thm}{\rm (Andreas Lietz)}\Label{p-ext-super-ext-3} For a cardinal $\kappa$ and for all $n\geq 1$ \tfae:
  \wassertof{a} $\kappa$ is $C^{(n)}$-extendible.\smallskip

  \wassert{b} $\kappa$ is super-$C^{(n)}$-extendible.\smallskip

  \wassert{b'} $\kappa$ is $C^{(n)+}$-extendible.\smallskip

  \wassert{c} $\kappa$ is supercompact for $C^{(n+1)}$.\ifextended\else\qed\fi
\end{Thm}
{\ifextended\extendedcolor
  A detailed exposition of Lietz's proof of the theorem will be given after the following preparations. 

  The next lemma follows immediately from the definition of 
  super-$C^{(n)}$-extendibility.

  \begin{LemmaA}\Label{p-ext-super-ext-4}\extendedcolor
    Suppose $m\leq n$. If $\kappa$ is super-$C^{(n)}$-extendible then it is 
    super-$C^{(m)}$-extendible. \qed 
  \end{LemmaA}
  
  \begin{PropA}\Label{p-ext-super-ext-5}\extendedcolor

    \wassert{1} If $\kappa$ is $C^{(n)}$-extendible then $\kappa\in C^{(m)}$, where $m:=\max\ssetof{3,n+2}$.
    \smallskip

    \wassert{2}
    If $\kappa$ is super-$C^{(n)}$-extendible then $\kappa\in C^{(m)}$, where $m$ is as in {\rm\assertof{1}}.
  \end{PropA}
  \prf \assertof{1}: By Proposition 3.4 in Bagaria \cite{bagaria-Cn} and its proof.  
  \smallskip

  \assertof{2}: We prove the claim of the proposition by induction on $n$.

  For $n=0$, $\kappa$ is super-$C^{(n)}$-extendible if and only if $\kappa$ is extendible. 
  Since $V_\kappa\prec_{\Sigma_3}\uniV$ (Proposition 23.10 in \cite{higher-inf}), the claim 
  of the proposition holds.

  Suppose that the claim of the proposition holds for $n=k$. We show that the claim also 
  holds for $n+1$. 

  Suppose that $\kappa$ is super-$C^{(k+1)}$-extendible. Note 
  that $\kappa$ is then super-$C^{(k)}$-extendible by \LemmaAof{p-ext-super-ext-4} and hence 
  we have $\kappa\in C^{(k+2)}$ by the induction hypothesis. Let $\psi(x_0\ctenten)$ be a
  $\Pi_{(k+1)+2}$ formula. Say, $\psi(x_0\ctenten)=\exists x\varphi(x,x_0\ctenten)$ where 
  $\varphi$ is a $\Pi_{k+2}$ formula. Then \ixitema[x-ext-super-ext-4-a] $\varphi$ is absolute over $V_\kappa$.

  In particular, if $V_\kappa\models\exists x\varphi(x,a_0\ctenten)$ for
  $a_0\ctenten\in V_\kappa$, then $\uniV\models\exists x\varphi(x,a_0\ctenten)$.

  Now suppose that $\uniV\models\exists x\varphi(x,a_0\ctenten)$ for
  $a_0\ctenten\in V_\kappa$. Let $\alpha>\kappa$ be 
  \st\ \ixitema[x-ext-super-ext-4-0] $V_\alpha\prec_{\Sigma_{k+2}}\uniV$ and there is $c\in V_\alpha$ \st\
  $\uniV\models\varphi(c,a_0\ctenten)$. Since $\kappa$ is super-$C^{(k+1)}$-extendible, 
  there are \ixitema[x-ext-super-ext-5] $\beta\in C^{(k+1)}$ and
  $\Elembed{j}{V_\alpha}{V_\beta}{\kappa}$ \st\ $j(\kappa)>\alpha$. 

  Since $j(a_0)=a_0\ctenten$  and by \xitemof{x-ext-super-ext-5}, we have
  $V_\beta\models\varphi(c,j(a_0)\ctenten)$. It follows that
  $V_\beta\models\exists x\in V_{j(\kappa)}\ \varphi(x,j(a_0)\ctenten)$. By elementarity 
  of $j$, $V_\alpha\models\exists x\in V_{\kappa}\ \varphi(x,a_0\ctenten)$.

  Let $c'\in V_\kappa$ be \st\ $V_\alpha\models \varphi(c',a_0\ctenten$). By 
  \xitemof{x-ext-super-ext-4-0}, $\uniV\models \varphi(c',a_0\ctenten)$.
  Thus, by \xitemof{x-ext-super-ext-4-a}, it follows that
  $V_\kappa\models \varphi(c',a_0\ctenten)$, and hence
  $V_\kappa\models \psi(a_0\ctenten)$. 
  \qedofProp\qedskip

  \noindent\prfof{\Thmof{p-ext-super-ext-3}}
  \assertof{b} $\Rightarrow$ \assertof{b'}: Suppose that $\kappa$ is 
  super-$C^{(n)}$-extendible. and $\Elembed{j}{V_\alpha}{V_\beta}{\kappa}$ is \st\ 
	\ixitema[x-ext-super-ext-6] $V_\alpha\prec_{\Sigma_n}\uniV$, \ixitemar[x-ext-super-ext-7]
  $V_\beta\prec_{\Sigma_n}\uniV$, and $j(\kappa)>\alpha$. It is enough to show that
  $j(\kappa)\in C^{(n)}$.

  By \Propof{x-ext-super-ext-5},\assertof{2}, \ixitema[x-ext-super-ext-8] $\kappa\in C^{(n)}$. Thus
  $V_\alpha\models\kappa\in C^{(n)}$ by \xitemof{x-ext-super-ext-6}. By elementarity 
  of $j$, it follows that \ixitema[x-ext-super-ext-9] $V_\beta\models j(\kappa)\in C^{(n)}$. 
  Thus $\uniV\models j(\kappa)\in C^{(n)}$ by \xitemof{x-ext-super-ext-7}. \smallskip

  \assertof{b'} $\Rightarrow$ \assertof{a}: is clear by definitions. \smallskip

  \assertof{a} $\Rightarrow$ \assertof{c}: Suppose that $\kappa$ is $C^{(n)}$-extendible. 
  Let $\lambda>\kappa$ and let $\alpha\in C^{(n+1)}\setminus\lambda+1$. 
  By $C^{(n)}$-extendibility of $\kappa$, there are $\beta$ and $j$ \st\
  $\Elembed{j}{V_\alpha}{V_\beta}{\kappa}$,\\  \ixitema[x-ext-super-ext-10]
  $j(\kappa)\in C^{(n)}$, and $j(\kappa)>\alpha$. 

  Let $U:=\setof{X\subseteq\Pkl{}{}}{j\imageof\lambda\in j(X)}$. Then $U$ is a normal 
  ultrafilter. Let $M$ be the Mostowski collapse of $\fnsp{\Pkl{}{}}{\uniV}/U$, and
  $\Elembed{i}{\uniV}{M}{\kappa}$ where $i:=j_U$. Let
  $\mapping{k}{{V_{i(\alpha)}}^M}{V_\beta}$; $[f]_U\mapsto j(f)(j\imageof{\lambda})$. Then
  $k\circ i\restr V_\alpha=j$,  $k$ is an elementary embedding and $crit(k)>\alpha$ (see 
  Proposition 22.11 and Lemma 22.12 in \cite{higher-inf}).

  The following claim shows that 
  $i$ witnesses the supercompactness for $C^{(n+1)}$ of $\kappa$.
  \begin{Claim}\Label{cl-ext-super-ext-2}
    $C^{(n+1)}\cap\lambda= (C^{(n+1)})^M\cap\lambda$.
  \end{Claim}

  \noindent\prfofClaim
  For $\gamma<\lambda$, $\gamma\in C^{(n+1)}$\ \   $\Leftrightarrow$\ \ 
  $V_{j(\kappa)}\models\gamma\in C^{(n+1)}$\\
  \mbox{}\hfill\scalebox{0.75}{[\![ ``$\Rightarrow$'': By 
    \xitemof{x-ext-super-ext-10}, and since ``$\gamma\in C^{(n+1)}$'' is $\Pi^{n+1}$, see 
    \cite{bagaria-Cn}, p.214. ``$\Leftarrow$'': 
    By $\gamma<\lambda<\alpha<j(\kappa)$. ]\!]}\medskip

  $\Leftrightarrow$\ \ 
  ${V_{i(\kappa)}}^M\models\gamma\in C^{(n+1)}$
  \hfill\scalebox{0.75}{[\![ By elementarity of $k$, and 
      since $crit(k)>\lambda>\gamma$. ]\!]}\medskip

  $\Leftrightarrow$\ \ $M\modelof{V_{i(\kappa)}\models\gamma\in C^{(n+1)}}$\ \
  $\Leftrightarrow$\ \  $M\models\gamma\in C^{(n+1)}$.\\ \mbox{}\hfill\scalebox{0.75}{[\![ 
      $M\modelof{i(\kappa)\mbox{ is }C^{(n)}\mbox{-extendible}}$ by elementarity of $i$. 
      Thus $M\models i(\kappa)\in C^{(n+2)}$ by \PropAof{p-ext-super-ext-5},\,\assertof{1}. ]\!]}\\
  \qedofClaim\qedskip

  \assertof{c} $\Rightarrow$ \assertof{b}: Suppose that $\kappa$ is supercompact for
  $C^{(n+1)}$. Let $\lambda_*>\kappa$. We have to show that there are $\lambda_0$, 
  $\mu_0$, $j_0$  \st\ $\lambda_0$, $\mu_0\in C^{(n)}$, $\lambda_0>\lambda_*$, 
  $\Elembed{j_0}{V_{\lambda_0}}{V_{\mu_0}}{\kappa}$, and $j_0(\kappa)>\lambda_0$.

  Let $\delta>\lambda>\lambda_*$ be \st\ $\lambda\in C^{(n+1)}$. By assumption on
  $\kappa$, there are $j$, $M\subseteq\uniV$ \st\ $\Elembed{j}{\uniV}{M}{\kappa}$,
  $j(\kappa)>\delta$, \ixitema[x-ext-super-ext-10-0] $\fnsp{\delta}{M}\subseteq M$, $j(\kappa)>\delta$, and 
\ixitema[x-ext-super-ext-11] $C^{(n+1)}\cap\delta=(C^{(n+1)})^M\cap\delta$. Note that 
  $\lambda$ is a strong limit (this follows from $V_\lambda\prec_{\Sigma_1}\uniV$ see e.g.\ 
  \cite{bagaria-Cn}), and hence $\cardof{V_\lambda}=\lambda<\delta$. Thus,
  $V_\lambda\in M$ by \xitemof{x-ext-super-ext-10-0}. We also have
  \ixitema[x-ext-super-ext-12] $M\modelof{\lambda\in C^{(n+1)}}$ by \xitemof{x-ext-super-ext-11}. 

  It follows that $\lambda$, $j(\lambda)$, $j\restr V_\lambda$ are witnesses of the
  $\Sigma_{n+1}$-statement $\sigma(\lambda_*,\kappa)$:
  \begin{xitemize}
  \xitemx[]
    $\exists\ul{\lambda}\,\exists\ul{\mu}\,\exists\ul{i}\ (\,\ul{\lambda},\ul{\mu}\in C^{(n)}
    \ \land\ \ul{\lambda}>\lambda_*\ \land\ \Elembed{\ul{i}}{V_{\ul{\lambda}}}{V_{\ul{\mu}}}{\kappa}
    \ \land\ \ul{i}(\kappa)>\ul{\lambda}\,)$ 
  \end{xitemize}
  in $M$. By \xitemof{x-ext-super-ext-12}, it follows that
  $M\modelof{V_\lambda\models\sigma(\lambda_*,\kappa)}$, and hence
  $V_\lambda\models\sigma(\lambda_*,\kappa)$. Since $\lambda\in C^{(n+1)}$, it follows that
\ixitema[x-ext-super-ext-13] $V\models\sigma(\lambda_*,\kappa)$. The 
  witnesses $\lambda_0$, $\mu_0$, $j_0$ of \xitemof{x-ext-super-ext-13} are as desired. 
  \qedof{\Thmof{p-ext-super-ext-3}}\\
\fi}

\section{Models with super-\texorpdfstring{$C^{(\infty)}$}{C(∞)}-extendible cardinals}
\Label{ext}
We prove that there are unboundedly many super-$C^{(\infty)}$-extendible cardinals  
in $V_\kappa$ below an almost-huge cardinal $\kappa$ (\Corof{p-ext-2}). 

For a cardinal $\kappa$, we say that $V_\kappa$ satisfies the {\It Second-order Vop\v{e}nka's 
Principle} if for any set $C\subseteq V_\kappa$ of structures of the same signature with
$C\not\in V_\kappa$ (which is not necessarily a definable subset of $V_\kappa$), there are 
non-isomorphic $\gmA$, $\gmB\in C$ \st\ we have $\elembed{i}{\gmA}{\gmB}$ 
for an elementary embedding $i$. 

The following is well-known (see e.g.\ Jech \cite{millennium-book}, Lemma 20.27), and attributed 
to William C.\ Powell.
\begin{Lemma}{\rm(W.C.\,Powell \cite{powell})}\Label{p-ext-0} If $\kappa$ is an almost-huge 
  cardinal then $V_\kappa$  
  satisfies the Second-order Vop\v{e}nka's Principle. 
\end{Lemma}
\memo{Scan\_2024-11-23--10.44 extendible - annotated, p.41}
\prf Suppose that $C\subseteq V_\kappa$ where $C$ is a set of structures of the same 
signature, and $C\not\in V_\kappa$. \Wolog, we may assume that \ixitem[x-ext-a] $C$ is closed \wrt\ 
isomorphism. Then it is enough to show that there are non-isomorphic
$\gmA$, $\gmB\in C$ \st\ $\gmA\prec\gmB$. Note that $\rank(C)=\kappa$. 

Let $\Elembed{j}{\uniV}{M}{\kappa}$ be an almost-huge elementary embedding (i.e.\ $M$ 
satisfies \ixitem[x-ext-0] $\fnsp{j(\kappa)\GT}{M}\subseteq M$\/).
Let $\gmA\in j(C)\setminus C$ --- note that $j(C)\setminus C\not=\emptyset$  since
$M\modelof{rank(j(C))=j(\kappa)>\kappa}$. Let $A$ be the underlying set of the structure
$\gmA$. 

We have \ixitem[x-ext-1] $M\models j(\gmA)\not\cong\gmA$ --- otherwise $M\models j(\gmA)\cong\gmA\in j(C)$,  and hence
$\uniV\models\gmA\in C$ by elementarity. This is a contradiction to the choice of $\gmA$. 

Let $\gmA':=j(\gmA)\restr j\imageof{A}$. 
\begin{Claim}
  \Label{cl-ext-0} \wassertof{1} $M\models\gmA'\in j(C)$. \smallskip

  \wassert{2} $M\models\gmA'\prec j(\gmA)$.
\end{Claim}
\prfofClaim  \assertof{1}: $\gmA'\in M$, and $M\models\gmA'\cong\gmA$ by \xitemof{x-ext-0}. 
Since $gmA\in j(C)$ by the choice of $\gmA$,   
\xitemof{x-ext-a} and the elementarity of $j$ imply $M\models\gmA'\in j(C)$. \smallskip

\assertof{2}: Working in $M$, we check that $\gmA'\subseteq j(\gmA)$ satisfy Vaught's criterion.

Suppose $a_0\ctentenc a_{n-1}\in j\imageof{A}$, $a\in j(A)$ and
$j(\gmA)\models\varphi(a,a_0\ctentenc a_{n-1})$. Let $a'_0\ctentenc a'_{n-1}\in A$ be \st\
$a_0=j(a'_0)\ctentenc a_{n-1}=j(a'_{n-1})$. Since
$M\models\exists\ul{a}\in j(A)\ j(\gmA)\models\varphi(\ul{a},j(a'_0)\ctentenc j(a'_{n-1}))$, 
it follows that 
$\uniV\models\exists\ul{a}\in A\ \gmA\models\varphi(\ul{a},a'_0\ctentenc a'_{n-1})$. Let
$a'\in A$ be \st\ $\uniV\models\gmA\models\varphi(a',a'_0\ctentenc a'_{n-1})$. Then
$j(a')\in j\imageof{A}$, and
$M\models j(\gmA)\models\varphi(j(a'),j(a'_0)\ctentenc j(a'_{n-1}))$ by elementarity, as desired. 
\qedofClaim
\qedskip

Now by \xitemof{x-ext-1} and \Claimof{cl-ext-0}, \assertof{2},
$M\modelof{\xmbox{there are non-isomorphic }\gmA,\gmB\in j(C)\mbox{ \st\ }\gmA\prec\gmB}$. By 
elementarity it follows that 
$\uniV\modelof{\xmbox{there are non-isomorphic }\gmA,\gmB\in C\mbox{ \st\ }\gmA\prec\gmB}$. 
\qedofLemma
\qedskip

To prove the existence of a super-$C^{(\infty)}$-Lg extendible cardinal with an 
appropriate Laver function (\Lemmaof{p-Lg-ext-0}), we need the existence of a cardinal 
stronger than a super-$C^{(\infty)}$-Lg 
extendible which we call below strongly super-$C^{(\infty)}$-Lg extendible. 
Similarly to the existence of a super-$C^{(\infty)}$-Lg extendible cardinal, the existence 
of this storger large cardinals also is not formalizable in the first-order logic.

For an inaccessible $\mu$ and $\kappa<\mu$, we say that $\kappa$ is strongly 
super-$C^{(\infty)}$-extendible, and denote this as 
{\It$V_\mu\modelof{\kappa\mbox{ is strongly super-}C^{(\infty)}\mbox{-extendible}}$},
if 
for any $\kappa<\lambda<\mu$, there are $\lambda<\lambda'<\lambda''$ and $j$ \st\
$V_{\lambda'}\prec V_{\lambda''}\prec V_\mu$, and 
$\Elembed{j}{V_{\lambda'}}{V_{\lambda''}}{\kappa}$.\memo{--- Gaifman's theorem 
  (\href{https://fuchino.ddo.jp/notes/math-notes-20.pdf\pdfpage{60}}{math-notes-20,p.60}) makes 
  some of the following arguments much simpler.}

Note that by \Propof{p-ext-super-ext-2},
$V_\mu\modelof{\kappa\mbox{ is strongly super-}C^{(\infty)}\mbox{-extendible}}$ 
implies
$V_\mu\modelof{\kappa\mbox{ is super-}C^{(\infty)}\mbox{-extendible}}$. 

\begin{Prop}\Label{p-ext-1}
  \memo{Check if Mahloness follows from second-order Vopěnka.}
  Suppose that $\mu$ is an inaccessible cardinal, and \ixitemr[x-ext-2] $V_\mu$ satisfies 
  the Second-order  
  Vopěnka's Principle. Then there are unboundedly many $\kappa<\mu$ \st\
  $V_\mu\modelof{\kappa\mbox{ is strongly super-}C^{(\infty)}\mbox{-extendible}}$.
\end{Prop}
\prf\memo{Scan\_2024-11-23--10.44 extendible - annotated, (p.38) p.40}
Suppose $\beta^*<\mu$. We want to show that there is $\beta^*\leq\kappa<\mu$ \st\
$V_\mu\modelof{\kappa\mbox{ is strongly super-}C^{(\infty)}\mbox{-extendible}}$.

Let ${\darkred C}:=\setof{\alpha<\mu}{V_\alpha\prec V_\mu}$, and
$\calV_\mu=\pairof{V_\kappa,\in, C}$. 

Let 
\begin{xitemize}
\xitemx[] 
  ${\darkred I}:=\setof{\alpha<\mu}{
  \alpha\xmbox{ is an }\omega\xmbox{-limit of ordinals }\eta\mbox{ \st\ }
  \calV_\mu\restr V_\eta\prec\calV_\mu}$.
\end{xitemize}
Note that $I$ is cofinal in $\mu$ and we have $\calV_\mu\restr V_\alpha\prec\calV_\mu$ for all $\alpha\in I$. 

For each $\alpha\in I$, let ${\darkred C_\alpha}\subseteq\alpha$ be a cofinal subset of $\alpha$ of 
order-type $\omega$ consisting of (increasing) 
$\eta^\alpha_n$, $n\in\omega$ with $\calV_\mu\restr V_{\eta^\alpha_n}\prec\calV_\mu$. 
In particular $C_\alpha\subseteq C$. 
Let  
\begin{xitemize}
\xitemx[] $\calC:=\setof{\pairof{V_{\alpha},\in, C_\alpha,\xi}_{\xi<\beta^*}}{{}
  \begin{array}[t]{@{}l}
    \alpha\in I}.
  \end{array}$
\end{xitemize}

By \xitemof{x-ext-2}, there are $\alpha$, $\beta\in I$, $\alpha<\beta$ \st\ letting 
$\calW_{\alpha}=\pairof{V_{\alpha},\in,C_\alpha,\xi}_{\xi<\beta^*}$, and 
$\calW_{\beta}=\pairof{V_{\beta},\in,,C_\beta,\xi}_{\xi<\beta^*}$,
there is an elementary 
embedding $\elembed{i}{\calW_{\alpha}}{\calW_{\beta}}$. Since 
$i\imageof{C_\alpha}=C_\beta$ by elementarity, $i\restr\alpha$ is not an identity mapping. 

Let 
${\darkred\kappa}:=\crit(i)$. Then $\beta^*\leq\kappa<\alpha$ by virtue of the constants
$\xi<\beta^*$ in the structures. 

For all $k\in\omega$ \st\ $\eta^\alpha_k>\kappa$, we have
$\Elembed{i\restr V_{\eta^\alpha_k}}{V_{\eta^\alpha_k}}{V_{\eta^\beta_k}}{\kappa}$. Thus 

\begin{xitemize}
\xitemx[] $\calV_\mu\modelof{\exists\ul{\eta}\exists\ul{i}\,
  (\ul{\eta}\in \symb{C}\ \land\ \Elembed{\ul{i}}{V_{\eta^\alpha_k}}{V_{\ul{\eta}}}{\kappa})
  }$.
\end{xitemize}
By elementarity (in the extended language), 
\begin{xitemize}
\xitemx[] 
  $\calV_\mu\restr V_\alpha\modelof{\exists\ul{\eta}\exists\ul{i}\,
    (\ul{\eta}\in \symb{C}\ \land\ \Elembed{\ul{i}}{V_{\eta^\alpha_k}}{V_{\ul{\eta}}}{\kappa})
  }$
\end{xitemize}
for all $k\in\omega$ \st\ $\eta^\alpha_k>\kappa$. 

It follows that 
\begin{xitemize}
\xitemx[] $\calV_\mu\restr V_\alpha\models\forall\ul{\nu}\exists\ul{\eta}_0\exists\ul{\eta}\exists\ul{i}\ 
  (\ul{\eta}_0>\ul{\nu},\kappa\ \land\ \ \ul{\eta}_0,\ul{\eta}\in\symb{C}\ \land\ 
  \Elembed{\ul{i}}{V_{\ul{\eta}_0}}{V_{\ul{\eta}}}{\kappa})$. 
\end{xitemize}

Thus by elementarity (in the extened language), we have:
\begin{xitemize}
\xitemx[] $\calV_\mu\models\forall\ul{\nu}\exists\ul{\eta}_0\exists\ul{\eta}\exists\ul{i}\ 
  (\ul{\eta}_0>\ul{\nu},\kappa\ \land\ \ \ul{\eta}_0,\ul{\eta}\in\symb{C}\ \land\ 
  \Elembed{\ul{i}}{V_{\ul{\eta}_0}}{V_{\ul{\eta}}}{\kappa})$. 
\end{xitemize}
This simply means $V_\mu\modelof{\kappa\mbox{ is strongly super-}C^{(\infty)}\mbox{-extendible}}$.
\qedofProp\qedskip

\begin{Cor}\Label{p-ext-2} Suppose that $\mu$ is almost-huge. Then there are 
  unboundedly many $\kappa<\mu$ \st\
  $V_\mu\modelof{\kappa\mbox{ is strongly super-}C^{(\infty)}\mbox{-extendible}}$. 
\end{Cor}
\prf By \Lemmaof{p-ext-0} and \Propof{p-ext-1}. \qedofThm\qedskip

The following can be proved by the same argument as that of \Propof{p-ext-super-ext-2}.

\begin{Lemma}\Label{p-ext-1-0} Suppose that $\mu$ is an inaccessible 
  cardinal, $\kappa<\mu$.  
  Then $V_\mu\modelof{\kappa\mbox{ is strongly super-}C^{(\infty)}\mbox{-extendible}}$ 
  is equivalent to the condition that, for any $\kappa<\lambda<\mu$ there are
  $\lambda<\lambda'<\lambda''$, $j$ and $M$ \ \st\  $\lambda''=j(\lambda')$,
  $V_{\lambda'}\prec V_{\lambda''}\prec V_\mu$, 
  $\Elembed{j}{V_{\mu}}{M}{\kappa}$, $M\subseteq V_\mu$, $j(\kappa)>\lambda'$ and
    $V_{\lambda''}\in M$.\qed
\end{Lemma}

\section{LgLCAs and 
        super-\texorpdfstring{$C^{(\infty)}$}{C(∞)}-LgLCAs for extendibility imply (almost) everything}\Label{everyth}
In this section, we prove that Laver-generic Large Cardinal Axioms (LgLCAs) and
super-$C^{(\infty)}$-LgLCAs for extendibility  
imply the strongest forms of resurrection axioms, maximality principles, and absoluteness 
known to be consistent under some large cardinal consistency. 
This 
was previously known to hold under LgLCAs for hugeness and super-$C^{(\infty)}$-LgLCAs for 
hyperhugeness. 

We begin with a short summary of definitions and known results around LgLCAs and
super-$C^{(\infty)}$-LgLCAs. 

Laver-generic large cardinals were introduced in \cite{sfetal-II}. For a class 
$\calP$ of \pos\ and a notion $\LC$ of large cardinal,  a cardinal 
$\kappa$ is said to be $\calP$-Laver generic $\LC$ if the statement about the existence of 
elementary embedding $\Elembed{j}{\uniV}{M}{\kappa}$ for $j$, $M\subseteq\uniV$ with the 
closedness condition $C_\LC$ of $M$ in the definition of the notion $\LC$ of large cardinal are 
replaced with the statement:
\begin{xitemize}
\xitem[x-everyth-0] for any $\poP\in\calP$, there is a $\poP$-name $\utpoQ$ \st\ 
$\forces{\poP}{\utpoQ\in\calP}$, and for $\poP\ast\utpoQ$-generic $\genH$ there are $j$,
  $M\subseteq\uniV[\genH]$ \st\ $\poP$, $\poP\ast\utpoQ$, $\genH\in M$, 
  $\Elembed{j}{\uniV}{M}{\kappa}$, and $M$ satisfies $C'_\LC$ which is the generic large cardinal 
  variant of the closedness property $C_\LC$ associated with the notion $\LC$ of large 
  cardinal.
\end{xitemize}

For supercompactness, the instance of \xitemof{x-everyth-0} for an iterable $\calP$ is as 
follows: a cardinal $\kappa$ is {\It$\calP$-Laver-generically supercompact\/} 
({\It $\calP$-Lg supercompact} for short) if,
\begin{xitemize}
\xitem[x-everyth-1] for any
  $\lambda>\kappa$, and for any $\poP\in\calP$, there is a $\poP$-name $\utpoQ$ \st\
  $\forces{\poP}{\poQ\in\calP}$, and, for any $(\uniV,\poP\ast\utpoQ)$-generic $\genH$, there 
  are $j$, $M\subseteq\uniV[\genH]$ 
  \st\ $\Elembed{j}{\uniV}{M}{\kappa}$, $j(\kappa)>\lambda$, $\poP$, $\poP\ast\utpoQ$,
  $\genH\in M$, $j\imageof{\lambda}\in M$.  
\end{xitemize}

Note that, in \xitemof{x-everyth-1}, the closure property
``$\fnsp{\lambda}{M}\subseteq M$'' in the usual definition of supercompactness is replaced 
with ``$j\imageof\lambda\in M$''. For a genuine elementary embedding introduced by some 
ultrafilter, these two conditions are equivalent (see e.g.\ Kanamori \cite{higher-inf}, 
Proposition 22.4, (b)). This equivalence is no more valid in general for generic 
embeddings. Nevertheless, the condition ``$j\imageof\lambda\in M$'' can be still 
considered as a certain closure property (see Lemma 3.5 in Fuchino-Rodrigues-Sakai 
\cite{sfetal-II}). 

We say that a $\calP$-Lg supercompact  cardinal $\kappa$ is {\It tightly
    $\calP$-Laver-generically supercompact\/} ({\It tightly $\calP$-Lg supercompact}, for 
  short) if additionally, we have 
  $\cardof{RO(\poP\ast\utpoQ)}\leq j(\kappa)$.
 
A (tightly) $\calP$-Lg supercompact cardinal is often decided uniquely as the 
cardinal ${\It\kappa_\refl}:=\sup(\ssetof{\continuum,\aleph_2})$. This is the case, if $\calP$ is 
the class of all $\sigma$-closed \pos.  Then \CH\ holds under the existence of 
a $\calP$-Lg supercompact $\kappa$ and $\kappa=\aleph_2$ ($=\kappa_\refl$). 

Similarly, if $\calP$ is either the class of all proper \pos\ or the class of all 
semi-proper \pos, the existence of a $\calP$-generically supercompact $\kappa$ implies
$\continuum=\aleph_2$ and again $\kappa=\kappa_\refl$. 

For the case that $\calP$ is the 
class of all ccc \pos, it is open whether a $\poP$-Lg supercompact 
cardinal is decided to be $\kappa_\refl$. However a tightly $\calP$-Lg 
supercompact cardinal under the present definition of tightness\footnote{In course of the 
  development of the theory of Laver-genericity, we strengthened the definition of 
  tightness. However, the modification is chosen so that it still holds in all the standard 
  models of Laver genericity.  } is the  
continuum ($=\kappa_\refl$) and, in this case,  the continuum is extremely large. There is 
a more general theorem which suggests that for a ``natural'' class $\calP$ of \pos, 
the existence of (tightly) $\calP$-Lg supercompact cardinal implies that 
the continuum is either $\aleph_1$ or $\aleph_2$ or else extremely large (see \cite{sfetal-II}, 
\cite{future}, \cite{janos}). 

The naming ``Laver-generic ...'' based on the fact that the standard models with this type of 
generic large cardinal is created by starting from a large cardinal, and then iterating along 
with a Laver function for the large cardinal with the support appropriate for the class of 
\pos\ in consideration. This is exactly the way to create models of 
\PFA\ and \MM. Actually, for $\calP$ being the class of all proper \pos\ or the class of 
all 
semi-proper \pos, the existence of a $\calP$-Lg supercompact cardinal implies 
the double-plus version of the corresponding forcing axiom (see \Thmof{p-square-3} below), 
and can be considered as an  
axiomatization of the standard models of such axioms.  

In the following, we call the axiom asserting that the cardinal $\kappa_\refl$ is a/the 
tightly $\calP$-Laver  
generic $\LC$,  the {\It$\calP$-Laver-generic large cardinal axiom for the notion of large 
cardinal $\LC$} ({\It the $\calP$-LgLCA for $\LC$}, for short). 

The instances of $\calP$-LgLCAs for other notions of large cardinal considered in
\cite{future}, \cite{janos}, \cite{FGP}, \cite{sfetal-II}, \cite{FuOt}, \cite{recurrence}, 
etc.\ are summarized in the 
following chart.\bigskip

{
\begin{tabular}{| l || l | l |}\hline
  \vbox{\hbox{\rule[-3pt]{0pt}{2.34ex}The $\calP$-LgLCA\ for}\hbox{ }}
  & \vbox{\hbox{\rule[-3pt]{0pt}{2.8ex}The condition\quad
      ``$j\imageof{\lambda}\in M$''\quad in the definition of}
    \hbox{``the $\calP$-LgLCA for super-compact'' is replaced with: }} \\
  \hline
  \rule{0pt}{2.16ex}hyperhuge & $j\imageof j(\lambda)\in M$\\[\jot] 
  ultrahuge & $j\imageof j(\kappa)\in M$ and ${V_{j(\lambda)}}^{\uniV[\genH]}\in M$\\[\jot] 
  superhuge & $j\imageof j(\kappa)\in M$\\[\jot] 
  super-almost-huge & $j\imageof \mu\in M$ for all $\mu<j(\kappa)$\\[\jot]
  \hline
\end{tabular}
}\bigskip

It has been proved that LgLCAs for sufficiently strong notions of large cardinal imply strongest 
forms of resurrection, maximality and absoluteness among the known consistent variants of 
resurrection, maximality and absoluteness:

\begin{xitemize}
\xitem[x-everyth-2] 
  In \cite{future}, it is proved that a boldface variant of Resurrection Axiom 
  by Hamkins and Johnstone (\cite{hamkfins-johnstone}, \cite{hamkins-johnstone2}) for
  $\calP$ and parameters from $\calH(\kappa_\refl)$ follows from the $\calP$-LgLCA for ultrahuge.\smallskip

\xitem[x-everyth-3] 
  In \cite{recurrence} or \cite{janos} (see \cite{FGP} for an improved version, or see also 
  \Thmof{p-Lg-RcA-0} below), it is 
  proved that the $\calP$-LgLCA for ultrahuge implies  
  a restricted form of Maximality Principle for $\calP$ and $\calH(\kappa_\refl)$. More 
  specifically, It is proved that under the $\calP$-LgLCA for ultrahuge implies
  $(\calP,\calH(\kappa_\refl))_{\Sigma_2}$-\RcAp\ holds --- see below for the definition of 
  this priniple (Theorem 21 in \cite{janos}).\memo{\textbackslash Tmof\{p-Lg-RcA-0\} in \cite{janos}}
\end{xitemize}

It can be shown that LgLCA type axiom formulated in a single formula is incapable of covering the 
full Maximal Principle (\cite{future}). The notion of super-$C^{(\infty)}$ LgLCAs is 
introduced in \cite{recurrence} to fill this gap.

For a notion \LC\ of large cardinal let $C'_\LC$ be the closedness property of the target 
model of the generic large cardinal corresponding to \LC. 
We call a cardinal $\kappa$ {\It tightly super-$C^{(\infty)}$-$\calP$-Laver generically \LC} 
if, for any $n\in\natnums$,  
$\lambda_0>\kappa$, and $\poP\in\calP$, there are $\lambda\geq\lambda_0$ and a $\poP$-name $\utpoQ$ \st\
$V_\lambda\prec_{\Sigma_n}\uniV$, $\forces{\poP}{\utpoQ\in\calP}$ and for any
$(\uniV,\poP\ast\utpoQ)$-generic $\genH$, there 
are $j$, $M\subseteq\uniV[\genH]$ \st\ ${V_{j(\lambda)}}^{\uniV[\genH]}\prec_{\Sigma_n}\uniV[\genH]$, 
$\Elembed{j}{\uniV}{M}{\kappa}$, $j(\kappa)>\lambda$,
$C'_\LC$, and $\cardof{RO(\poP\ast\utpoQ)}\leq j(\kappa)$.

The {\It super-$C^{(\infty)}$-$\calP$-Laver-generic large cardinal axiom for the notion 
\LC\ of large cardinal} ({\It the super-$C^{(\infty)}$-$\calP$-LgLCA\ for \LC}, for short) 
is the assertion that  
$\kappa_\refl$ is a/the tightly super $C^{(\infty)}$-$\calP$-Laver-generic extendible cardinal.

Note that tightly super-$C^{(\infty)}$-$\calP$-Laver generically \LC\ is not formalizable in 
general but the super-$C^{(\infty)}$-$\calP$-LgLCA\ for \LC\ is as an axiom scheme, since the generic large 
cardinal in the axiom is named as $\kappa_\refl$. 

It is shown that for a transfinitely iterable $\calP$, the consistency of 
the super-$C^{(\infty)}$-$\calP$-LgLCA for hyperhuge follows from a 2-huge cardinal 
(\cite{recurrence}, Lemma 2.6 and Theorem 2.8). \memo{\textbackslash Lemmaof\{p-Lg-RcA-2\},\textbackslash Thmof\{p-Lg-RcA-4\}}

\begin{xitemize}
\xitem[x-everyth-3-0] 
  The super-$C^{(\infty)}$-$\calP$-LgLCA for ultrahuge implies the full Maximality Principle for $\calP$ 
  and $\calH(\kappa_\refl)$ (\cite{recurrence}, \href{https://fuchino.ddo.jp/papers/recurrence-axioms-x.pdf\pdfpage{23}}{Theorem 4.10}).\smallskip
\end{xitemize}

\begin{xitemize}
\xitem[x-everyth-4] 
  In \cite{FGP}, it is proved that under $\BFA_{\LT\kappa_\refl}(\calP)$ and the $\calP$-LgLCA 
  for huge, a generalization of Viale's Absoluteness Theorem in Viale \cite{viale-revisited} 
  holds (see Theorem 5.7 in \cite{FGP}).
\end{xitemize}


For extendibility and super-$C^{(\infty)}$-extendibility, the natural 
Laver-generic versions of these notions of large cardinals should be the following: 
a cardinal $\kappa$ is {\It tightly $\calP$-Laver generically extendible} ({\It tightly
  $\calP$-Lg extendible}, for short)
if, for any
$\lambda>\kappa$, and for any $\poP\in\calP$, there is a $\poP$-name $\utpoQ$ \st\
$\forces{\poP}{\utpoQ\in\calP}$ and for any $(\uniV,\poP\ast\utpoQ)$-generic $\genH$, there 
are $j$, $M\subseteq\uniV[\genH]$ 
\st\ $\Elembed{j}{\uniV}{M}{\kappa}$, $j(\kappa)>\lambda$, 
${V_{j(\lambda)}}^{\uniV[\genH]}\in M$, and $\cardof{RO(\poP\ast\utpoQ)}\leq j(\kappa)$. 

The {\It $\calP$-Laver-generic large cardinal axiom} for the notion of 
extendibility ({\It the $\calP$-LgLCA\ for extendible}, for short) is the assertion that 
$\kappa_\refl$ is a/the tightly $\calP$-Lg extendible cardinal. 

A cardinal $\kappa$ is {\It tightly super-$C^{(\infty)}$-$\calP$-Laver generically extendible} 
({\It tightly super-$C^{(\infty)}$-Lg extendible}, for short)
if, for any $n\in\natnums$,  
$\lambda_0>\kappa$, and $\poP\in\calP$, there are $\lambda\geq\lambda_0$ and a $\poP$-name $\utpoQ$ \st\
$V_\lambda\prec_{\Sigma_n}\uniV$, $\forces{\poP}{\utpoQ\in\calP}$ and for any
$(\uniV,\poP\ast\utpoQ)$-generic $\genH$, there 
are $j$, $M\subseteq\uniV[\genH]$ \st\ ${V_{j(\lambda)}}^{\uniV[\genH]}\prec_{\Sigma_n}\uniV[\genH]$, 
$\Elembed{j}{\uniV}{M}{\kappa}$, $j(\kappa)>\lambda$,
${V_{j(\lambda)}}^{\uniV[\genH]}\in M$, and $\cardof{RO(\poP\ast\utpoQ)}\leq j(\kappa)$.

The {\It super-$C^{(\infty)}$-$\calP$-Laver-generic large cardinal axiom} for the notion of 
extendibility ({\It the super-$C^{(\infty)}$-$\calP$-LgLCA\ for extendible}, for short) is 
the assertion that  
$\kappa_\refl$ is a/the tightly super $C^{(\infty)}$-$\calP$-Lg extendible cardinal.
\bigskip

{
\begin{tabular}{| l || l | l |}\hline
  \vbox{\hbox{\rule[-3pt]{0pt}{2.34ex}The $\calP$-LgLCA\ for}\hbox{ }}
  & \vbox{\hbox{\rule[-3pt]{0pt}{2.8ex}The condition\quad
      ``$j\imageof{\lambda}\in M$''\quad in the definition of}
    \hbox{``the $\calP$-LgLCA for super-compact'' is replaced with: }} \\
  \hline
  \rule{0pt}{2.16ex}hyperhuge & $j\imageof j(\lambda)\in M$\\[\jot] 
  ultrahuge & $j\imageof j(\kappa)\in M$ and ${V_{j(\lambda)}}^{\uniV[\genH]}\in M$\\[\jot] 
  superhuge & $j\imageof j(\kappa)\in M$\\[\jot] 
  super-almost-huge & $j\imageof \mu\in M$ for all $\mu<j(\kappa)$\\[\jot]
  {\darkred extendible} & $\darkred{V_{j(\lambda)}}^{\uniV[\genH]}\in M$\\[\jot]
  \hline
\end{tabular}
}\bigskip

As we have shown an almost-huge cardinal produces a transitive model with cofinally many
super-$C^{(\infty)}$-extendible cardinals. In the next section, we show that we can generic 
extend such models to a model of the super-$C^{(\infty)}$ $\calP$-LgLCA for extendible for each 
reasonable (i.e. transfinitely iterable) class $\calP$ of \pos. Thus the super-$C^{(\infty)}$
$\calP$-LgLCA for extendible  
and the $\calP$-LgLCA for extendible are of relatively low consistency strength.

These LgLCAs are placed at the expected places in the web of implications of LgLCAs (see 
also the chart on page \pageref{the-diagram}):
\begin{Lemma}\Label{p-identitiy-crisis-0}
  Suppose that $\calP$ is an arbitrary class of \pos. \wassertof{1} 
  The $\calP$-LgLCA for hyperhuge implies the $\calP$-LgLCA for extendible. 
  The super-$C^{(\infty)}$-$\calP$-LgLCA for hyperhuge implies 
  the super-$C^{(\infty)}$-$\calP$-LgLCA for extendible. \smallskip

  \wassert{2} The $\calP$-LgLCA for extendible implies the $\calP$-LgLCA for supercompact. 
\end{Lemma}
\prf \assertof{1}: By definition.\smallskip

\assertof{2}: Asdumr that the $\calP$-LgLCA for extendible holds. Let $\kappa:=\kappa_\refl$. 
Suppose that $\lambda>\kappa$ is regular, and $\poP\in\calP$. Then there are $\utpoQ$,
$\genH$, $j$, $M\subseteq\uniV[\genH]$ \st\ $\forces{\poP}{\utpoQ\in\calP}$, $\genH$ is
$(\uniV,\poP\ast\utpoQ)$-generic, $j$, $M\subseteq\uniV[\genH]$, 
\ixitem[x-i-c-0] $\Elembed{j}{\uniV}{M}{\kappa}$, \\
\ixitem[x-i-c-1] $j(\kappa)>\lambda$, 
\ixitemc[x-i-c-2] $\cardof{RO(\poP\ast\utpoQ)}\leq j(\kappa)$, and 
\ixitemr[x-i-c-2-0] ${V_{j(\lambda)}}^{\uniV[\genH]}\in M$. 

$M\modelof{j(\lambda)\mbox{ is regular}}$ by elementarity \xitemof{x-i-c-0}. Thus 
$\uniV\modelof{j(\lambda)\mbox{ is regular}}$. By $j(\kappa)<j(\lambda)$ and 
\xitemof{x-i-c-2}, it follos that $\uniV[\genH]\modelof{j(\lambda)\mbox{ is regular}}$. 
Thus $\uniV[\genH]\models\cf(j\imageof{\lambda})\leq\lambda<j(\lambda)=\cf(j(\lambda))$, 
and $\uniV[\genH]\models j\imageof{\lambda}\in {V_{j(\lambda)}}^{\uniV[\genH]}\in M$.
Hence $j\imageof{\lambda}\in M$ by transitivity of $M$.

This shows that $j$'s for all regular $\lambda>\kappa$ and $\poP\in\calP$ witness that 
$\kappa$ is tightly $\calP$-Lg.\ supercompact. 
\qedofLemma
\qedskip

All of the results mentioned in \xitemof{x-everyth-2} $\sim$ \xitemof{x-everyth-4} can be 
proved under the assumption of the LgLCA for extendible instead of the stronger assumptions in 
the original results. Below we shall  state these results. In the extended 
version of the present paper, we shall include all the details of the proofs (though the 
proofs are practically identical with the original ones) for the convenience of the reader. 
  
The following boldface version of the Resurrection Axioms was studied by Hamkins and 
Johnstone in \cite{hamkins-johnstone2}: 
For a class $\calP$ of \pos\ and a definition $\mu^\bullet$ of a cardinal (e.g.\ as 
  $\aleph_1$, $\aleph_2$, $\continuum$, $(\continuum)^+$. etc.) the {\It Resurrection Axiom 
  in Boldface for 
  $\calP$ and $\calH(\mu^\bullet)$} is defined by:

\begin{xitemize}
\item[$\darkred\BfRA^\calP_{\calH(\mu^\bullet)}$ :] For any $A\subseteq\calH(\mu^\bullet)$ and any 
  $\poP\in\calP$, there is a $\poP$-name $\utpoQ$ 
  of \po\ \st\ $\forces{\poP}{\utpoQ\in\calP}$ and, for
  any $(\uniV,\poP\ast\utpoQ)$-generic $\genH$, there is
  $A^*\subseteq\calH(\mu^\bullet)^{\uniV[\genH]}$ \st\ 
  $(\calH(\mu^\bullet)^\uniV,A,{\in})\prec(\calH(\mu^\bullet)^{\uniV[\genH]},A^*,{\in})$. 
\end{xitemize}

\begin{Thm}\Label{T-Lg-RA-0}{\rm 
    (\href{https://fuchino.ddo.jp/papers/RIMS2022-RA-MP-x.pdf\pdfpage{32}}{Theorem 7.1} in 
    \cite{future} reformulated under the LgLCA for extendible)}
  For an iterable class $\calP$ of \pos, assume that the $\calP$-LgLCA for extendible holds. 
  Then 
  $\BfRA^\calP_{\calH(\kappa_\refl)}$ holds.  \ifextended\else\qed\fi

\end{Thm}
{\ifextended\extendedcolor
\prf 
Let $\kappa:=\kappa_\refl$ and assume that $\kappa$ is tightly $\calP$-Lg extendible. 
Suppose $A\subseteq\calH(\kappa)$ and $\poP\in\calP$. By the 
tightly $\calP$-Lg  extendibility  
of $\kappa$, there is a $\poP$-name $\utpoQ$ \smallskip of a \po\ with
$\forces{\poP}{\utpoQ\in\calP}$ \st, for $(\uniV,\poP\ast\utpoQ)$-generic $\genH$, there 
are $j$, $M\subseteq\uniV[\genH]$ with 
\begin{xitemize}
\xitemA[was-a] $\Elembed{j}{\uniV}{M}{\kappa}$, 
\xitemA[was-b] $j(\kappa)=\cardof{RO(\poP\ast\utpoQ)}$,
\xitemA[was-c] $\poP$, $\genH\in M$, and
\xitemA[was-d] ${V_{j(\lambda)}}^{\uniV[\genH]}\in M$.
\end{xitemize}

\Wolog, we may assume that the underlying set of $\poP\ast\utpoQ$ 
is $j(\kappa)$. 
Since $\crit(j)=\kappa$, $j(a)=a$ for all $a\in(\calH(\kappa))^\uniV$. 

By \xitemof{was-d}, we have 
  $\calH(j(\kappa))^{\uniV[\genH]}\subseteq M$, and hence 
  $j(\calH(\kappa))=\calH(j(\kappa))^M=\calH(j(\kappa))^{\uniV[\genH]}$.
 
Thus, we have
\begin{xitemize}
\item[] 
  $\Elembed{id_{\calH(\kappa)^\uniV}=j\restr\calH(\kappa)^\uniV}
  {\ (\calH(\kappa)^\uniV,A,{\in})\ }{\ (\calH(j(\kappa))^{\uniV[\genH]},j(A),{\in})}{}$.
\end{xitemize}
\qedofThm
\fi}

Recurrence Axioms are introduced in Fuchino and Usuba \cite{recurrence}. 

For an iterable class $\calP$ 
of \pos, a set $A$ (of parameters), and a set $\Gamma$ of $\Lin$-
formulas, {\It$\calP$-Recurrence Axiom$^+$ for formulas in $\Gamma$ with parameters from $A$} 
({\It$(\calP,A)_\Gamma$-\RcAp}, for short) is the following assertion expressed as an axiom 
scheme formulated in $\Lin$:
\begin{xitemize}
\item[{\It$(\calP,A)_\Gamma$-\RcAp\,:}] For any $\varphi(\overline{x})\in\Gamma$ and $\overline{a}\in A$, if
  $\forces{\poP}{\varphi(\overline{a})}$, then there is a $\calP$-ground 
  $\uniW$ of $\uniV$ \st\ $\overline{a}\in \uniW$ and $\uniW\models\varphi(\overline{a})$. 
\end{xitemize}
Here, an inner model $M$ of $\uniV$ is said to be a {\It$\calP$-ground of\/ $\uniV$}, if 
there are $\poP\in M$ and $\genG\in\uniV$ \st\ $M\modelof{\poP\in\calP}$, $\genG$ is an
$(M,\poP)$-generic filter, and $\uniV=M[\genG]$. If $M$ is a $\calP$-ground of $\uniV$ for 
the class $\calP$ of all \pos, we shall say that $M$ is a {\It ground of\/ $\uniV$}. 
The Recurrence Axiom {\It $(\calP,A)_\Gamma$-\RcA}  without $+$ is obtained 
when ``$\calP$-ground'' in the definition of $(\calP,A)_\Gamma$-\RcAp\ is replaced 
with ``ground''. 

If $\Gamma$ is the set of all $\Lin$-formulas, we drop the subscript $\Gamma$ and say  
simply {\It$(\calP,A)$-\RcAp} or {\It$(\calP,A)$-\RcA}.

As it is noticed in \cite{recurrence}, $(\calP,A)$-\RcAp\ is equivalent to the Maximality 
Principle $\MP(\calP, A)$ (see Proposition 2.2,  \assertof{2} in \cite{recurrence}).

When \cite{FGP} was written, we didn't consider the notion of LgLCA for extendible among 
the possible LgLCAs. This is why 
the following theorem was stated there under the assumption of the $\calP$-LgLCA for ultrahuge. 
However the proof given in \cite{FGP} works perfectly under the $\calP$-LgLCA for extendible 
without any change. 
\begin{Thm}{\rm (\href{https://fuchino.ddo.jp/papers/generic-absoluteness-revisited-x.pdf\pdfpage{40}}{Theorem 6.1} in Fuchino, Gappo and Parente \cite{FGP} reformulated under LgLCA for extendible)}
  \Label{p-Lg-RcA-0} Assume the $\calP$-LgLCA for extendible.
  Then
  $(\calP,\calH(\kappa_\refl))_{\Gamma}$-\RcAp\ holds where $\Gamma$ is the set of all 
  formulas which are conjunctions of a $\Sigma_2$-formula and a $\Pi_2$-formula.\ifprivate\else \qed\fi
\end{Thm}
{\ifextended\extendedcolor

For the proof of \Thmabove, we use the following lemma which should be a well-known fact. 
\begin{LemmaA}\extendedcolor
  \Label{p-Lg-RcA-0-0} If $\alpha$ is a limit ordinal and $V_\alpha$ satisfies a 
  sufficiently large finite fragment of \ZFC, then for any $\poP\in V_\alpha$ 
  and $(\uniV,\poP)$-generic $\genG$,  
  we have $V_\alpha[\genG]={V_\alpha}^{\uniV[\genG]}$. 
\end{LemmaA}
\memox{connect this Lemma with recurrence5? Corollary 25.
to Ikegami + Trang Theorem 1.8}
{
\ifextended
\prf ``$\subseteq$'': This inclusion holds without the condition  
on the fragment of \ZFC. Also the condition ``$\poP\in V_\alpha$'' is irrelevant for this 
inclusion. 

We show by induction on $\alpha\in\On$ that
$V_\alpha[\genG]\subseteq {V_\alpha}^{\uniV[\genG]}$ holds for all $\alpha\in\On$.

The induction steps for $\alpha=0$ and limit ordinals $\alpha$ are trivial. So we assume that
$V_\alpha[\genG]\subseteq {V_\alpha}^{\uniV[\genG]}$ holds and show that the same inclusion 
holds for $\alpha+1$. Suppose $a\in V_{\alpha+1}[\genG]$. Then $a=\uta^\genG$ for 
a $\poP$-name $\uta\in V_{\alpha+1}$. Since $\uta\subseteq V_\alpha$, each
$\pairof{\utb,\condp}\in\uta$ is an element of $V_\alpha$. By induction hypothesis, it 
follows that $\utb^\genG\in {V_\alpha}^{\uniV[\genG]}$. It follows that 
$\uta^\genG\subseteq{V_\alpha}^{\uniV[\genG]}$. Thus 
$a=\uta^\genG\in{V_{\alpha+1}}^{\uniV[\genG]}$. \smallskip 

``$\supseteq$'':  
Suppose that $a\in {V_\alpha}^{\uniV[\genG]}$. Note that we can choose the ``sufficiently 
large finite fragment of \ZFC'' which should hold in $V_\alpha$, \st\ this implies that 
{\ifextended\extendedcolor $(*)$\fi}
${V_\alpha}^{\uniV[\genG]}$ still satisfies a large enough fragment of \ZFC, although the 
fragment may be different from the one $V_\alpha$ satisfies. In particular we find a 
cardinal $\lambda>\cardof{\poP}$ in ${V_\alpha}^{\uniV[\genG]}$ (and hence it is also a 
cardinal in
$\uniV[\genG]$) \st\
$a\in\calH(\lambda)^{V_\alpha^{\uniV[\genG]}}
\subseteq\calH(\lambda)^{\uniV[\genG]}\subseteq {V_\alpha}^{\uniV[\genG]}$. 
{\ifextended\extendedcolor {[}\,Note that
    $\calH(\lambda)^{{V_\alpha}^{\uniV[\genG]}}
    =\setof{a}{\cardof{\trcl(a)}<\lambda}^{{V_\alpha}^{\uniV[\genG]}}
    \subseteq\setof{a}{\cardof{\trcl(a)}<\lambda}^{\uniV[\genG]}
    =\calH(\lambda)^{\uniV[\genG]}.$ {]}
\fi}

Let
$a^*\in\calH(\lambda)^{\uniV[\genG]}$ be a transitive set \st\ $a\in a^*$. Then $a^*$ 
can be coded by a subset of $\lambda$. We can find the subset of $\lambda$ 
in $\uniV[\genG]$ and this subset has a nice $\poP$-name which is an element of
${V_\alpha}^\uniV$ since $\poP\in V_\alpha$. 
This shows that $a^*\in V_\alpha[\genG]$ and hence also $a\in V_\alpha[\genG]$. 
\qedofLemmaA\qedskip\fi
}

\noindent
\prfof{\bfThmof{p-Lg-RcA-0}} 
Assume that $\kappa:=\kappa_\refl$ is tightly $\calP$-Laver generically ultrahuge for an 
iterable class $\calP$ of \pos. 

Suppose that $\varphi=\varphi(\overline{x})$ is $\Sigma_2$ formula (in $\Lin$), 
$\psi=\psi(\overline{x})$ is $\Pi_2$ formula (in $\Lin$), 
$\overline{a}\in\calH(\kappa)$, and  
$\poP\in\calP$ is \st\ 
\begin{xitemize}
\xitemA[x-Lg-RcA-a] 
  $\uniV\models\forces{\poP}{\varphi(\overline{a})\land \psi(\overline{a})}$. 
\end{xitemize}

Let $\lambda>\kappa$ be 
\st\ $\poP\in\uniV_\lambda$ and 
\begin{xitemize}
\xitemA[x-Lg-RcA-0] $V_\lambda\prec^{}_{\Sigma_n}\uniV\mbox{ for a sufficiently large }n$. 
\end{xitemize}
In particular, we may assume that we have chosen the $n$ above so that a sufficiently large 
fragment of \ZFC\ holds in 
$V_\lambda$ (``sufficiently large'' means here, in particular, in terms of 
\Lemmaof{p-Lg-RcA-0-0} and that the argument at the end of this proof is possible).  

Let $\utpoQ$ be a $\poP$-name \st\ $\forces{\poP}{\utpoQ\in\calP}$, and for
$(\uniV,\poP\ast\utpoQ)$-generic $\genH$, there are $j$, $M\subseteq\uniV[\genH]$ with 
\begin{xitemize}
  \xitemA[x-Lg-RcA-1] 
  $\Elembed{j}{\uniV}{M}{\kappa}$, 
\end{xitemize}
\begin{xitemize}
  \xitemA[x-Lg-RcA-1-0] 
  $j(\kappa)>\lambda$, 
\end{xitemize}
\begin{xitemize}
  \xitemA[x-Lg-RcA-1-1] 
  $\poP\ast\utpoQ,\ \poP,\ \genH,\ {V_{j(\lambda)}}^{\uniV[\genH]}\in M,\mbox{ and }$
\end{xitemize}
\begin{xitemize}
  \xitemA[x-Lg-RcA-1-2]
  $\cardof{\poP\ast\utpoQ}\leq j(\kappa)$. 
\end{xitemize}
By \xitemof{x-Lg-RcA-1-2}, we may assume that the underlying set 
of $\poP\ast\utpoQ$ is $j(\kappa)$ and $\poP\ast\utpoQ\in {V_{j(\lambda)}}^\uniV$. 

Let $\genG:=\genH\cap\poP$. 
Note that $\genG\in M$ by \xitemof{x-Lg-RcA-1-1} 
we have 
\begin{xitemize}
\xitemA[x-Lg-RcA-2] 
  ${V_{j(\lambda)}}^M\ubecause{=}{}{by \xitemof{x-Lg-RcA-1-1}}
  {V_{j(\lambda)}}^{\uniV[\genH]}
  \obecause{=}{1.44ex}{\hspace{20em}
    \vbox{\hbox{Since ${V_{j(\lambda)}}^M$ ($={V_{j(\lambda)}^{\uniV[\genH]}}$) satisfies a 
    sufficiently large fragment of 
    \ZFC\vspace{-0.18ex}}\hbox{by elementarity of $j$, and by 
        \Lemmaof{p-Lg-RcA-0-0}}}} 
  {V_{j(\lambda)}}^\uniV[\genH].$ 
\end{xitemize}
Thus, by \xitemof{x-Lg-RcA-1-1} and by the definability of grounds, we have 
${V_{j(\lambda)}}^\uniV\in M$ and ${V_{j(\lambda)}}^\uniV[\genG]\in M$. 
We may assume that $V_{j(\lambda)}^{\uniV}$ as a ground of $V_{j(\lambda)}^M$ satisfies a 
large enough fragment of \ZFC. 

\begin{Claim}
  \Label{cl-Lg-RcA-0}
  ${V_{j(\lambda)}}^\uniV[\genG]\models\varphi(\overline{a})\land\psi(\overline{a})$.  
\end{Claim}
\noindent
\prfofClaim
By \Lemmaof{p-Lg-RcA-0-0}, ${V_\lambda}^{\uniV}[\genG]={V_\lambda}^{\uniV[\genG]}$, and
${V_{j(\lambda)}}^{\uniV}[\genG]={V_{j(\lambda)}}^{\uniV[\genG]}$. 
By \xitemof{x-Lg-RcA-0}, both ${V_\lambda}^\uniV[\genG]$ and $V_{j(\lambda)}^\uniV[\genG]$ 
satisfy still large enough fragment of \ZFC. Thus, by \Lemmaof{p-Lg-RcA-0-0} below, 
it follows that 
\begin{xitemize}
  \xitemA[x-Lg-RcA-2-0] 
  ${V_\lambda}^\uniV[\genG]\prec_{\Sigma_1}{V_{j(\lambda)}}^\uniV[\genG]\prec_{\Sigma_1}V[\genG]$. 
\end{xitemize}
By \xitemof{x-Lg-RcA-a} and \xitemof{x-Lg-RcA-0}, we have
${V_\lambda}^\uniV[\genG]\models\varphi(\overline{a})$ and $\uniV[\genG]\models\psi(\ol{a})$. By 
\xitemof{x-Lg-RcA-2-0} and since 
$\varphi$ is $\Sigma_2$, and $\psi$ is $\Pi_2$, it follows that
${V_{j(\lambda)}}^\uniV[\genG]\models\varphi(\overline{a})\land\psi(\overline{a})$.  
\qedofClaim\qedskip

Thus we have 
\begin{xitemize}
  \xitemA[x-Lg-RcA-3] 
  $M\modelof{\mbox{there is a }
  \calP\mbox{-ground }N\mbox{ of }V_{j(\lambda)}\mbox{ with }N\models\varphi(\overline{a})\land\psi(\overline{a})}$.
\end{xitemize}

By the elementarity \xitemof{x-Lg-RcA-1}, it follows that 
\begin{xitemize}
  \xitemA[x-Lg-RcA-4] 
  $\uniV\modelof{\mbox{there is a }
  \calP\mbox{-ground }N\mbox{ of }V_{\lambda}\mbox{ with }N\models\varphi(\overline{a})\land\psi(\overline{a})}$.
\end{xitemize}

Now by \xitemof{x-Lg-RcA-0}, it follows that there is a $\calP$-ground $\uniW$ of $\uniV$ \st\\
$\uniW\models\varphi(\overline{a})\land\psi(\overline{a})$. 
  \qedofThm\qedskip
\fi}

\Thmof{p-Lg-RcA-0} has an important application (\Thmof{p-Lg-RcA-1}). For 
this theorem, we need the following facts about Recurrence Axioms.
\begin{Lemma}{\rm(Fuchino and Usuba \cite{recurrence}, see also \href{https://fuchino.ddo.jp/papers/reflection\_and\_recurrence-Janos-Festschrift-x.pdf\pdfpage{18}}{Lemma 20} in the 
    extended version of \cite{janos})}
  \Label{p-Lg-RcA-0-1} Assume that $\calP$ is an iterable class of \pos. \wassertof{1} If 
  $\calP$ contains a \po\ which adds a real (over 
  the universe), then 
  $(\calP,\calH(\kappa_\refl))_{\Sigma_1}$-\RcA\ implies $\neg\CH$.\smallskip

  \wassert{2} Suppose that $\calP$ contains a 
  \po\ which 
  forces ${\aleph_2}^\uniV$ to be equinumerous with ${\aleph_1}^\uniV$. Then 
  $(\calP,\calH(2^{\aleph_0}))_{\Sigma_1}$-\RcA\ implies $\continuum\leq\aleph_2$.
  \smallskip

  \wassert{2'} If $\calP$ contains a \pos\ which forces ${\aleph_2}^\uniV$ to be 
  equinumerous with 
  ${\aleph_1}^\uniV$, then $(\calP,\calH((\aleph_2)^+))_{\Sigma_1}$-\RcA\ does not hold. 
  \smallskip

  \wassert{3} If $(\calP,\calH(\kappa_\refl))_{\Sigma_1}$-\RcA\ holds then all 
  $\poP\in\calP$ preserve $\aleph_1$ and they are also stationary preserving.\smallskip

  \wassert{4} If $\calP$ contains a \po\ which 
  adds a real as well as a \po\ which collapses 
  ${\aleph_2}^\uniV$, then 
  $(\calP,\calH(\kappa_\refl))_{\Sigma_1}$-\RcA\ implies $\continuum=\aleph_2$.
  \smallskip

  \wassert{5} If $\calP$ contains a 
  \po\ which 
  collapses ${\aleph_1}^\uniV$, then $(\calP,\calH(2^{\aleph_0}))_{\Sigma_1}$-\RcA\ implies 
  \CH.\smallskip 

  \wassert{5'} If $\calP$ contains a \po\ which collapses ${\aleph_1}^\uniV$ then 
  $(\calP,\calH((\continuum)^+))_{\Sigma_1}$-\RcA\ does not hold. 

    \wassert{6} Suppose that all $\poP\in\calP$ preserve cardinals 
    and $\calP$ contains \pos\  
  adding at least $\kappa$ many reals for each $\kappa\in\Card$. Then
  $(\calP,\emptyset)_{\Sigma_2}$-\RcAp\ implies that $\continuum$ is very large.

  \wassert{6'} Suppose that $\calP$ is as in \assertof{6}. Then 
  $(\calP,\calH(\continuum))_{\Sigma_2}$-\RcAp\ implies that $\continuum$ is a limit 
  cardinal. Thus if $\continuum$ is regular in addition, then $\continuum$ is weakly 
  incaccessible.
  \qed
\end{Lemma}

\begin{Thm}\Label{p-Lg-RcA-1}
  Suppose that $\calP$-LgLCA for extendible holds. Then we have:\smallskip\\
  \wassert{1} Elements of $\calP$ are stationary preserving.\smallskip

  \wassert{2}
  For all classes $\calP$ of \pos\ covered by \Lemmaof{p-Lg-RcA-0-1}, $\calP$-LgLCA for extendible 
  implies that the continuum is either $\aleph_1$ or $\aleph_2$ or very large. 
\end{Thm}
\memo{Since under the constraint of \Thmof{p-Lg-RcA-1},\assertof{1}, cases of 
\Lemmaof{p-Lg-RcA-0-1} are quite exhaustive (perhaps except that the cases with classes of \pos\ 
with various $\mu$-closedness properties for $\mu>\aleph_1$ are not considered), we can see 
\Thmof{p-Lg-RcA-1},\assertof{2} as the trichotomy theorem for the size of the continuum.}
\prf \assertof{1}: By \Thmof{p-Lg-RcA-0} and \Lemmaof{p-Lg-RcA-0-1},\,\assertof{3}. 
\smallskip

\assertof{2}: By \Thmof{p-Lg-RcA-0} and the rest of \Lemmaof{p-Lg-RcA-0-1}. \qedofThm\qedskip

The proof of the following theorem is almost identical with the older proof in 
\cite{recurrence}. Nevertheless I would like to repeat the proof here since this 
proof shows how the notion of super-$C^{(\infty)}$ version of the LgLCA is incorporated into 
the whole picture. 

\begin{Thm}{\rm(Fuchino and Usuba \cite{recurrence}, \href{https://fuchino.ddo.jp/papers/recurrence-axioms-x.pdf\pdfpage{23}}{Theorem 4.10})}
  \LabelR{p-Lg-RcA-5} Suppose that $\calP$ is an iterable class of \pos\ and the 
  super-$C^{(\infty)}$-$\calP$-LgLCA for extendible holds. Then 
  $\MP(\calP,\calH(\kappa_\refl))$ holds. \iftrue\else\qed\fi
\end{Thm}
{\iftrue
\prf It is enough to show that $(\calP,\calH(\kappa_\refl))$-\RcAp\ holds. For this, 
a modification of the proof of \Thmof{p-Lg-RcA-0} works. 

Suppose that $\kappa:=\kappa_\refl$ is tightly 
super-$C^{(\infty)}$-$\calP$-Lg extendible,
$\poP\in\calP$, and $\forces{\poP}{\varphi(\overline{a})}$ for 
an $\Lin$-formula $\varphi$ and $\overline{a}\in\calH(\kappa)$. We want to show that
$\varphi(\overline{a})$ holds in some $\calP$-ground of $\uniV$.

Let $n$ be a sufficiently large natural number $\GE1$ \st\ the following arguments go through. In 
particular, we assume that ${V_\alpha}^\uniV\prec_{\Sigma_n}\uniV$ implies that
``$\varphi(\overline{x})$'' and ``$\forces{\cdot}{\varphi(\overline{x})}$ are 
absolute between ${V_\alpha}^\uniV$ and $\uniV$, 
and ${V_\alpha}^\uniV\prec_{\Sigma_n}\uniV$ also implies that a sufficiently large fragment 
of \ZFC\ holds in $V_\alpha$. 

Let $\utpoQ$ be a $\poP$-name \st\ $\forces{\poP}{\utpoQ\in\calP}$ and, for
$(\uniV,\poP\ast\utpoQ)$-generic $\genH$, there are a $\lambda>\kappa$ 
with 
\begin{xitemize}
\xitem[x-Lg-RcA-9] 
  $V_\lambda\prec_{\Sigma_n}\uniV$, 
\end{xitemize}
and   $j$, $M\subseteq \uniV[\genH]$ \st
\begin{xitemize}
\xitem[x-Lg-RcA-9-a] 
  \assertof{a}\ \,$\Elembed{j}{\uniV}{M}{\kappa}$,\quad
  \assertof{b}\ \,$j(\kappa)>\lambda$,\quad\assertof{c}\ \,$\poP$, $\genH$,
  ${V_{j(\lambda)}}^{\uniV[\genH]}\in M$,\\
  \assertof{d}\ \,$\cardof{RO(\poP\ast\utpoQ)}\leq j(\kappa)$  ($<j(\lambda)$), and 
  \quad\assertof{e}\ \,${V_{j(\lambda)}}^{\uniV[\genH]}\prec_{\Sigma_n}\uniV[\genH]$. 
\end{xitemize}

By \xitemof{x-Lg-RcA-9-a},\,\assertof{c} and \assertof{d}, we may assume that
$\poP\ast\utpoQ\in {V_{j(\lambda)}}^\uniV$ by replacing $\poP\ast\utpoQ$ by an 
appropriate isomorphic \po\ (and replacing $\genH$ by  
corresponding filter), 

By the choice of $n$, we 
have $V_\lambda\models\forces{\poP}{\varphi(\overline{a})}$.
$j({V_\lambda}^\uniV)={V_{j(\lambda)}}^M\prec_{\Sigma_n}M$ by elementarity of $j$,   
and 
\begin{xitemize}
\xitem[x-Lg-RcA-9-0] 
  ${V_{j(\lambda)}}^M={V_{j(\lambda)}}^{\uniV[\genH]}$
\end{xitemize}
by the closedness of $M$. 
Since $V_\lambda\prec_{\Sigma_n}\uniV$, we have
$V_\lambda[\genH]\prec_{\Sigma_{n_0}}\uniV[\genH]$ for a still large 
enough $n_0\leq n$. Since
${V_{j(\lambda)}}^{\uniV[\genH]}\prec_{\Sigma_n}\uniV[\genH]$, it  
follows 
that ${V_\lambda}^{\uniV[\genH]}\prec_{\Sigma_{n_0}}{V_{j(\lambda)}}^{\uniV[\genH]}$. 
Thus 
\begin{xitemize}
\xitem[x-Lg-RcA-10] 
  ${V_\lambda}^\uniV\prec_{\Sigma_{n_1}}{V_{j(\lambda)}}^\uniV$
\end{xitemize}
for a still large enough
$n_1\leq n_0$. 

In particular, we have
${V_{j(\lambda)}}^\uniV\models\forces{\poP}{\varphi(\overline{a})}$,  and hence
$V_{j(\lambda)}[\genG]\models\varphi(\overline{a})$ where $\genG$ is the $\poP$-part of
$\genH$. Note that by \xitemof{x-Lg-RcA-9} and \xitemof{x-Lg-RcA-10}, $V_{j(\lambda)}$ 
satisfies a sufficiently large fragment of \ZFC. 

Thus we have
$V_{j(\lambda)}[\genH]\modelof{\,\mbox{there is a }
  \calP\mbox{-ground satisfying }\varphi(\overline{a})}$, and hence 
\begin{xitemize}
\item[] ${V_{j(\lambda)}}^{\uniV[\genH]}\modelof{\,\mbox{there is a }
  \calP\mbox{-ground satisfying }\varphi(\overline{a})}$
\end{xitemize}
by \LemmaAof{p-Lg-RcA-0-0}. By 
\xitemof{x-Lg-RcA-9-0} and elementarity, it follows that
\begin{xitemize}
\item[] ${V_{\lambda}}\modelof{\,\mbox{there is a }
  \calP\mbox{-ground satisfying }\varphi(\overline{a})}$.
\end{xitemize}
Finally, this implies ${\uniV}\modelof{\,\mbox{there is a }
  \calP\mbox{-ground satisfying }\varphi(\overline{a})}$ by \xitemof{x-Lg-RcA-9}.\\
\qedofThm\qedskip\fi}

For an ordinal $\alpha$, let
$\alpha^{(+)}:=\sup(\setof{\cardof{\beta}^+}{\beta<\alpha})$. Note that
$\alpha^{(+)}=\alpha$ if $\alpha$ is a cardinal. Otherwise, we have
$\alpha^{(+)}=\cardof{\alpha}^+$. 

\begin{Thm}{\rm(\href{https://fuchino.ddo.jp/papers/generic-absoluteness-revisited-x.pdf\pdfpage{39}}{Theorem 5.7} in Fuchino, Gappo and Parente \cite{FGP} restated under LgLCA 
    for extendible)}\Label{p-genabs-Laver-1} For an iterable class $\calP$ of \pos, assume 
  that   
  $\calP$-LgLCA for extendible holds.   
  Then, for $\kappa:=\kappa_\refl$
  \begin{xitemize}
  \xitem[x-genabs-Laver-a] 
    for any $\poP\in\calP$ \st\  $\forces{\poP}{\BFA_{\LT\kappa}(\calP)}$, 
    $\calH(\mu^+)^\uniV\prec_{\Sigma_2}\calH(\mu^+)^{\uniV[\genG]}$\quad
    holds 
    for all $\mu<\kappa$ and for $(\uniV,\poP)$-generic $\genG$.\footnotemark
  \end{xitemize}
  Thus, we have
  $\calH(\kappa)^\uniV\prec_{\Sigma_2}\calH((\kappa^{(+)})^{\uniV[\genG]})^{\uniV[\genG]}$
  for $\genG$ as above. 
  \ifextended\else\qed\fi
\end{Thm}
\footnotetext{$\mu^+$ is $\calH(\mu^+)^{\uniV[\genG]}$ is actually $(\mu^+)^{\uniV[\genG]}$.}
{\ifextended\extendedcolor
\prf Suppose that $\forces{\poP}{\calH(\mu^+)\models\varphi(\ol{a})}$ for $\poP\in\calP$ 
with $\forces{\poP}{\BFA_{\LT\kappa}(\calP)}$, $\mu<\kappa$, $\Sigma_2$-formula $\varphi$ and
for $\ol{a}\in\calH(\mu^+)$. Let $\genG$ be a $(\uniV,\poP)$-generic filter. Then we have
\begin{xitemize}
\xitemA[x-genabs-Laver-a-0] 
  $\uniV[\genG]\modelof{\BFA_{\LT\kappa}(\calP)\land\calH(\mu^+)\models\varphi(\ol{a})}$. 
\end{xitemize}
Let $\varphi=\exists y\psi(\ol{x},y)$ where $\psi$ is a $\Pi_1$-formula in $\Lin$. Let
$b\in\calH((\mu^+)^{\uniV[\genG]})^{\uniV[\genG]}$. be \st\
$\calH((\mu^+)^{\uniV[\genG]})^{\uniV[\genG]}\models\psi(\ol{a},b)$. 

Let $\utpoQ$ be a $\poP$-name with $\forces{\poP}{\utpoQ\in\calP}$ \st, for
$(\uniV,\poP\ast\utpoQ)$-generic $\genH$ with 
\begin{xitemize}
\xitemA[x-genabs-Laver-a-1] 
  $\genG\subseteq\genH$ (under the 
  identification $\poP\circleq\poP\ast\utpoQ$), 
\end{xitemize}
there are $j$, $M\subseteq\uniV[\genH]$ \st\
$\Elembed{j}{\uniV}{M}{\kappa}$, 
\begin{xitemize}
\xitemA[x-genabs-Laver-0] $\cardof{\poP\ast\utpoQ}\leq j(\kappa)$\qquad (by tightness), 
\xitemA[x-genabs-Laver-1] $\poP$, $\poP\ast\utpoQ$, $\genH\in M$ and
\xitemA[x-genabs-Laver-2] $j\imageof{j(\kappa)}\in M$. 
\end{xitemize}
\memo{[25.05.19(Mo09:40(JST))] \xitemof{x-genabs-Laver-2} can be replaced by ${V_{j(\lambda)}}^{\uniV[\genH]}\in M$. 
  Thus the assumtion of the theorem can be replaced by $\calP$-LgLCA for extendible.}

By \xitemof{x-genabs-Laver-a-0}, \xitemof{x-genabs-Laver-a-1} and Bagaria's Absoluteness 
Theorem (applied to $V[\genG]$), we have $\uniV[\genH]\modelof{\psi(\ol{a},b)}$  and hence
$\uniV[\genH]\modelof{\calH(\mu^+)\models\psi(\ol{a},b)}$.

By \xitemof{x-genabs-Laver-0}, and \xitemof{x-genabs-Laver-2}, 
there is a $\poP$-name of $b$ in $M$. By \xitemof{x-genabs-Laver-1}, it follows that
$b\in M$. By similar argument, we have $\calH((\mu^+)^{\uniV[\genH]})^{\uniV[\genH]}\subseteq M$ and 
hence $\calH((\mu^+)^{\uniV[\genH]})^{\uniV[\genH]}=\calH((\mu^+)^M)^M\in M$. Thus 
$M\modelof{\calH(\mu^+)\models\psi(\ol{a},b)}$.

By elementarity, it follows that
$\uniV\modelof{(\exists\ul{b}\in\calH(\mu^+))\ \calH(\mu^+)\models\psi(\ol{a},\ul{b})}$, and  
hence $\uniV\modelof{\calH(\mu^+)\models\varphi(\ol{a})}$ as desired. 
\smallskip

Suppose now that $\poP$, $\mu$, $\varphi$,
$\ol{a}$ are as above and assume that 
$\uniV\modelof{\calH(\mu^+)\models\varphi(\ol{a})}$ holds. For a $\Pi_1$-formula $\psi$ as 
above, let $b\in\calH(\mu^+)^\uniV$ be \st\ $V\modelof{\calH(\mu^+)\models\psi(\ol{a},b)}$.

Since $\uniV\models\BFA_{\LT\kappa}(\calP)$ by assumption, 
it follows that $\uniV[\genG]\models\psi(\ol{a},b)$ by Bagaria's Absoluteness, and hence
$\uniV[\genG]\modelof{\calH(\mu^+)\models\varphi(\ol{a})}$.  

The last assertion of the theorem follows from this.
\qedofThm\fi}

\section{Consistency of LgLCAs and 
        super-\texorpdfstring{$C^{(n)}$}{C(∞)}-LgLCAs for extendibility}
\LabelR{Lg-ext}
We examine first that extendible and super-$C^{(\infty)}$-extendible cardinals are endowed with 
Laver functions. 

A function $\mapping{f}{\kappa}{V_\kappa}$ for an extendible cardinal $\kappa$ is said to 
be a {\It Laver fuction for extendibility} (of $\kappa$) if for any $\lambda>\kappa$ and any set $a$, 
there are $j$, $M\subseteq\uniV$ \st\ $\Elembed{j}{\uniV}{M}{\kappa}$, $j(\kappa)>\lambda$,
$V_{j(\kappa)}\in M$, and $j(f)(\kappa)=a$.

Suppose that $\mu$ is an inaccessible cardinal and
$V_\mu\modelof{\kappa\xmbox{ is a super-}C^{(\infty)}\xmbox{-ex\-tend\-ible}}$. Then 
$\mapping{f}{\kappa}{V_\kappa}$ is said to be a {\It super-$C^{(\infty)}$-extendible Laver 
function for $\kappa$ in $V_\mu$} if, for any $x\in V_\mu$, $n\in\omega$ and for any $\kappa<\lambda<\mu$ with
$\lambda\in C^{(n)}$, there are $j$,
$M\subseteq V_\mu$ \st\ $\Elembed{j}{V_\mu}{M}{\kappa}$, $j(\kappa)>\lambda$,
$j(\lambda)\in C^{(n)}$, $V_{j(\lambda)}\in M$ and $j(f)(\kappa)=x$. 

\begin{Lemma}\Label{p-Lg-ext-0} \wassertof{1} Suppose that $\kappa$ is 
  extendible. Then there is a Laver function
  $\mapping{f}{\kappa}{V_\kappa}$ for extendibility. \smallskip


  \wassert{2} 
  Suppose that, for an inaccessible cardinals $\mu$ and a cardinal $\kappa<\mu$,  
  we have $V_{\mu}\modelof{\kappa\mbox{ is strongly super-}C^{(\infty)}\mbox{-extendible}}$. 
  Then there is a super-$C^{(\infty)}$-extendible Laver function $\mapping{f}{\kappa}{V_{\kappa}}$ for 
   $\kappa$ in $V_{\mu}$. 
\end{Lemma}
\memo{Scan\_2024-11-23--10.44 extendible - annotated p.82-83}
\prf \assertof{1}: is known previously (see e.g.\ Corraza \cite{corraza}). The proof of 
\assertof{2} below can be easily modified and slimmed down to a proof of \assertof{1}.\smallskip

\assertof{2}: 
Assume, towards a contradiction, that there is no super-$C^{(\infty)}$-extendible Laver function
$\mapping{f}{\kappa}{V_\kappa}$ for $\kappa$ in $V_\mu$. 

Let $C:=\setof{\alpha<\mu}{V_\alpha\prec V_\mu}$. 

For $n\in\omega$, let $\varphi_n(\ul{f}, \ul{x})$ be the formula 
\begin{xitemize}
\xitemx[] 
  $\forall\ul{\delta}\,\forall\ul{\delta}'\,
  \forall\ul{j}\,\forall\alpha\ \varphi^*_n(\ul{f},\ul{\alpha},\ul{x},\ul{\delta},\ul{\delta'},\ul{j})$ 
\end{xitemize}
  \memo{see \href{https://fuchino.ddo.jp/notes/math-notes-20.pdf\pdfpage{58}}{math-notes-20 L.13.18 (3)}}
  where $\varphi^*_n(\ul{f},\ul{\alpha},\ul{x},\ul{\delta},\ul{\delta'},\ul{j})$ is the
  $\Sigma_n$-formula:
{\small\begin{itemize}
\xitemx[] ${((\,{}}\mapping{\ul{f}}{\ul{\alpha}}{V_{\ul{\alpha}}}
  \ \land\ \ul{\alpha}<\ul{\delta}\ \land\ \ 
  \obecause{V_{\ul{\delta}}\prec_{\Sigma_n}\uniV}{}{$\Pi_n$-formula}
  \ \land\ \obecause{V_{\ul{\delta}'}\prec_{\Sigma_n}\uniV}{}{$\Pi_n$-formula}
   \land\ \ \Elembed{\ul{j}}{V_{\ul{\delta}}}{V_{\ul{\delta}'}}{\ul{\alpha}}\ \land\ 
   \ul{j}(\ul{\alpha})>\ul{\delta}
   \,)$\\
   \mbox{}\hfill
   $\rightarrow\ \ \ul{j}(\ul{f})(\ul{\alpha})\not\equiv\ul{x}\,).$
\end{itemize}}
If $V_\mu\models\exists\ul{x}\varphi_n(f,\ul{x})$ for some $n\in\omega$ and $f\in V_\mu$, then
let $x_{n,f}\in V_\mu$ be a witnesses for 
$\ul{x}$ in $\varphi_n(f)$, and let $\mu_{n,f}:=\rank(x_{n,f})$. 
If $V_\mu\models\exists\ul{x}\varphi_n(f,\ul{x})$ does not hold, we let $x_{n,f}:=\emptyset$ 
and $\mu_{n,f}:=0$.

$x_{n,f}$ might not be determined 
uniquely. However, we can choose $x_{n,f}$ 
\st\ $\mu_{n,f}$ is minimal and, in this way, $\mu_{n,f}$ is determined uniquely to 
each $n$ and $f$.

By assumption, we have 
\begin{xitemize}
\xitem[x-Lg-ext-0] 
  $V_\mu\models\exists \ul{x}\,\varphi_n(f,\ul{x})$ for some $n\in\omega$ holds, for all $\mapping{f}{\kappa}{V_\kappa}$.
\end{xitemize}

Since $\kappa$ is strongly super-$C^{(\infty)}$-extendible in $V_\mu$, we can  
find $\lambda<\mu$ and $j^*$, $M^*\subseteq V_\mu$  by 
\Lemmaof{p-ext-1-0} \st\  
\begin{xitemize}
\xitem[x-Lg-ext-0-a-0] 
  $\lambda\geq\lambda_0
  :=\max\setof{\mu_{n,f}}{n\in\omega,\, 
  \mapping{f}{\alpha}{V_\alpha}\xmbox{ for an inaccessible }\alpha\leq\kappa}$,  
\xitem[x-Lg-ext-0-0] $\lambda\in C$, 
\xitem[x-Lg-ext-0-1] $\Elembed{j^*}{V_\mu}{M^*}{\kappa}$ with\ \ 
  \assertof{a}\ \,$j^*(\kappa)>\lambda$,\ \ \assertof{b}\ \,$j^*(\lambda)\in C$, and\ \ 
  \assertof{c}\ \,$V_{j^*(\lambda)}\in M^*$.
\end{xitemize}

Let
\begin{xitemize}
\xitemx[] 
  ${\darkred A}:=\setof{\alpha<\kappa}{{}
  \begin{array}[t]{@{}l}
    \alpha\mbox{ is inaccessible, and for all }\mapping{f}{\alpha}{V_\alpha}\mbox{, there is }
    n\in\omega\\
    \mbox{\st\ }
    V_\mu\models\exists{\ul{x}}\varphi_n(f,\ul{x})\ }.
  \end{array}$
\end{xitemize}

\begin{Claim}\Label{cl-Lg-ext-a}
$M\models \kappa\in j^*(A)$.   
\end{Claim}
\prfofClaim For each $\mapping{f}{\kappa}{V_\kappa}$, there is $n\in\omega$ \st\
$V_\mu\models\exists\ul{x}\varphi_n(f,\ul{x})$ (see \xitemof{x-Lg-ext-0}). It follows that
$V_{j^*(\lambda)}\models\exists\ul{x}\varphi_n(f,\ul{x})$ by 
\xitemof{x-Lg-ext-0-1},\,\assertof{b}. By \xitemof{x-Lg-ext-0-0} and elementarity,
$V_{j^*(\lambda)}\prec_{\Sigma_m} M$ holds for all $m\in\omega$. Thus 
$M\models\exists\ul{x}\varphi_n(f,\ul{x})$. This implies
$M\models\kappa\in j^*(A)$. 
\qedofClaim\qedskip

Let $\mapping{f^*}{\kappa}{V_\kappa}$ be defined by
  
\begin{xitemize}
\item[] 
  $f^*(\alpha):=\left\{
  \begin{array}{@{}ll}
    x_{n, f^*\restr\alpha}, &\mbox{if }\alpha\in A\mbox{, }
    \mapping{f^*\restr\alpha}{\alpha}{V_\alpha}\\[-1pt]
    &\mbox{and }n\mbox{ is minimal with }x_{n,f^*}\not=\emptyset\,;\\[\jot]
    \emptyset, &\mbox{otherwise}.
  \end{array}
  \right.$
\end{xitemize}

Let $x^*:=j^*(f^*)(\kappa)$. By definition of $f^*$, by 
\Claimof{cl-Lg-ext-a}, and elementarity, we have 
 $x^*\not=\emptyset$ and $x^*$ 
witnesses $\exists\ul{x}\varphi_n(f^*,\ul{x})$ for some $n\in\omega$. 
($\mu_{n,f^*}$ is just as chosen before since it is uniquely determined. 
$x^*$ may be different from $x_{n,f^*}$ but this does not matter.)\smallskip

In particular, $x^*\smashx{\ubecause{\not=}{}{\qquad\qquad\qquad by the property of $j(f^*)$ 
    inherited from the 
  definition of $f^*$}}j(f^*)(\kappa)\smash{\obecause{=}{}{by definition of $x^*$ }}x^*$. 
This is a contradiction. 
\qedofLemma\qedskip\memo{add Scan\_2024-11-23--10.44 extendible - annotated pp.51-52}


\begin{Thm}\Label{p-Lg-ext-1} Suppose that $\calP$ is $\omega_1$ preserving transfinitely 
  iterable class of \pos\ \st\ either $\calP$ contains a collapsing of $\omega_2$ or adding 
  a new reals (or both). \smallskip

  \wassert{1} If $\kappa$ is extendible, and $\calP$ is $\Sigma_2$-definable, then
  there is a \po\ $\poP_\kappa\in\calP$ \st\ 
  $\forces{\poP_\kappa}{\kappa=\kappa_\refl\xmbox{ and }\calP\xmbox{-LgLCA for 
      extendible holds}}$.\smallskip

  \wassertof{2} Suppose that $\kappa$ is super-$C^{(n^*)}$-extendible for $n\in\natnums$ 
  and $n^*=\max\ssetof{n,2}$. If super-$C^{(n^*)}$-extendible Laver function for $\kappa$ 
  exists, and $\calP$ is 
  $\Sigma_{n+1}$-definable, then there is a \po\ $\poP_\kappa\in\calP$ \st\ 
  $\forces{\poP_\kappa}{\kappa=\kappa_\refl
    \xmbox{ and super-}C^{(n)}\xmbox{-}\calP\xmbox{-LgLCA for extendible holds}}$.\smallskip

  \wassertof{3} Suppose
  $V_\mu\modelof{\kappa\xmbox{ is strongly super-}C^{(\infty)}\xmbox{ extendible}}$ for 
  an inaccessible cardinal $\mu$ and $\kappa<\mu$. Then there is
  $\poP_\kappa\in{\calP}^{V_\mu}$ \st\ 
  $V_\mu\models\forces{\poP_\kappa}{\kappa=\kappa_\refl
    \xmbox{ and super-}C^{(\infty)}\xmbox{-}\calP\xmbox{-LgLCA for extendible holds}}$.
\end{Thm}
Note that many natural classes of \pos\ including the classes of all ccc \pos,  
all $\sigma$-closed \pos, all proper \pos,\ all semi-proper \pos, etc., 
are $\Sigma_2$-definable, and they also satisfy all the conditions stated at the beginning of the theorem. \qedskip 

\noindent
\prfof{\bfThmof{p-Lg-ext-1}} \assertof{1}: We show the assertion for the 
case that $\calP$ is the class of all proper \pos. The proof for the general case can be 
done by replacing the CS-iteration in the following proof by the type of iteration for which the class $\calP$ is 
transfinitely iterable. Let $f$ be a Laver function for 
extendiblity of $\kappa$ ($f$ exists by \Lemmaof{p-Lg-ext-0},\,\assertof{1}).

Let $\seqof{\poP_\alpha,\utpoQ_\beta}{\alpha\leq\kappa, \beta<\kappa}$ be an CS-iteration 
of elements of $\calP$ \st\
\begin{itemize}
\item[] $\utpoQ_\beta:={}${$\left\{\,
  \begin{array}{@{}ll}
    f(\beta), &\mbox{if }f(\beta)\mbox{ is 
      a }\poP_\beta\mbox{-name}\\ &\mbox{and }\forces{\poP_\beta}{f(\beta)\in\calP};\\[\jot] 
    \poP_\beta\mbox{-name of the trivial forcing}, &\mbox{otherwise}.
  \end{array}\right.$}
\end{itemize}\smallskip

We show that $\forces{\poP_\kappa}{\kappa\mbox{ is tightly }\calP\mbox{-Lg extendibile}}$.
Note that once this has been established, it follows that
$\uniV[\genG_\kappa]\models\kappa=\kappa_\refl$ 
by Section 5 of \cite{sfetal-II} (see also \cite{future}, 
\protect\href{https://fuchino.ddo.jp/papers/RIMS2022-RA-MP-x.pdf\pdfpage{13}}{Theorem 3.5}), 
and hence $\uniV[\genG_\kappa]\models\calP\mbox{-LgLCA for extendible}$. 
\smallskip


Let $\genG_\kappa$ be a $(\uniV,\poP_\kappa)$-generic filter. 
  In $\uniV[\genG_\kappa]$, suppose that $\poP\in\calP$ and let $\utpoP$ be 
  a $\poP_\kappa$-name for $\poP$. 
  
Suppose that $\lambda>\kappa$. 
Let $\lambda^*>\lambda$ be a cardinal \st\ $\calH(\lambda^*)=V_{\lambda^*}$.
Note that we have $V_{\lambda^*}\prec_{\Sigma_1}\uniV$.

Let $\Elembed{j}{\uniV}{M}{\kappa}$ be \st\ \quad\ixitem[x-1*] 
$j(\kappa)>\lambda^*$,\quad
\ixitemc[x-2*] $V_{j(\lambda^*)}\in M$, and\quad\\
\ixitem[x-3*] $j(f)(\kappa)=\utpoP$. 
The last condition is possible since $f$ is a Laver function for the extendible $\kappa$.
By \xitemof{x-2*}, we have
$\calH(j(\lambda^*))^M={V_{j(\lambda^*)}}^M={V_{j(\lambda^*)}}^\uniV=\calH(j(\lambda^*))^\uniV$,
and hence \ixitem[x-Lg-ext-1-0] $V_{j(\lambda^*)}\prec_{\Sigma_1}M$ and
$V_{j(\lambda^*)}\prec_{\Sigma_1}\uniV$. 

In $M$, there is a $\poP_\kappa\ast\utpoP$-name $\utpoQ$ \st\
\begin{xitemize}
\xitemx[] 
  $M\models\forces{\poP_\kappa\ast\utpoP}{{}
  \begin{array}[t]{@{}l}
    \utpoQ\in\calP\mbox{ and }\utpoQ\xmbox{ is
      the direct limit of CS-iteration of small}\\
    \mbox{\pos\ in }\calP
    \mbox{ of length }j(\kappa)},
  \end{array}$\\
  \phantom{$M\models{}$}and $\poP_\kappa\ast\utpoP\ast\utpoQ\sim j(\poP_\kappa)$.
\end{xitemize}

By 
\xitemof{x-Lg-ext-1-0}, and since ``$\underline{\poP}\in\calP$'' is $\Sigma_2$, the same 
statement holds in $\uniV$ with the same $\utpoQ$.\smallskip 
\memo{see Scan\_2024-11-23--10.44 extendible - annotated p.49}

We have $j(\poP_\kappa)/\genG_\kappa\ \sim\ \poP\ast\utpoQ$ where 
we identify $\utpoQ$ with a corresponding $\poP$-name.

Let $\genH$ be $(\uniV,j(\poP_\kappa))$-generic filter with $\genG_\kappa\subseteq\genH$. 
The lifting $\Elembed{\tilde{j}}{\uniV[\genG_\kappa]}{M[\genH]}{\kappa}$;
$\utilde{a}[\genG_\kappa]\mapsto j(\utilde{a})[\genH]$ witnesses that $\kappa$
  is tightly $\calP$-Laver generic extendible in $\uniV[\genG_\kappa]$.
  $\cardof{RO(j(\poP_\kappa/\genG_\kappa))}\leq j(\kappa)$ follows from
  $\LT j(\kappa)$-c.c.\ of $j(\poP_\kappa)/\genG_\kappa$,
  $M\models j(\kappa)^{\LT j(\kappa)}=j(\kappa)$, and \xitemof{x-2*}. So it is enough to  
  show: 

\begin{Claim}\Label{cl-Lg-ext-0}
  ${V_\alpha}^{\uniV[\genH]}\in M[\genH]$ for all $\alpha\leq j(\lambda)$.
\end{Claim}
\prfofClaim By induction on $\alpha\leq j(\lambda)$. The successor step 
  from $\alpha<j(\lambda)$ to $\alpha+1$ 
  can be 
  proved by showing that $\poP_\kappa$-names of subsets of ${V_\alpha}^{\uniV[\genH]}$ can be 
chosen as elements of $M$. This is the case because of 
\xitemof{x-2*} and $\kappa$-c.c.\ of $\poP_\kappa$.\qedofClaim\qedskip

{\iftrue
\assertof{2}: Let $f$ be a Laver function for 
super-$C^{(n)}$
extendible $\kappa$. 

Let $\seqof{\poP_\alpha,\utpoQ_\beta}{\alpha\leq\kappa, \beta<\kappa}$ be an iteration 
  of elements of $\calP$ with the support appropriate for $\calP$ \st\
  \begin{itemize}
  \item[] $\utpoQ_\beta:={}${$\left\{\,
    \begin{array}{@{}ll}
      f(\beta), &\mbox{if }f(\beta)\mbox{ is 
        a }\poP_\beta\mbox{-name}\\ &\mbox{and }\forces{\poP_\beta}{f(\beta)\in\calP};\\[\jot] 
      \poP_\beta\mbox{-name of the trivial forcing}, &\mbox{otherwise}.
    \end{array}\right.$}
  \end{itemize}\smallskip

We show that $\forces{\poP_\kappa}{\kappa\mbox{ is tightly }\calP\mbox{-Lg extendibile}}$.
As in the proof of \assertof{1}, this is enough to show.

Let $\genG_\kappa$ be a $(\uniV,\poP_\kappa)$-generic filter. 
  In $\uniV[\genG_\kappa]$, suppose that $\poP\in\calP$ and let $\utpoP$ be 
  a $\poP_\kappa$-name for $\poP$. 
  
Suppose that $\lambda>\kappa$. Then, there is 
$\lambda^*>\lambda$, $\lambda^*\in C^{(n^*)}$ with $\Elembed{j}{\uniV}{M}{\kappa}$ \st\ \quad\ixitem[x-1*+] 
$j(\kappa)>\lambda^*$,  \quad\ixitemc[x-1*+-0] $V_{j(\lambda^*)}\prec_{\Sigma_{n^*}}\uniV$,\\
\ixitemc[x-2*+] $V_{j(\lambda^*)}\in M$, and\quad
\ixitem[x-3*+] $j(f)(\kappa)=\utpoP$. 
The last condition is possible since $f$ is a Laver function for extendibility of $\kappa$.


In $M$, there is a $\poP_\kappa\ast\utpoP$-name $\utpoQ$ \st\
$\forces{\poP_\kappa\ast\utpoP}{\utpoQ\in\calP\mbox{ and }\utpoQ\xmbox{ is
    the direct limit of the iteration specified for }\calP\mbox{ of small \pos\ in }\calP
  \mbox{ of length }j(\kappa)}$, and $\poP_\kappa\ast\utpoP\ast\utpoQ\sim j(\poP_\kappa)$. 
By 
\xitemof{x-1*+-0}, 
\xitemof{x-2*+}, and since ``$\underline{\poP}\in\calP$'' is $\Sigma_{n+1}$, the same situation holds in $\uniV$.\smallskip
\memo{see Scan\_2024-11-23--10.44 extendible - annotated p.49}

We have $j(\poP_\kappa)/\genG_\kappa\ \sim\ \poP\ast\utpoQ$ where 
we identify $\utpoQ$ with a corresponding $\poP$-name.

Let $\genH$ be $(\uniV,j(\poP_\kappa))$-generic filter with $\genG_\kappa\subseteq\genH$. 
Then the lifting $\Elembed{\tilde{j}}{\uniV[\genG_\kappa]}{M[\genH]}{\kappa}$;
$\utilde{a}[\genG_\kappa]\mapsto j(\utilde{a})[\genH]$ witnesses that $\kappa$
is tightly super-$C^{(n)}$-$\calP$-Laver generic extendible in $\uniV[\genG_\kappa]$:  
${V_{j(\lambda)}}^{\uniV[\genH]}\in M[\genH]$ holds by \Claimof{cl-Lg-ext-0}.
\smallskip

\assertof{3}: Similarly to \assertof{2}. \fi}
\qedofThm
\qedskip

What we obtained so far can be put together with the results we are going to mention 
in the next two sections into the following diagram: \bigskip

\noindent
\mbox{}\hspace*{0em}
\includegraphics[width=36em]{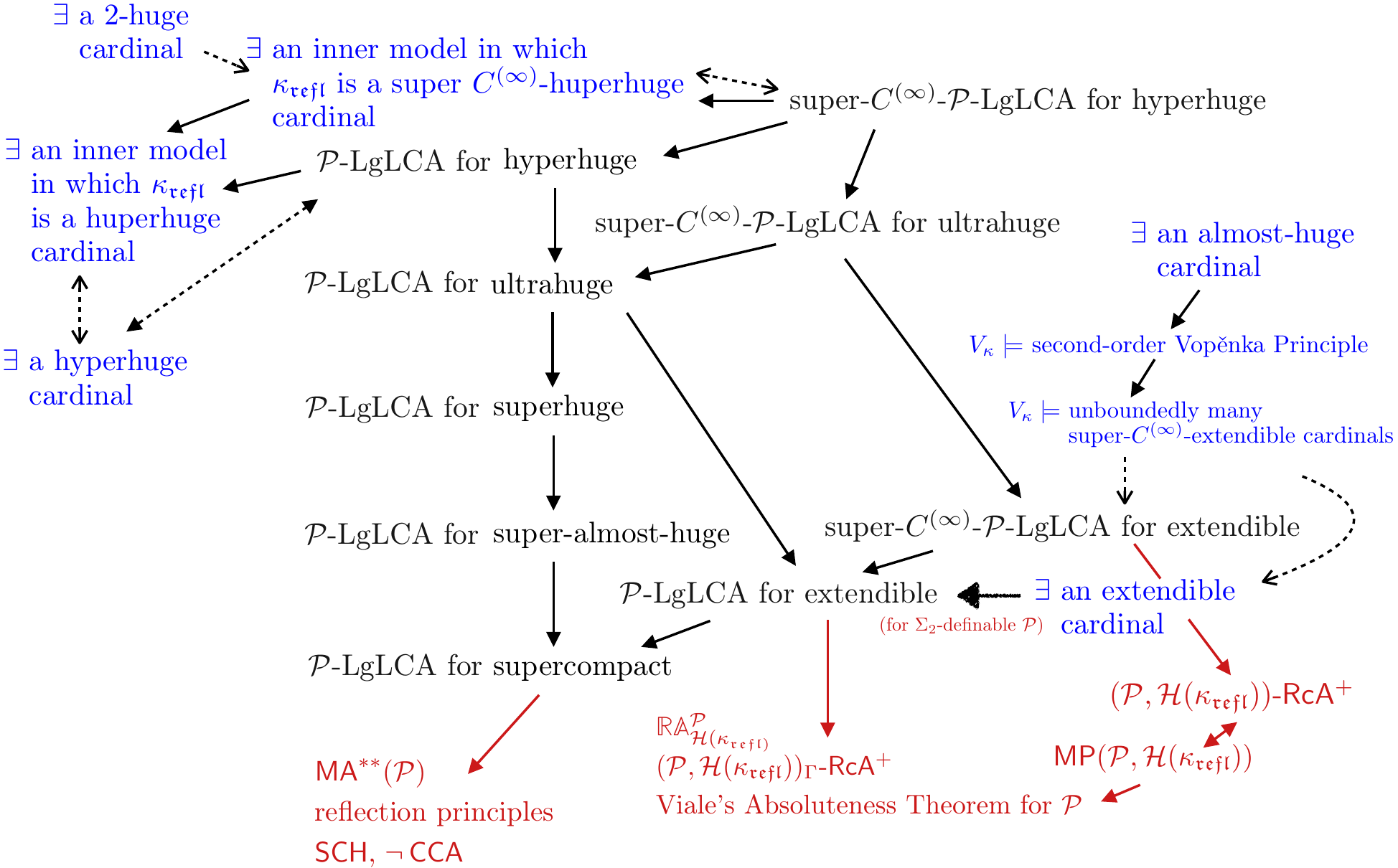}\hfill\mbox{}
\vspace{-30ex}

\noindent
\mbox{}\scalebox{0.6}{\cite{recurrence},\,\href{https://fuchino.ddo.jp/papers/recurrence-axioms-x.pdf\pdfpage{24}}{Section 5}}
\vspace{7ex}

\mbox{}\hfill{}\scalebox{0.65}{\color{cyan}\Thmof{p-Lg-RcA-5}}
\vspace{8ex}

\mbox{}\hfill{}\scalebox{0.6}{\cite{FGP},\,\href{https://fuchino.ddo.jp/papers/generic-absoluteness-revisited-x.pdf\pdfpage{35}}{Theorem 4.1}}\quad\quad\mbox{}
\vspace{5ex}

\phantomsection
\Label{the-diagram}
\smallskip

\noindent
\scalebox{0.6}{\parbox{\textwidth}{\raisebox{-0.0ex}{\includegraphics[width=3.6em]{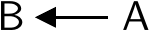}}\ : 
  ``\,the axiom A implies the axiom B''\\
  \raisebox{-0.0ex}{\includegraphics[width=3.6em]{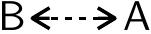}}\ : 
  ``the axioms A and B are equi-consistent.''\\
  \raisebox{-0.0ex}{\includegraphics[width=3.6em]{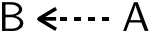}}\ : 
  ``the consistency of A implies the consistency of B \rlap{but not the other way around.''}\\
  \raisebox{-0.0ex}{\includegraphics[width=3.7em]{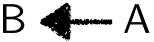}}\ : 
  ``the consistency of A implies the consistency of B \rlap{but the equi-consistency is not 
    (yet?) known.''
}}}\medskip

{\ifextended\extendedcolor
\section{Simultaneous and diagonal reflections, \\and total failure of square principles\\ under LgLCAs}\Label{square}
In this section, we return to the narration that LgLCAs are strong form of reflection 
principles. Since square principles imply the existence of structures with strong 
non-reflection properties, a natural guess is that LgLCAs imply the total failure of square 
principles, at least above a certain cardinal. \Corof{p-square-1} below confirms this intuition. 

We consider the following variant of the Reflection Property \RP\ in Jech 
\cite{millennium-book}:
\begin{xitemize}
\item[(\RP$^*$): ]  For any uncountable $\lambda$ and 
  $S_\xi\subseteq[\lambda]^{\aleph_0}$, $\xi<\eta$ for some $\eta<\kappa_\refl$ 
  \st\ $S_\xi$ is stationary in $[\lambda]^{\aleph_0}$ for all $\xi<\eta$, there is
  $X\in[\lambda]^{\LT\kappa_\refl}$ \st\ $S_\xi\cap[X]^{\aleph_0}$ is stationary in
  $[X]^{\aleph_0}$ for all $\xi<\eta$.
\end{xitemize}
  
\begin{Thm}\Label{p-square-0}
  \wassertof{1} Suppose that $\calP$ is proper. Then $\calP$-LgLCA for supercompact implies
  $\RP^*$. \smallskip 

  \wassert{2} Suppose that $\calP$ is stationary preserving  (i.e. preserving stationarity 
  of stationary subsets of $\omega_1$),   
  and $\calP$ contains all $\sigma$-closed \pos. Then $\calP$-LgLCA for supercompact implies $\RP^*$. 
\end{Thm}
\memo{A.Suzuki20250514132458.pdf p.18〜}
\prf \assertof{1}: Assume that $\calP$ is proper, and $\calP$-LgLCA for 
supercompact holds. Let $\calS=\seqof{S_\xi}{\xi<\eta}$ for some $\eta<\kappa_\refl$ \st\ 
each $S_\xi$, $\xi<\eta$ is a stationary subset of $[\lambda]^{\aleph_0}$. 
\Wolog, we may assume $\lambda\geq\kappa$, 

Let $\kappa:=\kappa_\refl$. 
$\poP\in\calP$ be arbitrary, and let $\poP\ast\utpoQ$, $\genH$, $j$, $M$ be 
\st\ $\forces{\poP}{\utpoQ\in\calP}$, $\genH$ is a $(\uniV,\poP\ast\utpoQ)$-generic set,
$\Elembed{j}{\uniV}{M}{\kappa}$, $j(\kappa)>\lambda$, $\poP$, $\poP\ast\utpoQ$, $\genH$,
$j\imageof{\lambda}\in M$, and $\cardof{RO(\poP\ast\utpoQ)}\leq j(\kappa)$. 

Then $j\imageof{S_\xi}=j(S_\xi)\cap [j\imageof{\lambda}]$ and 
$\uniV[\genH]\modelof{j\imageof{S}\mbox{ is stationary in }[j\imageof{\lambda}]^{\aleph_0}}$
for all $\xi<\eta$, the latter holds since $\poP\ast\utpoQ$ is proper. Thus
$M\modelof{j\imageof{S}\mbox{ is stationary in }[j\imageof{\lambda}]^{\aleph_0}}$ for all
$\xi<\eta$.

Since $j(\calS)=\seqof{j(S_\xi)}{\xi<\eta}$, we have
\begin{xitemize}
\item[] $M\modelof{\,
  \begin{array}[t]{@{}l}
    \mbox{there is }X\in[j(\lambda)]^{\LT j(\kappa)}
    \mbox{\st\ for all components }S\mbox{ of }j(\calS)\\
    S\cap[X]^{\aleph_0}\mbox{ is stationary in }[X]^{\aleph_0}
  }.\end{array}$
\end{xitemize}

By elementarity, it follows that 
\begin{xitemize}
\item[] $M\modelof{\,
  \begin{array}[t]{@{}l}
    \mbox{there is }X\in[\lambda]^{\LT \kappa}
    \mbox{\st\ for all components }S\mbox{ of }\calS\\
    S\cap[X]^{\aleph_0}\mbox{ is stationary in }[X]^{\aleph_0}
  }.\end{array}$
\end{xitemize}

\assertof{2}: Suppose that $\calP$ is a class consisting of stationary preserving \pos\ 
\st\ all $\sigma$-closed \pos\ are elements of $\calP$. 

Suppose that $\lambda$ and $\calS=\seqof{S_\xi}{\xi<\eta}$ are as above. Let
$\poP:=\Col(\aleph_1, \lambda^+)$ and $\utpoQ$ be a $\calP$-name with $\genH$, $j$, $M$ 
\st\ $\forces{\poP}{\utpoQ\in\calP}$, $\genH$ is a $(\uniV,\poP\ast\utpoQ)$-generic set,
$\Elembed{j}{\uniV}{M}{\kappa}$, $j(\kappa)>\lambda$, $\poP$, $\poP\ast\utpoQ$, $\genH$,
$j\imageof{\lambda}\in M$, and $\cardof{RO(\poP\ast\utpoQ)}\leq j(\kappa)$. 

Let $\genG$ be the $\poP$ part of $\genH$. Since $\poP$ is proper, we have
$\uniV[\genG]\modelof{S_\xi\xmbox{ is stationary in }[\lambda]^{\aleph_0}}$ where we have
$\uniV[\genG]]\models\cardof{\lambda}=\aleph_1$. Thus 
$\uniV[\genH]\modelof{j\imageof{S}\mbox{ is stationary in }[j\imageof{\lambda}]^{\aleph_0}}$, 
since $\utpoQ[\genG]$ is stationary preserving in $\uniV[\genG]$.

Now the final part of the proof of \assertof{1} can be repeated here to obtain the same 
conclusion as in \assertof{1}. 
\qedofThm\qedskip

\Thmof{p-square-0},\,\assertof{2} will be surpassed by \Thmof{p-square-3}: the assumptions 
of \Thmof{p-square-0},\,\assertof{2} imply \MApp($\sigma$-closed) by \Thmof{p-square-3}. 
Cox \cite{cox} proved that \MApp($\sigma$-closed) implies the strongest form of 
Diagonal Reflection Principle and our $\RC^*$ follows from from this. 

Actually, with almost the same proof we can also prove \Thmof{p-square-0} with $\RP^*$ replaced 
by the following Diagonal Reflection Principle which implies $\RP^*$:
\begin{xitemize}
\item[(\RP$^{**}$): ]  For any uncountable $\lambda\geq\kappa$, $\eta<\kappa_\refl$ and 
  $S_\xi\subseteq[\lambda]^{\aleph_0}$, $\xi<\lambda$ 
  \st\ $S_\xi$ is stationary in $[\lambda]^{\aleph_0}$ for all $\xi<\lambda$, there are 
  stationarily many 
  $X\in[\lambda]^{\LT\kappa_\refl}$ \st\ $\eta\subseteq X$, and  $S_\xi\cap[X]^{\aleph_0}$ is stationary in
  $[X]^{\aleph_0}$ for all $\xi\in X$.
\end{xitemize}

\begin{Cor}\Label{p-square-1} Assume that $\calP$ is a class of \pos\ which is either 
  proper, or stationary preserving containing all $\sigma$-closed \pos.  
  Then $\calP$-LgLCA for supercompact implies $\neg\square_\mu$ for all $\mu\geq\kappa_\refl$. Further, if
  $\kappa_\refl\leq\aleph_2$, then $\calP$-LgLCA for supercompact implies the total failure of 
  square principles.
\end{Cor}
\prf $\RP^*$ implies that all stationary
$S\subseteq E_{\omega_0}^{\mu^+}:=\setof{\alpha<\mu^+}{\cf(\alpha)=\omega_0}$ for all $\mu$ 
with $\mu^+\geq\kappa_\refl$ reflects. This implies $\neg\square_{\mu}$.
\qedofCor\qedskip

``$\mu>\kappa_\refl$'' in the Corollary above is optimal.
\begin{Lemma}\Label{p-square-2} Suppose that $\calP$ is transfinitely iterable class of 
  \pos\ \st\ all elements of $\calP$ are cardinality preserving. Then $\calP$-LgLCA for 
  supercompact is compatible with the statement that there are cofinally many $\mu<\kappa_\refl$ 
  with $\square_\mu$. 
\end{Lemma}
\prf Let $\kappa$ be a supercompact cardinal. By Easton support iteration of 
length $\kappa$ we can force that cofinally many $\mu<\kappa$ satisfy $\square_\mu$ while 
$\kappa$ is kept to be supercompact. If we force $\calP$-LgLCA by the standard construction 
starting form this model, $\square_\mu$-sequences of the cofinally many $\mu<\kappa$ 
remain $\square_\mu$-sequences in the generic extension since $\mu^+$ in the model at the 
start remains $\mu^+$ in the generic extension. \qedofLemma\qedskip

If $\calP$ is either the class of all proper \pos, or the class of all semi-proper \pos, the 
reason of the reflection properties available under $\calP$-LgLCA for supercompact is 
simply that $\calP$-LgLCA implies \PFApp\ and \MMpp\ respectively.
More generally, we have the following.

We consider the following axiom for a class $\calP$ of \pos:
\newcommand{\MAxx}{\mbox{\sf MA$^{\ast\ast}$}}

\begin{xitemize}
\item[($\MAxx(\calP)$): ] For any $\poP\in\calP$ and $\seqof{D_\alpha}{\alpha<\mu}$,
  $\seqof{\utS_\beta}{\beta<\nu}$ where $\mu$, $\nu<\kappa_\refl$, $D_\alpha\subseteq\poP$ 
  is dense subset of $\poP$ for all $\alpha<\mu$ and $\utS_\beta$ is a $\poP$-name of a 
  stationary subset of $[\lambda_\beta]^{\aleph_0}$ for some uncountable
  $\lambda_\beta<\kappa_\refl$, for all $\beta<\nu$, there is a filter $\genG\subseteq\poP$ 
  \st\ $\genG\cap D_\alpha\not=\emptyset$ for all $\alpha<\mu$ and $\utS_\beta[\genG]$ is a 
  stationary subset of $[\lambda_\beta]^{\aleph_0}$ for all $\beta<\nu$. 
\end{xitemize}
Note that, if $\kappa_\refl=\aleph_2$, $\MAxx(\calP)$ is simply equivalent to the 
usual $\MApp(\calP)$. 

\begin{Thm}{\rm(Theorem 5.7 in \cite{sfetal-II} restated under LgLCA for 
    extendible)}\Label{p-square-3} Suppose that an iterable class of \pos\ $\calP$ is 
  either  
  proper, or semi-proper and contains all $\sigma$-closed \pos. Then $\calP$-LgLCA for 
  supercompact implies $\MAxx(\calP)$.
\end{Thm}
\prf We show the proof for the case that $\calP$ is proper. The proof for the other case is 
similar using the fact that $\calP$-LgLCA 
implies $\kappa_\refl=\aleph_2$ for $\omega_1$ preserving $\poP$ containing 
all $\sigma$-closed \pos\ (see \cite{sfetal-II}).
\memo{see also talk-nagoya-2019-05-31.pdf p.3}

Thus, let us assume that $\calP$ is proper, and let $\kappa:=\kappa_\refl$. 

Let $\poP\in\calP$ be arbitrary,  and 
$\calD:=\seqof{D_\alpha}{\alpha<\mu}$,
$\calS:=\seqof{\utS_\beta}{\beta<\nu}$, $\seqof{\lambda_\beta}{\beta<\nu}$ be as in the definition of
$\MAxx(\calP)$. We show that there is a $\genG\subseteq\poP$ as in the definition of
$\MAxx(\calP)$. for these sequences.

\Wolog, we may assume that the underlying set of $\poP$ is a cardinal $\lambda>\kappa$, and 
$\utS_\beta$ are nice $\poP$-names (e.g.\ in the sense of Kunen \cite{kunen-2011}).

Let $\utpoQ$ be a $\poP$-name \st\ 
\begin{xitemize}
\xitem[x-square-0] $\forces{\poP}{\utpoQ\in\calP}$
\end{xitemize}
and that, for $(\uniV,\poP\ast\poQ)$-generic $\genH$, there 
are $j$, $M\subseteq\uniV[\genH]$ 
\st\ $\Elembed{j}{\uniV}{M}{\kappa}$, $j(\kappa)>\lambda$, 
$\poP$,
$\poP\ast\utpoQ$, $\genH$, $j\imageof{\lambda}\in M$, and
$\cardof{RO(\poP\ast\poQ)}\leq j(\kappa)$.

Let $\genH$, $j$, $M$ be as above and let $\genG$ be the $\poP$ part of $\genH$. 
Note that $j(\calD)=\seqof{j(D_\alpha)}{\alpha<\mu}$, and
$j(\calS)=\setof{j(\utS_\beta)}{\beta<\nu}$. 

By the closure properties of $M$, we have $j\imageof{\genG}\in M$, 
$j\imageof{\genG}\subseteq j(\poP)$. $j\imageof{\genG}$ is a subset of $j(\poP)$ with 
finite intersection property by elementarity of $j$. Let $\genG^*$ be the filter 
$\subseteq j(\poP)$ generated by $j\imageof{\genG}$. Then $\genG^*\in M$, and 
$\genG^*\cap j(D_\alpha)\supseteq \genG^*\cap j\imageof{D_\alpha}\not=\emptyset$ for all
$\alpha<\mu$.

We also have
$\uniV[\genH]\modelof{j(\utS_\beta)[\genG^*]=\utS_\beta[\genG]
  \mbox{ is a stationary subset of }[\lambda_\beta]^{\aleph_0}}$ 
for $\beta<\nu$ by the choice of $\calS$ and \xitemof{x-square-0}. It follows that 
\begin{xitemize}
\xitemx[] 
  $M\modelof{j(\utS_\beta)[\genG^*]=\utS_\beta[\genG]
    \mbox{ is a stationary subset of }[\lambda_\beta]^{\aleph_0}}$.
\end{xitemize}

Thus,
\begin{xitemize}
\xitemx[] 
  $M\modelof{\exists\ul{G}\,(\,
  \begin{array}[t]{@{}l}
    \ul{G}\mbox{ is a }j(\calD)\mbox{-generic filter on }j(\poP)\mbox{, and for each component }\utS\\
    \mbox{of }
    j(\calS)\mbox{ with the corresponding }\lambda\mbox{, }\utS[\ul{G}]\mbox{ is a 
      stationary subset}\\
    \mbox{of }[\lambda]^{\aleph_0}
    )}.
  \end{array}$
\end{xitemize}
By elementarity of $j$, it follows that 
\begin{xitemize}
\xitemx[] 
  $\uniV\modelof{\exists\ul{G}\,(\,
  \begin{array}[t]{@{}l}
    \ul{G}\mbox{ is a }\calD\mbox{-generic filter on }\poP\mbox{, and for each component }\utS\\
    \mbox{of }
    \calS\mbox{ with the corresponding }\lambda\mbox{, }\utS[\ul{G}]\mbox{ is a 
      stationary subset}\\
    \mbox{of }[\lambda]^{\aleph_0}
    )}
  \end{array}$
\end{xitemize}
as desired.
\qedofThm\qedskip

Fodor-type Reflection Principle (\FRP) is a reflection principle with reflection point
$\LT\aleph_2$. \FRP\ has many ``mathematical'' characterizations. The following is one of 
such characterizations in terms of non-metrizability of topological spaces.

\begin{Prop}\Label{p-square-4}{\rm (see \cite{more}, Theorem 2.8)}
  \FRP\ is equivalent to the following statement:\vspace{-0.9ex}
  \begin{xitemize}
  \xitem[x-square-1] For any locally countably compact topological space $X$ if $X$ is not 
    metrizable then there is a subspace $Y$ of $X$ of cardinality $\LT\aleph_2$ which is 
    not metrizable. 
  \end{xitemize}
\end{Prop}
\memo{see Lemma 1.3 in \cite{fjetal}}

In the following, we shall not go into the combinatorial definition of \FRP. Instead, we 
use the following basic facts about \FRP\ to sestablish \Thmof{p-square-6} showing the 
connection of \FRP\ to LgLCA.

\begin{Lemma}\Label{p-square-5} \wassertof{1} {\rm(\cite{fjetal}, 
    Theorem 2.5)} \RP$^*$ implies \FRP.\smallskip
  
  \wassert{2} {\rm(\cite{fjetal}, Proposition 1.5)}
  Non-reflecting stationary set $S\subseteq E_{\omega_0}^\lambda$ for 
  any uncountable $\lambda$ creates a counter-example to \xitemof{x-square-1}. \smallskip

  \wassert{3} {\rm(\cite{fjetal}, Theorem 3.4)} \FRP\ is preserved by 
  c.c.c.\ forcing. \qed
\end{Lemma}

\begin{Thm}\Label{p-square-6} \wassertof{1} Suppose that $\calP$ stationary preserving 
  containing all $\sigma$-closed \pos. then $\calP$-LgLCA for supercompact implies 
  \FRP.\smallskip

  \wassertof{2} If $\calP$ is c.c.c.\ 
  and transfinitely iterable class of \pos, then \FRP\ is independent over $\calP$-LgLCA for supercompact. 
\end{Thm}
\prf \assertof{1}: By \Thmof{p-square-0},\,\assertof{2} and 
\Lemmaof{p-square-5},\,\assertof{1}.\smallskip

\assertof{2}: Assume that $\calP$ is c.c.c. (remember that by the convention we introduced 
in \sectionof{intro}, this means that all $\poP\in\calP$ are c.c.c.).

To establish the consistency of $\FRP$ over $\calP$-LgLCA for supercompact, we start from a 
model with two supercompact cardinals. We use the first supercompact to force \FRP\ \st\ 
the second supercompact servives. Using this remaining supercompact cardinal $\kappa$ we 
force $\calP$-LgLCA for supercompact over the first generic extension. This can be done by 
a c.c.c.\ forcing and thus \FRP\ survives the generic extension by 
\Lemmaof{p-square-5},\,\assertof{3}.

To show the consistency of $\neg\FRP$ over $\calP$-LgLCA for supercompact, we now 
start from a model with a supercompact cardinal $\kappa$. We force $\square_\mu$ for some
$\mu<\kappa$ with a small forcing. Over this generic extension in which $\kappa$ remains 
supercompact, we force $\calP$-LgLCA for supercompact using $\kappa$ by a forcing 
in $\calP$.
$\square_\mu$-sequence survives the second generic extension since $\mu^+$ does not change 
by the second generic extension. 
Now $\square_\mu$ implies the existence of a non-reflecting stationary subset of
$E_{\omega_0}^{\mu^+}$. Thus by \Lemmaof{p-square-5},\,\assertof{2}, $\neg\FRP$ holds in 
the second generic extension.\qedofThm\qedskip

There are several other reflection of properties for which their preservation in the 
generic extensions is known only for special classes of \pos, and thus the corresponding 
LgLCAs are also rather special.

{\It Game Reflection Principle} (\GRP, which is called the global Game Reflection Principle 
and denoted by \GRPp\ in König \cite{koenig}) is the statement that the non-existence 
of a winning strategy of Player II of the game of sertain kind is reflected down to a 
subgame over an underlying set of cardinality $\LT\aleph_2$.  Since non-existence of a 
winning strategy is preserved by $\sigma$-closed \pos\ it is easy to prove 
that $\calP$-LgLCA for supercompact for the class $\calP$ of all $\sigma$-closed \pos\ implies 
\GRP\ by an argument similar to the proof of \Thmof{p-square-0}.

However, König proved in \cite{koenig} 
that \GRP\ is characterized by $\aleph_2$ being generic supercompact by $\sigma$-closed 
forcing which is apparently a weakening of $\calP$-LgLCA for $\calP$ as above. It is also proved in 
\cite{koenig} that \GRP\ implies Rado's Conjecture. Thus, we have
\begin{Thm}
  Let $\calP$ be the class of all $\sigma$-closed \pos. Then $\calP$-LgLCA for supercompact 
  implies \GRP. In particular, $\calP$-LgLCA for supercompact implies Rado's Conjecture. \qed
\end{Thm}

Non-freeness of a structure is preserved by c.c.c. forcing (see e.g. \cite{potential}, 
Theorem 2.1). Using this fact, an argument similar to the proof of \Thmof{p-square-0} shows 
the following: 
\begin{Thm} Suppose that $\calP$ is c.c.c. Then $\calP$-LgLCA for supercompact implies the 
  following reflection theorem:\vspace{-0.9ex}
  \begin{xitemize}
  \xitemx[] For any non-free algebra $A$ (in a universal algebraic class of strucrures), 
    there is a non-free subalgebra $B$ of $A$ of cardinality $<\continuum$.\qed
  \end{xitemize}
\end{Thm}

Let us call $\calP$ (and elements of $\calP$) Cohen, if
\begin{xitemize}
\xitemx[] 
  $\calP=\setof{\poP}{\poP \mbox{ is forcing equivalent to }\Fn(\kappa,2)\mbox{ for some infinite }\kappa}$ 
\end{xitemize}
where $\Fn(\kappa,2)$ is the usual \po\ adding $\kappa$ many Cohen reals by finite conditions.

Dow, Tall and  Weiss \cite{dtw} proved that Cohen \pos\ preserve the 
non-metrizability of topological spaces. In a supercompact elementary embedding context 
with $\Elembed{j}{\uniV}{M}{\kappa}$, $j\imageof{\calO}$ generates the subspace topology of
$j\imageof{X}$, in $(j(X),j(\calO))$ for a topological space $X=(X,\calO)$ with character
$\LT\kappa$.

Thus an argument  
similar to the proof of \Thmof{p-square-0} establishes: 
\begin{Thm} Suppose that $\calP$ is Cohen. Then $\calP$-LgLCA for supercompact implies the 
  following reflection theorem:\vspace{-0.9ex}
  \begin{xitemize}
  \xitemx[] If $X$ is a non-metrizable topological space with character $\LT\continuum$ 
    then there is a subspace $Y$ of $X$ of cardinality $\LT\continuum$ which is also 
    non-metrizable. \qed
  \end{xitemize}
\end{Thm}

Singular Cardinal Hypothesis (\SCH) can be seen as a reflection principle --- see e.g.\ Fuchino, and 
Rinot \cite{rinot}, \href{https://fuchino.ddo.jp/papers/scepin10-x.pdf\pdfpage{13}}{Theorem 2.4} 
in which it is shown that Shelah's Strong Hypothesis (\SSH), a generalization of \SCH, is 
characterized as reflection of a topological property. 

In \cite{rinot}, it is shown that \FRP\ implies \SSH. Thus \Thmof{p-square-6} above implies 
that some instances of LgLCA for supercompact imply \SSH. By virtue of the definition of 
LgLCA in the present paper, we can reformulate 
\href{https://fuchino.ddo.jp/papers/SDLS-II-x.pdf\pdfpage{13}}{Proposition 2.8} in 
\cite{sfetal-II} as:
\begin{Thm}\Label{p-square-7}{\rm 
    (\href{https://fuchino.ddo.jp/papers/SDLS-II-x.pdf\pdfpage{13}}{Propsition 2.8} in 
    \cite{sfetal-II})} 
  For any iterable class $\calP$ of \pos, $\calP$-LgLCA for supercompact implies \SCH.\qed
\end{Thm}

\section{Separation of some axioms under Continuum Coding Axiom}\Label{CCA}
In this section, we show the separation of some of the axioms and principles we considered 
in the previous sections. 

We start by showing that \MMpp, as well as the property in \xitemof{x-genabs-Laver-a} for 
all stationary preserving \pos\ $\calP$, does not imply LgLCAs for supercompact. 
\memo{HKP20250514132350.pdf p.6〜}
The main tools of the separation are the Ground Axiom (\cite{reitz}) and some other related 
axioms.

Recall that a ground is an inner model $M$ of the universe $\uniV$ from which we can return to 
the universe by a set generic extension (for a \po\ in $M$). 
{\It Ground Axiom} ({\It\GA}) is the axiom asserting that there is no non-trivial ground.
If a ground $M$ satisfies the Ground Axiom then we call $M$ the bedrock. {\It Bedrock Axiom} 
({\It\BA}) is the assertion that the bedrock exists.

{\It Continuum Coding Axiom} ({\It\CCA}) is the axiom saying that each set is coded class 
many times by the pattern $\setof{\mu^+}{\alpha<\mu^+<\beta,2^{\mu^+}=\mu^{++}}$. It is 
easy to see that \CCA\ implies global choice and \GA. On the other hand, \CCA\ negates any 
form of $\calP$-LgLCAs if $\calP$ is non-trivial (see below).  

\imemox{Definition of \CCA\ \BA\ \GA}

\begin{Thm}\Label{p-CCA-0} \wassertof{1}
  \MMpp\ does not imply $\calP$-LgLCA for 
  supercompact for any non-trivial class $\calP$ of \pos.\smallskip

  \wassert{2} The conclusion of Viale's Absoluteness Theorem, namely the assertion
  \xitemof{x-genabs-Laver-a} for all stationary preserving \pos\ $\calP$
  does not imply $\calP$-LgLCA for supercompact for any non-trivial class $\calP$ of \pos.
\end{Thm}
\prf \assertof{1}: In \cite{cox2}, it is shown that \MMpp\ is preserved by $\LT\omega_2$-directed closed 
forcing. 
Thus starting from a model of \MMpp\ and class many supercompact cardinals. 
We can class force with $\LT\omega_2$-directed closed class forcing to obtain a model of 
\MMpp\ $+$ \CCA. If we start from a model with class many supercompact cardinals. We first 
force \MMpp\ using the least supercompact cardinal without destroying other class many 
supercompact cardinals. We then force Laver indestructivity for the rest of the 
supercompact cardinals. Then this generic extension models \CCA\ (see the lemma below). 

On the other hand, under \CCA, $\calP$-LgLCA for 
supercompact cannot  not hold for non-trivial $\calP$: Suppose that $\poP\in\calP$ is a 
nontrivial \po. For $\poP$-name $\utpoQ$ and $(\uniV,\poP*\utpoQ)$-generic $\genH$, if 
there are $j$, $M\subseteq\uniV[\genH]$ with $\Elembed{j}{\uniV}{M}{\kappa}$ and
$\genH\in M$.  Then $\genH$ witnesses $M\models\neg\CCA$. This is a contradiction to
$\uniV\models\CCA$ and the elementarity of $j$. \smallskip

\assertof{2}: By Viale's Absoluteness Theorem, the model mentioned in the proof of 
\assertof{1} satisfies \xitemof{x-genabs-Laver-a} while it does not satisfy $\calP$-LgLCA as we saw 
above. 
\qedofThm\qedskip

\begin{Lemma}\Label{p-CCA-0-a}{\rm (G.\,Goldberg)} If there are class many 
  (directed closed) indestructible supercompact cardinals then \CCA\ holds. 
\end{Lemma}
\prf Let $\gamma\in\On$ be a limit ordinal and $S\subseteq\gamma$ be arbitrary. For $\alpha\in\On$, let
$S_{\alpha,+\gamma}:=\setof{\alpha+\xi}{\xi\in S}$. For any $\beta\in\On$ we show that there is
$\alpha^*>\beta$ \st\ $S_{\alpha^*,+\gamma}=S^*_{\alpha^*,+\gamma}$
where
$S^*_{\alpha^*,+\gamma}:=\setof{\xi}{\alpha^*\leq\xi<\alpha^*+\gamma,\,2^{\aleph_{\xi+1}=\aleph_{\xi+2}}}$.

Let $\kappa>\alpha$, $\beta$, $\gamma$ be indestructible supercompact cardinal. Then we easily find 
a $\LT\kappa$-direcrted closed $\poP$ \st, for $(\uniV,\poP)$-generic $\genG$, we have
$\uniV[\genG]\models S_{\kappa,+\gamma}=S^*_{\kappa,+\gamma}$.
Since $\kappa$ is still supercompact in $\uniV[\genG]$ and hence $\Sigma_2$-correct, it 
follows that there is $\beta<\delta<\kappa$, \st\
$\uniV[\genG]\models S_{\delta,+\gamma}=S^*_{\delta,+\gamma}$. Since $\poP$ does not change 
the continuum function below $\kappa$ it follows that
$\uniV\models S_{\delta,+\gamma}=S^*_{\delta,+\gamma}$. 
\qedofLemma\qedskip

Some LgLCAs are separated by their consistency strength. The following theorem has actually 
a generalization for tight $\calP$-generic hyperhuge cardinals, and the 
generalization also implies Usuba's Main Theorem in \cite{usuba} asserting that a 
hyper-huge cardinal implies \BA. Usuba proved that this theorem in \cite{usuba} already 
holds under the existence of an extendible cardinal \cite{usuba2}. At the moment it is open if 
corresponding improvement of the theorem below down to LgLCA for extendible is possible. 

\begin{Thm}{\rm (Fuchino and Usuba \cite{recurrence})}\Label{p-CCA-0-0}
  \wassertof{1} For any class $\calP$ 
  of \pos, if $\calP$-LgLCA for hyperhuge is consistent, then \tfae:\smallskip

  \wassert{a} \ZFC\ $+$ $\calP$-LgLCA for hyperhuge;

  \wassert{b} \ZFC\ $+$ ``there is a hyperhuge cardinal'';

  \wassert{c} \ZFC\ $+$ \BA\ $+$ $\kappa_\refl$ is hyperhuge in the bedrock. \smallskip

  \wassert{2} For any class $\calP$ 
  of \pos, if super-$C^{(\infty)}$-$\calP$-LgLCA for hyperhuge is consistent, then \tfae:\smallskip
 
  \wassert{d} \ZFC\ $+$ super $C^{(\infty)}$-$\calP$-LgLCA for hyperhuge;

  \wassert{e} \ZFC\ $+$ \BA\ $+$ $\kappa_\refl$ is super-$C^{(\infty)}$ hyperhuge in the bedrock. \qed
\end{Thm}

\begin{Thm}{\rm(Tsaprounis \cite{tsaprounis2}, Proposition 3.2)}\Label{tsaprounis-0}
 If $\kappa$ is an 
  ultrahuge cardinal then there exists a normal measure one many superhuge $\alpha<\kappa$.\qed
\end{Thm}

\begin{Thm}\Label{p-CCA-1} \wassertof{1} Suppose that $\calP$ is a transfinitely iterable
  $\Sigma_2$-definable class of \pos\ and 
  $\ZFC$ $+$ $\calP$-LgLCA for superhuge is 
  consistent. Then $\ZFC$ $+$ $\calP$-LgLCA for superhuge does not prove $\calP$-LgLCA for 
  hyperhuge.\smallskip

  \wassert{2} Suppose that $\calP$ is a transfinitely iterable class of \pos\ (its 
  definition may be more complex than $\Sigma_2$). 
  $\ZFC$ $+$ super-$C^{(\infty)}$-$\calP$-LgLCA for extendible is 
  consistent. Then $\ZFC$ $+$ super-$C^{(\infty)}$-$\calP$-LgLCA for extendible does not prove $\calP$-LgLCA for 
  hyperhuge. 
\end{Thm}
\prf \assertof{1}: Suppose otherwise. Then, since a hyperhuge cardinal is ultrahuge, 
\Thmof{p-CCA-0-0},\,\assertof{1} and \Thmof{tsaprounis-0} imply
\begin{xitemize}
\xitem[] \ZFC\ $+$ $\calP$-LgLCA for superhuge\ \ $\vdash$\ \ 
  $consis(\ZFC\ +\ \mbox{``\,$\exists$ a superhuge cardinal''})$. 
\end{xitemize}
Since $\calP$ is transfinitely iterable and $\Sigma_2$-definable, if follows (much like in the proof of 
\Thmof{p-Lg-ext-1}, \assertof{1})
that \memo{Check that the condition ``$\calP$ is $\Sigma_2$'' can be eleminated by 
  replacing superhuge by super-$C^{(n)}$ superhuge for sufficiently large $n$. Show that 
  Tsaprounis Theorem can be extended for super-$C^{(n)}$ superhuge.}
\begin{xitemize}
\xitem[] \ZFC\ $+$ $\calP$-LgLCA for superhuge $\vdash$
  $consis(\ZFC\ +\ \mbox{$\calP$-LgLCA for superhuge})$. 
\end{xitemize}
By the Second Incompleteness Theorem, it follows that ``\ZFC\ $+$ $\calP$-LgLCA for 
superhuge'' is inconsistent. This is a contradiction to our assumption.
\smallskip

\assertof{2}: Similarly to \assertof{1} using \Lemmaof{p-ext-0}, \Propof{p-ext-1} as well 
as an analogue of \Thmof{p-Lg-ext-1}, \assertof{2} in place of \Thmof{tsaprounis-0} and 
\Thmof{p-Lg-ext-1}, \assertof{1}.
\qedofThm

\begin{Lemma}\Label{p-CCA-1-0}
  Suppose that $(\calP,\emptyset)_{\Pi_2}$-\RcA\ holds for a class of \pos\ $\calP$
  \st\ either $\calP$ collapses arbitrary large cardinal making it an ordinal of small cardinality, or $\calP$ adds 
  more reals than any given cardinal.

  If there is an inaccessible cardinal. Then there are class 
  many inaccessible cardinals.
\end{Lemma}
\prf Assume towards a contradiction, that there are at least one but only set many 
inaccessible cardinals. By assumption on $\calP$, there is a $\poP\in\calP$ \st\
$\forces{\poP}{\mbox{there is no inaccessible cardinal}}$. Since the statement is $\Pi_2$ 
without any parameter, $(\calP,\emptyset)_{\Pi_2}$-\RcA\ implies that there is a ground $M$ 
without any inaccessible cardinal. But this is a contradiction since no inaccessible 
cardinal can resurrect by a set forcing. \qedofLemma

\begin{Prop}\Label{p-CCA-2} Suppose that $\calP$ is a transfinitely 
  iterable $\Sigma_2$-definable 
  class of \pos\ \st\ either $\calP$ collapses arbitrary large cardinal making it an ordinal of 
  small cardinality, or $\calP$ adds more reals than any given 
  cardinal. 
  
  Then 
  $\calP$-LgLCA for supercompact does not imply $\calP$-LgLCA for 
  extendible. 
\end{Prop}
\prf Start from a universe with a supercomact $\kappa$ and inaccessible $\lambda$ above it. We 
may assume that $\lambda$ is the unique inaccessible cardinal above $\kappa$ if there is 
another inaccessible cardinal above $\lambda$ then we can take the least $\lambda_0$ among such 
and consider $V_{\lambda_0}$ to be our universe. By assumption on $\calP$. We can construct 
a generic extension $\uniV[\genG]$ \st\
\begin{xitemize}
\xitemx[] 
  $\uniV[\genG]\models \kappa=\kappa_\refl\ \land\ \calP\xmbox{-LgLCA for supercompact\\
  \phantom{$\uniV[\genG]\models$}}
  \land\ \mbox{``}\lambda\xmbox{ is the unique inaccessible 
    cardinal above }\kappa\mbox{''}$.  
\end{xitemize}
On the other hand, by \Thmof{p-Lg-RcA-0} and \Lemmaof{p-CCA-1-0}, we have
\begin{xitemize}
\xitemx[] $\uniV[\genG]\models \neg\,(\calP\mbox{-LgLCA for extendible})$.\qedofProp
\end{xitemize}

\memo{Open: LgLCA for super-almost huge implies LgLCA for extendible? Does LgLCA for 
  extendible imply GA?}
\fi}

\section{Laver-generic Large Cardinal Axioms for all \pos}\Label{LgLCAA}
In this section, we want to examine the principles which we will call the 
{\It Laver-generic Large Cardinal Axioms for all \pos}\ and abreviate it as {\It LgLCAAs}:

For a notion  LC of large cardinal, {\It LgLCAA for LC\/} is the following assertion: 
\begin{xitemize}
\xitem[x-LgCAA-0] For any $\lambda>\continuum$ and for any \po\ $\poP$, there is 
  a $\poP$-name $\utpoQ$ of a \po\ with the property that, for 
  a $(\uniV,\poP\ast\utpoQ)$-generic $\genH$,  
  there are $j$, $M\subseteq\uniV[\genH]$ \st\ \assertof{a}\ \,$\Elembed{j}{\uniV}{M}{\continuum}$,
  \assertof{b}\ \,$j(\continuum)>\lambda$, \assertof{c}\ \,$\poP$, $\poP\ast\utpoQ$,
  $\genH\in M$, \assertof{d}\ \,$\cardof{RO(\poP\ast\utpoQ)}\leq j(\continuum)$, and 
  \assertof{e}\ \,$M$ satisfies the closure property corresponding to LC.
\end{xitemize}

This version of Laver-genericity has been already discussed in \cite{recurrence} and 
\cite{janos}. Below, we examine it in connection with extendibility.

\begin{Lemma}\Label{p-LgCAA-0} LgLCAA for any notion of large cardinal implies \CH.
\end{Lemma}
\prf Otherwise $j$ as in \xitemof{x-LgCAA-0} sends $\omega_1$ to itself. However $\poP$ can 
collapse $\omega_1$ and in such a case we would have
$M\modelof{\omega_1\mbox{ is countable}}$ by \xitemof{x-LgCAA-0},\,\assertof{c} and 
elementarity. This is a contradiction. 
\qedofLemma
\smallskip

Note that by the lemma above, the LgLCAA is different from ``the LgLCA for all \pos'' in the sense 
of previous sections in that  ``the LgLCA'' refers to $\kappa_\refl$ being the critical point of 
generic elementary embeddings while the present LgLCAA refers to
$\continuum=\aleph_1<\kappa_\refl$ being the critical point. Note also that by 
\Thmof{p-Lg-RcA-1}, ``the LgLCA for all \pos'' in the sense of previous sections is 
inconsitent.

For a notion  LC of large cardinal, we define {\It super-$C^{(\infty)}$-LgLCAA for LC\/} to 
be the following assertion: 
\begin{xitemize}
\xitem[x-LgCAA-1] For any $n\in\natnums$, $\lambda>\continuum$ with
  $\lambda\in C^{(\infty)}$, and for any \po\ $\poP$, there is  
  a $\poP$-name $\utpoQ$ of a \po\ with the property that, for 
  a $(\uniV,\poP\ast\utpoQ)$-generic $\genH$,  
  there are $j$, $M\subseteq\uniV[\genH]$ \st\ \assertof{a}\ \,$\Elembed{j}{\uniV}{M}{\continuum}$,
  \assertof{b}\ \,$j(\continuum)>\lambda$, \assertof{c}\ \,$\poP$, $\poP\ast\utpoQ$,
  $\genH\in M$, \assertof{d}\ \,$\cardof{RO(\poP\ast\utpoQ)}\leq j(\continuum)$,  
  \assertof{e}\ \,$M$ satisfies the closure property corresponding to LC, and
  $j(\lambda)\in C^{(n)}$.
\end{xitemize}

The following theorem can be proved similarly to \Thmof{p-Lg-ext-1} by constructing 
$\poP_\kappa$ as the direct limit of a FS iteration of length $\kappa$. 

\begin{Thm}\Label{p-LgCAA-1} \wassertof{1} If $\kappa$ is extendible, there is a \po\ $\poP_\kappa$ \st\ 
  $\forces{\poP_\kappa}{\kappa=\continuum\xmbox{ and LgLCAA for 
      extendible holds}}$.\smallskip

  \wassertof{2} Suppose
  $V_\mu\modelof{\kappa\xmbox{ is strongly super-}C^{(\infty)}\xmbox{ extendible}}$ for 
  an inaccessible cardinal $\mu$ and $\kappa<\mu$. Then there is
  $\poP_\kappa\in V_\mu$ \st\ 
  $V_\mu\models\forces{\poP_\kappa}{\kappa=\continuum
    \xmbox{ and super-}C^{(\infty)}\xmbox{-LgLCAA for extendible 
      holds}}$.\ifextended\else\qed\fi 
\end{Thm}
{\ifprivate\privatecolor\prf \assertof{1}: Let $f$ be a Laver function for 
extendible $\kappa$ ($f$ exists by \Lemmaof{p-Lg-ext-0},\,\assertof{1}).

Let $\seqof{\poP_\alpha,\utpoQ_\beta}{\alpha\leq\kappa, \beta<\kappa}$ be an FS-iteration 
of \pos\ \st\
\begin{itemize}
\item[] $\utpoQ_\beta:={}${$\left\{\,
  \begin{array}{@{}ll}
    f(\beta), &\mbox{if }f(\beta)\mbox{ is 
      a }\poP_\beta\mbox{-name of a \po};\\[\jot] 
    \poP_\beta\mbox{-name of the trivial forcing}, &\mbox{otherwise}.
  \end{array}\right.$}
\end{itemize}\smallskip

We show that $\forces{\poP_\kappa}{\kappa=\continuum\mbox{ and }
\kappa\mbox{ is tightly }\calP\mbox{-Lg extendibile}}$. for the class $\calP$ of all \pos. 

$\forces{\poP_\kappa}{\kappa=\continuum}$ is clear since there are cofinially many 
$\beta<\kappa$ \st\ $f(\beta)$ is $\poP_\beta$ name of a \po\ adding a real.

Let $\genG_\kappa$ be a $(\uniV,\poP_\kappa)$-generic filter. 
  In $\uniV[\genG_\kappa]$, suppose that $\poP$ is a \po\ and let $\utpoP$ be 
  a $\poP_\kappa$-name for $\poP$. 
  
Suppose that $\lambda>\kappa$.

Let $\Elembed{j}{\uniV}{M}{\kappa}$ be \st\ \quad\ixitem[x-1**] 
$j(\kappa)>\lambda$,\quad
\ixitemc[x-2**] $V_{j(\lambda)}\in M$, and\quad\\
\ixitem[x-3**] $j(f)(\kappa)=\utpoP$. 
The last condition is possible since $f$ is a Laver function for the extendible $\kappa$.

In $M$, there is a $\poP_\kappa\ast\utpoP$-name $\utpoQ$ \st\
\begin{xitemize}
\xitemx[] 
  $M\models\forces{\poP_\kappa\ast\utpoP}{{}
  \begin{array}[t]{@{}l}
    \utpoQ\xmbox{ is
      the direct limit of FS-iteration of small}\\
    \mbox{\pos\ of length }j(\kappa)\mbox{ and }\poP_\kappa\ast\utpoP\ast\utpoQ\sim j(\poP_\kappa)}.
  \end{array}$\\
  
\end{xitemize}
Then the same 
statement holds in $\uniV$ with the same $\utpoQ$.\smallskip 

We have $j(\poP_\kappa)/\genG_\kappa\ \sim\ \poP\ast\utpoQ$ where 
we identify $\utpoQ$ with a corresponding $\poP$-name.

Let $\genH$ be $(\uniV,j(\poP_\kappa))$-generic filter with $\genG_\kappa\subseteq\genH$. 
The lifting $\Elembed{\tilde{j}}{\uniV[\genG_\kappa]}{M[\genH]}{\kappa}$;
$\utilde{a}[\genG_\kappa]\mapsto j(\utilde{a})[\genH]$ witnesses that $\kappa$
  is tightly $\calP$-Laver generic extendible in $\uniV[\genG_\kappa]$ for the 
  class $\calP$ of all \pos. 
  $\cardof{RO(j(\poP_\kappa/\genG_\kappa))}\leq j(\kappa)$ follows from
  $\LT j(\kappa)$-c.c.\ of $j(\poP_\kappa)/\genG_\kappa$,
  $M\models j(\kappa)^{\LT j(\kappa)}=j(\kappa)$, and \xitemof{x-2**}. So it is enough to  
  show: 

\begin{Claim}\Label{cl-Lg-ext-0*}\ifprivate\privatecolor\fi
  ${V_\alpha}^{\uniV[\genH]}\in M[\genH]$ for all $\alpha\leq j(\lambda)$.
\end{Claim}
\prfofClaim By induction on $\alpha\leq j(\lambda)$. The successor step 
  from $\alpha<j(\lambda)$ to $\alpha+1$ 
  can be 
  proved by showing that $\poP_\kappa$-names of subsets of ${V_\alpha}^{\uniV[\genH]}$ can be 
chosen as elements of $M$. This is the case because of 
\xitemof{x-2**} and $\LT\kappa$-c.c.\ of $\poP_\kappa$.\qedofClaim\qedskip

\assertof{2}: In the following we work mainly in $V_\mu$ but our view-point is often set 
outside of $V_\mu$. Let $\calP$ be the class of all \pos, and  
let $f$ be a super-$C^{(\infty)}$-extendible-Laver function 
for $\kappa$ (which exists by \Lemmaof{p-Lg-ext-0},\,\assertof{2}). 

Let $\seqof{\poP_\alpha,\utpoQ_\beta}{\alpha\leq\kappa, \beta<\kappa}$ be a FS-iteration 
  of elements of \pos\ \st\
  \begin{itemize}
  \item[] $\utpoQ_\beta:={}${$\left\{\,
    \begin{array}{@{}ll}
      f(\beta), &\mbox{if }f(\beta)\mbox{ is 
        a }\poP_\beta\mbox{-name of a \po};\\[\jot] 
      \poP_\beta\mbox{-name of the trivial forcing}, &\mbox{otherwise}.
    \end{array}\right.$}
  \end{itemize}\smallskip

We show that
$V_\mu\models\forces{\poP_\kappa}{\kappa=\continuum\ \land\ \kappa\mbox{ is tightly }\calP\mbox{-Lg extendibile}}$. 

$V_\mu\models\forces{\poP_\kappa}{\kappa=\continuum}$ is clear since there are cofinally many 
(actually stationary many) $\beta<\kappa$ \st\ $f(\beta)$ is $\poP_\beta$-name of \po\ 
adding a real. 

Let $\genG_\kappa$ be a $(V_\mu,\poP_\kappa)$-generic filter. 
In $V_\mu[\genG_\kappa]$, suppose that $\poP$ is a \po\ and let $\utpoP$ be 
a $\poP_\kappa$-name for $\poP$. 
Suppose that $n\in\omega$ and $\mu>\lambda>\kappa$.

Then, there is 
$\mu>\lambda^*>\lambda$, $\lambda^*\in (C^{(n^*)})^{V_\mu}$ for sufficiently large $n^*\in\omega$ 
with $\Elembed{j}{V_\mu}{M}{\kappa}$ \st\ 
\ixitem[x-1*++] 
$j(\kappa)>\lambda^*$, 
  \,\ixitemc[x-1*+-0+] $V_{j(\lambda^*)}\prec_{\Sigma_{n^*}}V_\mu$, 
\,\ixitemr[x-2*+*] $V_{j(\lambda^*)}\in M$, and\ \,
\ixitem[x-3*+*] $j(f)(\kappa)=\utpoP$. 
The last condition is possible since $f$ is a Laver function. 


In $M$, there is a $\poP_\kappa\ast\utpoP$-name $\utpoQ$ \st\
\begin{xitemize}
\xitemx[] 
  $M\models\forces{\poP_\kappa\ast\utpoP}{{}
  \begin{array}[t]{@{}ll}
    \utpoQ\in\calP\mbox{ and }\utpoQ\mbox{ is
      the direct limit of FS-iteration of length }\\
    j(\kappa)\mbox{ of \pos\ of size}<j(\kappa)
    \mbox{, and }\poP_\kappa\ast\utpoP\ast\utpoQ\sim j(\poP_\kappa)}.
  \end{array}$ 
\end{xitemize}
By 
\xitemof{x-1*+-0+} and 
\xitemof{x-2*+*}, 
the same situation holds in $V_\mu$.\smallskip 
\memo{see Scan\_2024-11-23--10.44 extendible - annotated p.49}

We have $V_\mu[\genG_\kappa]\models j(\poP_\kappa)/\genG_\kappa\ \sim\ \poP\ast\utpoQ$ where 
we identify $\utpoQ$ with a corresponding $\poP$-name.

Let $\genH$ be $(V_\mu,j(\poP_\kappa))$-generic filter with $\genG_\kappa\subseteq\genH$. 
Then the lifting $\Elembed{\tilde{j}}{V_\mu[\genG_\kappa]}{M[\genH]}{\kappa}$;
$\utilde{a}[\genG_\kappa]\mapsto j(\utilde{a})[\genH]$ for each $\poP$-mame $\uta$ witnesses that $\kappa$
is tightly super-$C^{(n)}$-$\calP$-Laver generic extendible in $V_\mu[\genG_\kappa]$:  
${V_{j(\lambda)}}^{\uniV[\genH]}\in M[\genH]$ holds by an argument practically identical 
with that of \Claimof{cl-Lg-ext-0}.
\qedofThm\fi}

the LgLCAA for extendible and the super-$C^{(\infty)}$-LgLCAA for extendible satisfy 
corresponding to Theorems \ref{T-Lg-RA-0}, \ref{p-Lg-RcA-0}, \ref{p-Lg-RcA-5}, and 
\ref{p-genabs-Laver-1}. The proof of the following theorems can be obtained practically by 
changing the phrase like ``let $\kappa:=\kappa_\refl$'' by let ``$\kappa:=\continuum$''.

\begin{Thm} {\rm(A variation of \Thmof{T-Lg-RA-0})}
  Asshume that the LgLCAA for extendible holds. 
  Then for the class $\calP$ of all \pos, 
  $\BfRA^\calP_{\aleph_1}$ holds.  \qed
\end{Thm}
\begin{Thm}{\rm (A variation of \Thmof{p-Lg-RcA-0})}  Assume the LgLCAA for extendible.
  Then for the class $\calP$ of all \pos, 
  $(\calP,\calH(\aleph_1)_{\Gamma}$-\RcAp\ holds where $\Gamma$ is the set of all 
  formulas which are conjunctions of a $\Sigma_2$-formula and a $\Pi_2$-formula.\qed
\end{Thm}
\begin{Thm}{\rm(A variation of \Thmof{p-Lg-RcA-5})} 
  Suppose that the 
  super-$C^{(\infty)}$-LgLCAA for extendible holds. Then for the class $\calP$ of all \pos, 
  $\MP(\calP,\calH(\aleph_1))$ holds. \qed
\end{Thm}
\begin{Thm}{\rm(A variation of \Thmof{p-genabs-Laver-1})} Assume 
  that LgLCAA for extendible holds.   
  Then, for any \po\ $\poP$, 
  $\calH(\aleph_1)^\uniV\prec_{\Sigma_2}\calH(\aleph_1)^{\uniV[\genG]}$\quad
  holds.  \qed
\end{Thm}


\end{document}